\documentclass[matann2]{svjour}

\usepackage{amsmath}
\usepackage{amsfonts}
\usepackage{latexsym}
\usepackage{mathrsfs}
\usepackage{times}
\usepackage{array}
\usepackage{amscd}
\usepackage{a4}
\usepackage[mathcal]{euscript}
\usepackage{amssymb}
\usepackage{pstricks}
\usepackage{amsrefs}

\input{boxedeps} %% windows
\SetEPSFDirectory{} %% windows
\HideDisplacementBoxes
\SetRokickiEPSFSpecial  %% dvips by Tom Rokicki
\DeclareMathAlphabet{\ams}{U}{msb}{m}{n}
\DeclareMathAlphabet{\goth}{U}{euf}{m}{n}

\def\po{\text{PO}}
\def\sl{\text{SL}}

\def\gl{\text{GL}}

\def\vol{\text{vol}\,}
\def\covol{\text{covol}\,}
\def\im{\text{im}\,}
\def\ker{\text{ker}\,}
\def\aut{\text{Aut}}
\def\isom{\text{Isom}\,}
\def\endo{\text{End}}

\def\ov{\overline}

\def\II{\mathscr I}

\def\AA{\mathcal A}
\def\CC{\mathcal C}

\def\JJ{\mathcal J}

\def\TT{\mathcal T}

\def\SS{\goth{S}}

\def\aa{\alpha}
\def\ww{\omega}

\def\ss{\sigma}
\def\vphi{\varphi}
\def\wvphi{\widehat{\varphi}}

\def\ve{\varepsilon}
\def\Om{\Omega}

\def\Z{\ams{Z}}\def\E{\ams{E}}
\def\H{\ams{H}}\def\R{\ams{R}}
\def\Q{\ams{Q}}
\def\F{\ams{F}}

\def\v{\mathbf{v}}

\def\u{\mathbf{u}}
\def\x{\mathbf{x}}

\def\0{\mathbf{0}}

\def\quo{/\kern -.45em\sim}

\newpsobject{showgrid}{psgrid}
            {gridcolor=red,subgridcolor=red}

\def\Langle{\langle\kern -2pt\langle}
\def\Rangle{\rangle\kern -1.9pt\rangle}
\title{Weyl groups, lattices and geometric manifolds}
\author{Brent Everitt and
Robert B. Howlett
\thanks{Some of the results of this paper were 
obtained while the first author was visiting
the Department of Mathematics, University
of Sydney and
the Institute for Geometry and its Applications, %Department of Mathematics, 
University of Adelaide.
He is grateful for their hospitality and financial
assistance. The research done at the University of Sydney was supported by
a grant from the Australian Research Council.}
}

\institute{{\sc Brent Everitt:}
Department of Mathematics, University of York, York
YO10 5DD, United Kingdom. \email{bje1@york.ac.uk}. %(Brent Everitt), 
\hspace{1em}{\sc Robert Howlett:} Department of Mathematics, University of
Sydney, Sydney NSW 2006, Australia. 
\email{R.Howlett@maths.usyd.edu.au}.
}

\titlerunning{}
\authorrunning{Brent Everitt and Robert B. Howlett}

\begin{document}

\maketitle

%%%%%%%%%%%%%%%%%%%%%%%%%%%%%%%%%%%%%%%%%%%%%%%%%%%%%%%%%%%%%%%%%%%%%%%%%

\begin{abstract}
By studying the action of the Weyl group of a simple Lie algebra
on its root lattice, we construct torsion free subgroups of small and
explicitly determined index in a large infinite class of Coxeter groups. 
One spin-off is the construction of hyperbolic manifolds of very small
volume in up to $8$ dimensions.
\end{abstract}

%%%%%%%%%%%%%%%%%%%%%%%%%%%%%%%%%%%%%%%%%%%%%%%%%%%%%%%%%%%%%%%%%%%%%%%%%

\section*{Introduction}

The starting point of this paper is the celebrated theorem that every irreducible
representation of the Weyl group of a simple Lie algebra
can be realized over $\Q$. The finishing point
is the construction of hyperbolic manifolds of very small volume that provide
partial solutions to the Siegel problem in low dimensions.

The Weyl groups are examples of 
Coxeter groups, and it is these, together with the study
of their torsion free subgroups, that forms the bridge between the two 
sentences of the previous paragraph. % In particular the 
% classical theorem of Selberg \cite{Selberg60} implies that any Coxeter group
% has a torsion free subgroup of finite index. 
%
The realizability over $\Q$ of a linear action $W\rightarrow\gl(V)$
of a Weyl group $W$ is equivalent to the existence of a $W$-invariant
lattice $L$ in the ambient space $V$ of the representation.
In this paper we
use the finite semi-direct product $L\rtimes W$ induced by this action to
construct torsion free subgroups 
inside a large class of Coxeter groups to which the Weyl group $W$
is closely associated. Moreover, the indices of these subgroups are
determined explicitly in terms of well known data associated to the
Weyl group. As far as we are aware, this is the first such construction of
this type.

A feature in recent years 
of the study of torsion free subgroups of Coxeter groups has been the application
to the cnstruction of hyperbolic manifolds. 
As with any large and seemingly intractable collection
of mathematical objects, they are best observed through the eyes of invariants.
For hyperbolic manifolds, the most fundamental and important invariant
is hyperbolic volume, which %has the remarkable property of ``discretising''
``discretises'' the situation: the set of possible volumes of hyperbolic manifolds in a fixed
dimension form a well-ordered subset of $\R$, discrete even, if the
dimension is not three. Moreover, the map $\vol:\II_n\rightarrow\R$ from 
the set of isometry classes of finite volume hyperbolic $n$-manifolds
is finite to one when $n\not=2$. 

As one of the first general results along these lines was obtained by Carl Ludwig Siegel
\cite{Siegel45} in 1945, the \emph{Siegel\/} (respectively \emph{weak Siegel\/})
\emph{problem\/} in dimension $n$ is: what is the minimum possible volume obtained 
by a complete hyperbolic $n$-orbifold (resp. $n$-manifold) without boundary?
Our interest is exclusively in the weak problem, which, despite the diminutive
(which we will drop from now on),
is one with a long history. It can probably be most concisely
summarized by saying that its solution in $2$-dimensions is classical;
in $3$-dimensions a solution was announced recently in \cite{Gabai07}
(at least in the closed 
orientable case);
in $4$-dimensions it was first done in \cite{Ratcliffe00} (see also
\cite{Everitt04a}) and a solution in $6$-dimensions was recently announced 
in \cite{Everitt05}. %There are only conjectured solutions in $3$ and $5$-dimensions.

Many of these solutions have a heavy computational element.
Now it turns out that among the Coxeter groups in the large class to which our
construction applies, there are a number of geometric examples, ie: Coxeter
groups generated by reflections in the faces of a Coxeter polytope in
hyperbolic space. When specialized to these cases, our results yield
manifolds that %, while not providing solutions to the Siegel problem, 
are of very small volume, often only twice the volume of a possible solution to the
Siegel problem. They also have the virtue of being constructed algebraically rather 
than computationally.
% More importantly, they are constructed without recourse to a
% computer.

The paper has been written so as to include in its audience workers in 
the field of geometric manifolds who may be less familiar with the
technical details of Weyl groups and Coxeter 
groups. The first section is therefore a summary of the basic definitions
and some of the results we will need later. Sections 
\ref{section:lattices}, \ref{section:homomorphisms} and
\ref{section:extended}, while the technical heart of the paper, can 
be glossed over on the first reading by those impatient to get to the
geometrical consequences in \S\ref{section:manifolds}. The main results are Theorems 
\ref{section:homomorphisms:result300}-\ref{section:extended:result500}.

\section{Preliminaries on reflection groups}\label{section:preliminaries}

%For the readers convenience 
%we collect together the basic results on Coxeter groups.
%we will need for the readers convenience. 
General references for this section %for Coxeter groups and reflection groups 
are \cite{Bourbaki02,Humphreys90,Kane01,Vinberg85}; the 
terminology below mostly follows \cite{Humphreys90}.

\subsection{Coxeter groups and reflection groups}\label{subsection:coxeter_groups}

Let $W$ be a group and $S\subset W$ finite. The pair $(W,S)$ is
called a Coxeter group if $W$ admits a
presentation with generators the
$s\in S$, and relations,
\begin{equation}\label{equation:coxeter_relations}
(st)^{m_{st}}=1,
\end{equation}
where $m_{st}=m_{ts}\in\Z^+\cup\{\infty\}$ and $m_{st}=1$ if and
only if $s=t$.
%It is customary to omit relations
%for which $m_{st}=\infty$.
%
Associated to any Coxeter group is its symbol $\Gamma$, with nodes the
$s\in S$, and where nodes $s$ and $t$ are joined by an edge
labeled $m_{st}$ if
$m_{st}\geq 4$ and an unlabeled edge if 
$m_{st}=3$ %and a dotted edge if $m_{st}=\infty$ 
(and so nodes
$s$ and $t$ not connected %by an edge 
correspond to $m_{st}=2$). 
Conversely, any symbol $\Gamma$ with edges so labeled gives rise to Coxeter group
$W(\Gamma)$, with generators the nodes of the symbol and relations 
(\ref{equation:coxeter_relations}) with $m_{st}$ the label of the edge
connecting nodes $s$ and $t$. %We will tend to write our Coxeter groups
%as $W(\Gamma)$ when we want to emphasise the symbol. 
%An alternative symbol convention is to join nodes $s$ and $t$ by
%$m_{st}-2$ unlabeled edges, and

If $T\subset S$ then the subgroup generated by the $s\in T$ is also a 
Coxeter group with symbol $\Delta\subset\Gamma$ obtained by removing the nodes
not in $T$ and any edges incident with them. Write $W_T$ or $W(\Delta)$ for this, a
\emph{visible\/} subgroup of $W(\Gamma)$ (with parabolic, special parabolic and visual being
alternative terms found in the literature). 
The Coxeter group $W(\Gamma)$ is
irreducible when $\Gamma$ is connected, %ie: there is a path of edges connected
%any two nodes, where the $m_{st}>2$ for all edges in the path. 
and if $\Gamma$ has connected components $\Gamma_i\subset\Gamma$, then 
the visible $W(\Gamma_i)$ are called the irreducible components of $W(\Gamma)$, 
with $W(\Gamma)$ the direct product of
the $W(\Gamma_i)$.

A Coxeter group $W(\Gamma)$ can be realized as a reflection group
via the \emph{reflectional\/}
(or geometric, or Tits) representation:
let $V$ be the real vector space with basis 
$\{v_s\,|\,s\in S\}$ and symmetric bilinear form defined by,
$$
B(v_s,v_{t})=-\cos\frac{\pi}{m_{st}}.
$$
When $m_{st}=\infty$, 
set the product to be some real number $c_{st}$ with
$c_{st}\leq -1$.
For $u\in V$, define $\ss_u:V\rightarrow V$
by $\ss_u(v)=v-2B(v,u)u$,
and then the map $s\mapsto\ss_s:=\ss_{v_s}$ extends to a faithful 
representation
%\sidecomment{this needs to be checked for the more general reflectional representation}
$\ss:W(\Gamma)\rightarrow\gl(V)$, irreducible when $W(\Gamma)$ is 
irreducible and $B$ non-degenerate. %(see, eg: \cite{Humphreys90}*{Corollary 5.4}).
We will abbreviate $\ss(g)(v)$ to $g(v)$.
%The reflectional representation turns out to be irreducible \cite{Bourbaki02}*{\S?}.
%and as we shall see below, realizable over $\Q$. 

The two most important classes of Coxeter groups for us will be the finite and
hyperbolic: $W(\Gamma)$ is finite if and only if the form $B$ 
is positive
definite, in which case $V$ is a Euclidean
space with $W(\Gamma)$-invariant inner product $(v,u):=B(v,u)$. 
%The finite Coxeter groups were classified by Coxeter \cites{Coxeter35,Coxeter34}.
%The $W(\Gamma)$-orbit of
%the vectors $v_s (s\in S)$ then forms a root system $\Phi\subset V$:
%a finite set of non-zero vectors such that, (i). if $v\in\Phi$ then
%$\lambda v\in\Phi$ if and only if $\lambda=\pm 1$, and (ii). if $u,v\in\Phi$ then
%$\ss_u(v)\in\Phi$. 
%The Coxeter symbol $\Gamma$ is crystallographic 
%(and $W(\Gamma)$ is a \emph{Weyl group\/}) if
%\begin{equation}\label{equation:crystallographic}
%\langle u,v\rangle\stackrel{\text{def}}{=}\frac{2(u,v)}{(u,u)}\in\Z,
%\end{equation}
%for all $u,v\in\Phi$. We will have much more to say about Weyl groups in
%\S\ref{subsection:weyl_groups}.

A hyperbolic Coxeter polytope is an intersection
$P=\bigcap_{s\in S}H_s^-$,
of closed half spaces $H_s^-$ 
in $n$-dimensional hyperbolic space $\H^n$, bounded by hyperplanes $H_s$,
any two of which are disjoint or subtend a dihedral angle of $\pi/m_{st}$.
The group generated by reflections in the faces of $P$ is then a Coxeter
group, and an abstract $W(\Gamma)$ is $n$-dimensional hyperbolic in this sense
if there exist
$c_{st}$'s such that the form $B$ above has 
signature
%\footnote{The form has signature $(p,q)$ if the matrix
%with $(s,t)$-th entry $B(v_s,v_t)$ has $p$ eigenvalues that
%are $<0$ and $q$ that are $>0$.} 
$(1,n)$. In this case $c_{st}=-\cosh\eta_{st}$, for $\eta_{st}$ the distance between non-intersecting
$H_s$ and $H_t$.

%The configuration of the %resulting collection of 
%bounding hyperplanes of $P$
%can then be described as follows. If $m_{st}\in\Z^+$ labels the edge
%of the symbol connecting nodes $s,t$, then the
%hyperplanes $H_s,H_{t}$ intersect with dihedral angle
%$\pi/m_{st}$; 
%%\marginlabel{\small{$\blob$ check this angle with the scaled form}}
%if $c_{st}=-1$ they are parallel, that is, intersect at
%infinity on the boundary $\partial\H^n$ of hyperbolic space; 
%and if $c_{st}<-1$ they are ultraparallel (they do not
%intersect in the closure of $\H^n$) with a common perpendicular geodesic
%of length $\eta_{st}$, where $-\cosh\eta_{st}=c_{st}$. 
%Because of these last two, some authors prefer to embellish the Coxeter
%symbol, replacing the dotted edges by thick solid ones if
%$c_{st}=-1$, or by ones labelled $\eta_{st}$ when $c_{st}<-1$.

\begin{table}
\begin{tabular}{c}
\hline
\begin{pspicture}(0,0)(14.3,4)
%\showgrid
\rput(0,-2){
\rput(0,2){
\rput(2,3){\BoxedEPSF{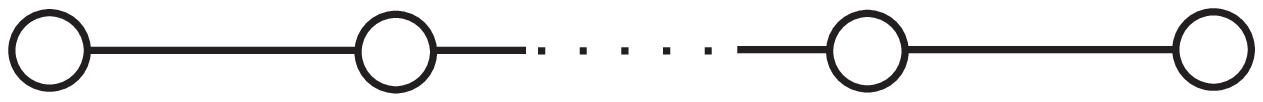 scaled 300}}
\rput(0,0){
\rput(.2,2.65){$1$}\rput(1.275,2.65){$2$}
\rput(2.7,2.65){$n-1$}\rput(3.8,2.65){$n$}
}
\rput(2,3.5){$A_n\,(n\geq 1)$}
\rput(7,3){\BoxedEPSF{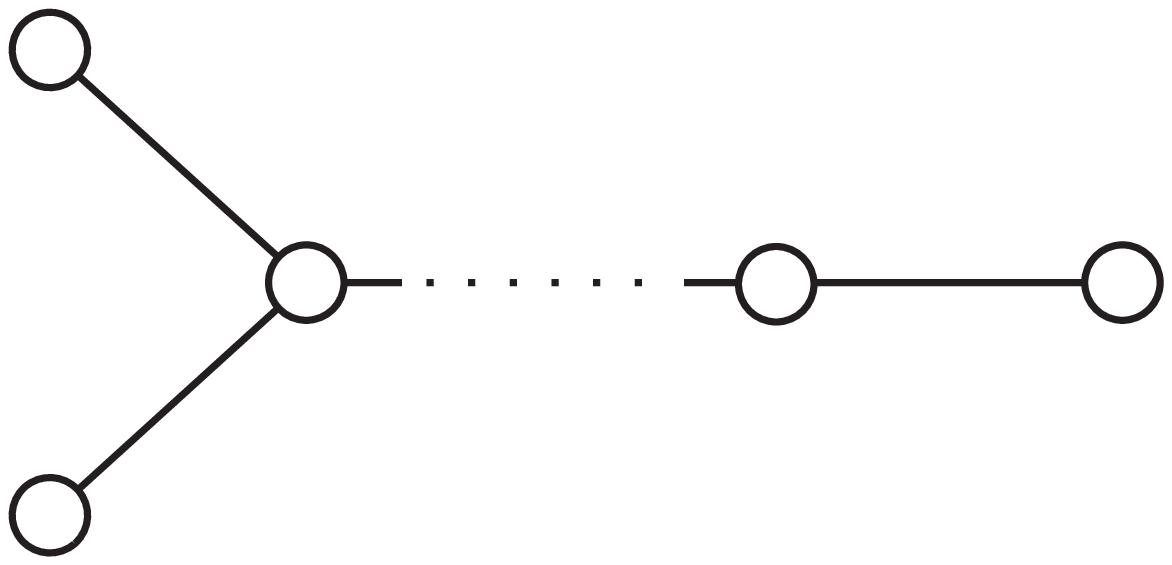 scaled 300}}
\rput(4.85,0){
\rput(.2,3.7){$n$}\rput(-.1,2.3){$n-1$}\rput(1.6,2.65){$n-2$}
\rput(2.725,2.65){$2$}\rput(3.8,2.65){$1$}
}
\rput(7.2,3.5){$D_n\,(n\geq 4)$}
\rput(12,3){\BoxedEPSF{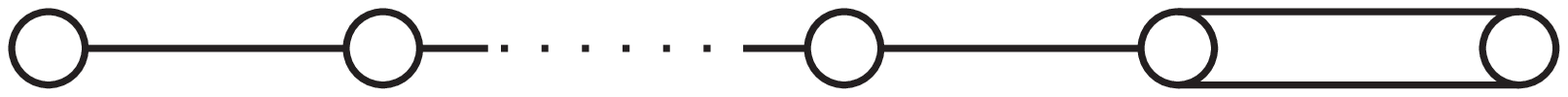 scaled 300}}
\rput(9.5,0){
\rput(.2,2.65){$1$}\rput(1.25,2.65){$2$}
\rput(2.7,2.65){$n-2$}\rput(3.7,2.65){$n-1$}
\rput(4.8,2.65){$n$}
}
\rput(11.8,3.5){$B_n\,(n\geq 2)$}%\rput(13.75,3.2){$4$}
\rput(-.225,-0.2){\rput(3.25,1){\BoxedEPSF{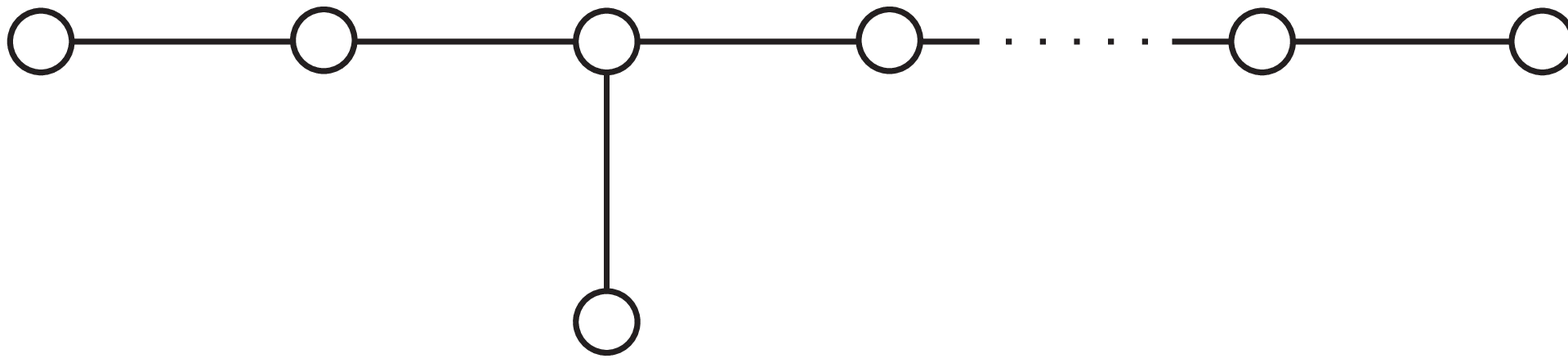 scaled 300}}}
%\rput(-1,0){\rput(3.25,1){\BoxedEPSF{E6.eps scaled 300}}
\rput(-.95,-.15){
\rput(1.15,1.8){$1$}\rput(2.225,1.8){$2$}\rput(3.275,1.8){$3$}
\rput(4.35,1.8){$4$}\rput(5.7,1.8){$n-2$}
\rput(6.75,1.8){$n-1$}\rput(2.95,.45){$n$}
}
\rput(4,.6){$E_n\,(n=6,7,8)$}
%\rput(-1.5,0){\rput(10.25,.975){\BoxedEPSF{E7.eps scaled 300}}
%\rput(11,1){$E_7$}}
}
\rput(0,0){\rput(12.5,3.3){\BoxedEPSF{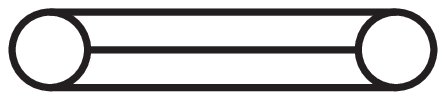 scaled 300}}
\rput(12.5,2.7){$G_2$}}%\rput(13.5,2.7){$6$}
\rput(10.8,1.85){
\rput(1.15,1.8){$1$}\rput(2.225,1.8){$2$}
}
%\rput(1,0){\rput(3.25,1){\BoxedEPSF{E8.eps scaled 300}}
%\rput(4,1){$E_8$}}
\rput(5.5,0){\rput(3.25,3.3){\BoxedEPSF{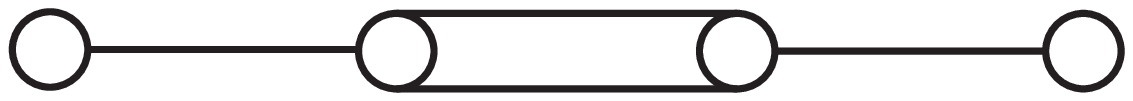 scaled 300}}
\rput(0.5,1.85){
\rput(1.15,1.8){$1$}\rput(2.225,1.8){$2$}
\rput(3.275,1.8){$3$}\rput(4.35,1.8){$4$}
}
\rput(3.2,2.7){$F_4$}%\rput(3.3,1.2){$4$}
}}
\end{pspicture}\\
\hline
\end{tabular}\caption{The irreducible Weyl groups}\label{table:roots1}
\end{table}

\subsection{Crystallographic groups and Weyl groups}\label{subsection:weyl_groups}

A Coxeter group $W(\Gamma)$
is crystallographic if in the reflectional representation
$W(\Gamma)\rightarrow\gl(V)$ there is a $W(\Gamma)$-invariant lattice 
$L\subset V$ (one often sees in the literature a somewhat larger class
of Coxeter groups referred to as crystallographic, namely the Weyl groups of
Kac-Moody Lie algebras).
%(the $\Z$-span
%of a basis for $V$, or a free $\Z$-module of rank $\text{dim}\,V$) that
%is stablized by $W(\Gamma)$. The crystallographic groups have been
%classified \cite{Humphreys90}*{\S 6.6} (although the reader should beware that,
%as remarked in \cite{Humphreys90}*{\S 6.6}, this definition of crystallographic
%Coxeter group differs slightly from\ other versions in the literature, where
%they are taken to be the ``Weyl groups'' of Kac-Moody Lie algebras).
%
A finite crystallographic Coxeter group is a \emph{Weyl\/} group,
and if $(u,v)$ is the inner product on $V$ from \S\ref{subsection:coxeter_groups}, let
$$
\langle u,v\rangle:=\frac{2(u,v)}{(u,u)}.
$$
The symbols for the irreducible Weyl groups are given in Table 
\ref{table:roots1}, using the alternative symbol convention
where nodes $v_s$ and $v_t$ are joined by 
$\langle v_s,v_t\rangle\langle v_t,v_s\rangle$ unlabelled edges.
We will tend to use this version when some additional labeling of the nodes
might otherwise lead to an overly cluttered picture.\footnote{The Weyl groups constitute 
\emph{almost\/} all the finite 
Coxeter groups, with the only omissions being the dihedral groups
and the symmetry groups of the $3$-dimensional dodecahedron/icosahedron and the
$4$-dimensional 120/600-cell. %(as these will not concern is in this paper,
%see \cite{Humphreys90}*{\S2.13} for their symbols)
}
We have fixed a numbering of $S$ in the table, and throughout
the paper we will write $v_i:=v_{s_i}$ for the basis of the ambient space
$V$ of the reflectional representation.
%The Weyl groups have been classified: 
%$W(\Gamma)$ is an irreducible Weyl group if and only
%if $\Gamma$ is one of the symbols in Table
%\ref{table:roots1}, %in which case $W(\Gamma)$ is the
%direct product of these, the irreducible Weyl groups 
%and a general Weyl group is a direct product of irreducible ones
%\cite{Humphreys90}*{\S 2}.
%The rank $\rk\Gamma$ of the Weyl group $W(\Gamma)$ is the number of nodes of
%the symbol $\Gamma$.

We use the Killing-Cartan notation from the classification of
simple Lie algebras to
denote the symbols for the irreducible Weyl groups. 
There are three infinite \emph{classical\/} families,
$A_n (n\geq 1), B_n (n\geq 2)$ and $D_n (n\geq 4)$,
together with five \emph{exceptional\/} groups $G_2, F_4, E_6, E_7$ and $E_8$.
%A quick perusal of Table \ref{table:roots1} reveals 
%another pair of families into which these groups fall, namely
%the \emph{simply-laced\/} groups, for which all the $m_{st}=2$ or $3$, 
%and the remaining, non-simply laced groups (the terminology comes from ...).
The subscripts give the rank of the group, which is just the number of nodes
of the symbol.

%The Weyl groups can be alternatively defined as being those finite Coxeter groups
%for which there is a $W(\Gamma)$-invartiant $\Z$-lattice in the space $V$ of 
%the reflectional representation. Indeed, the \emph{root lattice\/} $L$ of a Weyl group
%$W(\Gamma)$ is the $\Z$-span in $V$, the ambient space of the reflectional
%representation, 
%of the root system $\Phi$ of \S\ref{subsection:coxeter_groups},
%or equivalently, the $\Z$-span of the $v_s$, as the $\Phi$ are $\Z$-linear 
%combinations of the $v_s$ when $\Gamma$ is crystallographic. The result
%is a free $\Z$-module of rank $\rk\Gamma$.

There are a number of $W(\Gamma)$-invariant lattices in $V$, with the
most fundamental being 
the \emph{root lattice\/}, $L(\Gamma)=\Z$-span$\{x_s\,|\,s\in S\}$,
where $x_s=v_s$ in types $A,D$ and $E$.
Otherwise, %if the nodes of $B_n, F_4$
%and $G_2$ in Table \ref{table:roots1} are labelled $\{s_1,\ldots,s_n\}$ from
%left to right, then let 
for $L(B_n)$, let
$x_i=v_{i} (i<n)$ and $x_n=\sqrt{2}v_{n}$; 
in $L(F_4)$ we have
$x_1=v_{1},x_2=v_{2},x_3=\sqrt{2}v_{3},x_4=\sqrt{2}v_{4}$
and $L(G_2)=\Z$-span$\{x_1=v_{1},x_2=\sqrt{3}v_{2}\}$.
%There is nothing unique about these definitions. From now on, when we speak of the
%root lattice $L$ of an irreducible Weyl group, it will be these versions that we
%mean.
%
%One spin-off from the construction of these lattices, absolutely
%central to this paper, is that
%$B(x_i,x_j)\in\Q$ for $x_i,x_j$ basis vectors for $L$, and so
%we realise the geometric 
The $W(\Gamma)$-action on the root lattice realizes the
reflectional representation over $\Q$, and remarkably, any
irreducible representation of a Weyl group can be so realized 
(see \cite{Springer78}).

There is another interesting $W(\Gamma)$-invariant lattice in 
$V$: the \emph{weight lattice\/} $\widehat{L}=\widehat{L}(\Gamma)$ is
given by
$\widehat{L}=\{v\in V\,|\,\langle u,v\rangle
\in\Z\text{ for all }u\in L\}$.
One can also define $\widehat{L}$ indirectly as the dual of the coroot lattice.
% see \cite{Kane01}*{\S9.4}, 
% and indeed show %\cite{Kane01}*{\S 9.5} 
% that any $W(\Gamma)$-invariant
% lattice in $V$ is sandwiched in between
% the root and weight lattices. 
% \sidecomment{check, eg: $\Lambda_i$?}
% Writing $S=\{s_1,\ldots,s_n\}$, $n=\rk\Gamma$, and
% $v_i$ for $v_{s_i}$, 
The weight lattice has basis the simple weights
$\ww_s$ defined by 
$\langle v_s,\ww_s\rangle=1$ and $\langle v_t,\ww_s\rangle=0$
for $t\not= s$. 
In particular
$\widehat{L}$ is a $\Z$-lattice of the same rank as $L$ and with % the 
% crystallographic condition
% forcing the inclusion
$L\subset\widehat{L}$. Thus the root lattice is a subgroup
of finite index in the weight lattice, with this index called
the index of connection (see Table \ref{table:weyl_data}).

The classical Weyl groups have alternative descriptions
as groups of permutations of an orthonormal basis for $V$ 
(or a space closely related to it). For type $A_n$
and basis $\{v_1,\ldots,v_n\}$, let $\widehat{V}$ be 
an $(n+1)$-dimensional Euclidean 
space with orthonormal basis $\{u_1,\ldots,u_{n+1}\}$, and embed $V\hookrightarrow
\widehat{V}$ via,
$$
v_i=\frac{1}{\sqrt{2}}(u_i-u_{i+1})\text{ for }1\leq i\leq n.
$$
The $W(A_n)$-action on $V$ extends to $\widehat{V}$ with
the $\{u_i\}$ a $W(A_n)$-invariant subset and
the map $s_i\mapsto (u_i,u_{i+1})$ inducing an isomorphism
$W(A_n)\rightarrow\SS(u_1,\ldots,u_{n+1})\cong\SS_{n+1}$, the
symmetric group of degree $n+1$.

For types $B$ and $D$, let $\SS^\pm(X)$ be
the group of {\em signed permutations\/} of the set $X$,
ie: $\SS^\pm(X)=\{\ss\in\SS(X\,\cup\,-X)\,|\,\ss(-x)=-(\ss x)\}$. %It has 
%a subgroup of index two, $\SS^\pm_\circ(X)$, the even signed permutations,
%consisting of the signed permutations
%that send an even number of the $X$ to $-X$.
The vectors
$u_1,\ldots,u_n\in V$ defined by
$$
v_i=\frac{1}{\sqrt{2}}(u_i-u_{i+1})\text{ for }i<n,
\text{ and }v_n=u_n,
$$
form an orthonormal basis for $V$ and the 
map $s_i\mapsto (u_i,u_{i+1})(-u_i,-u_{i+1})$, $1\leq i<n$
and $s_n\mapsto (u_n,-u_n)$ induces an isomorphism $W(B_n)\rightarrow
\SS^\pm(u_1,\ldots,u_n)$. Taking instead $v_n=\frac{1}{\sqrt{2}}(u_{n-1}+u_n)$
and $s_n\mapsto(u_{n-1},-u_n)(u_n,-u_{n-1})$, induces an isomorphism
$W(D_n)\hookrightarrow\SS^\pm(u_1,\ldots,u_n)$ onto a subgroup of index two,
the \emph{even\/} signed permutations $\SS^\pm_\circ(u_1,\ldots,u_n)$.

%\centerline
\begin{table}
\begin{tabular}{cccccc}
\hline
$\Gamma$&basis $\{x_i\}$ for $L$&index of connection&$(-1)$-type
&Coxeter number&exponents\\
\hline
$A_n$&$\{v_i\}$&$n+1$&$n=1$&$n+1$&$1,2,\ldots,n$\\
$B_n$&$\{v_i (i<n),\sqrt{2}v_n\}$&$2$&$\checkmark$&$2n$&$1,3,\ldots,2n-1$\\
$D_n$&$\{v_i\}$&$4$&$n$ even&$2(n-1)$&$1,3,\ldots,2n-3,n-1$\\
$G_2$&$\{v_1,\sqrt{3}v_2\}$&$1$&$\checkmark$&$6$&$1,5$\\
$F_4$&$\{v_1,v_2,\sqrt{2}v_3,\sqrt{2}v_4\}$&$1$&$\checkmark$&$12$&$1,5,7,11$\\
$E_6$&$\{v_i\}$&$3$&$\times$&$12$&$1,4,5,7,8,11$\\
$E_7$&$\{v_i\}$&$2$&$\checkmark$&$18$&$1,5,7,9,11,13,17$\\
$E_8$&$\{v_i\}$&$1$&$\checkmark$&$30$&$1,7,11,13,17,19,23,29$\\
\hline
\end{tabular}
\caption{The Weyl group data discussed in
\S\ref{subsection:weyl_groups}}
\label{table:weyl_data}
\end{table}

%We end this section with a number of other fundamental concepts concerning 
%Weyl groups. First, 
Every Coxeter group has a length function
$\ell:W(\Gamma)\rightarrow\Z^{\geq 0}$, where $\ell(g)$ is the smallest $r$
for which an expression $g=s_1s_2\ldots s_r$ exists, with the $s_i\in S$.
In particular %\cite{Humphreys90}*{Theorem 1.8}, 
if $W(\Gamma)$ is finite, then there is a unique $w_\Gamma^{}$ 
with maximal length, which, significantly for us, turns out to be
an involution.
If the $\Gamma_i$ are the connected components of $\Gamma$ and $\ell_i$
the length functions on the $W(\Gamma_i)$, then $\ell=\sum\ell_i$ is the length
function on $W(\Gamma)$.
%whose length is maximal, and called the element of longest length
%and turns out to be half the number of roots in the root system
%$\Phi$ for $W$ 
%\cite{Humphreys90}*{Theorem 1.4}.

A Coxeter group $W(\Gamma)$ is of \emph{$(-1)$-type\/} if in the reflectional 
representation there is an element $g\in W$ acting on $V$ as the 
antipodal map, ie: $g(v)=-v$ for all $v\in V$
(and then $g$ is the element $w_\Gamma^{}$ of longest 
length)\footnote{Geometrically, every Weyl group is the symmetry group
of a solid polytope, indeed a regular polytope in types $A,B$ and $F$.
Then $W$ has $(-1)$-type when there is polytope whose vertices occur
in antipodal pairs. These polytopes are all familiar except possibly
in type $E$ (see
\cite{Coxeter88}).}. 
An irreducible group has $(-1)$-type precisely when it has
non-trivial center; every $W$ of $(-1)$-type is finite; and
$W$ is of $(-1)$-type if and only if its
irreducible components are of $(-1)$ type 
\cite{Bourbaki02}*{Chapter V, \S 4, Exercises 2(b) and 3(c)} or \cite{Richardson82}*{1.9}. 

%This can be interpreted geometrically: every
%finite Coxeter group $W(\Psi)$ is the group of symmetries of a
%$\rk\Psi$-dimensional solid (although amongst the Weyl groups, 
%only types $A,B,F$ and $G$ are the symmetry
%groups of \emph{regular\/} solids) and the group has $(-1)$-type
%precisely when one can find such a solid with vertices occuring 
%in antipodal pairs (these solids are well known except possibly in type 
%$E$; see \cite{Coxeter88}). 
%%Unfortunately we will have no further use of this pleasant geometrical fact.
%Thus, type $A$ has $(-1)$-type only for $n=1$, as this solid is the
%$n$-simplex, whereas type $B$ (the $n$-cube) has $(-1)$-type
%for all $n$. The full list is given in Table \ref{table:weyl_data}
%(see \cite{Richardson82}*{1.12}).

A product involving all the 
elements of $S$ precisely once is called a Coxeter element,
and if the symbol $\Gamma$ is a tree (in
particular, if we have a Weyl group) then any two Coxeter elements are
conjugate.
In particular, the Coxeter elements all have the same order,
the Coxeter number $h$. 
In the reflectional representation, the Coxeter elements have
eigenvalues, 
$\zeta^{m_1},\zeta^{m_2},\ldots,\zeta^{m_n}$,
with $\zeta$ a primitive $h$-th root of unity, and 
$0\leq m_1\leq m_2\leq\cdots\leq m_n<h$ 
the \emph{exponents\/} of $W(\Gamma)$. %and their explicit calculation
%is somewhat delicate (see \cite{Humphreys90}*{\S\S 3.19-3.20} 
%and Table \ref{table:weyl_data}). One of the truly remarkable properties of the 
%Weyl groups is that by increasing all the exponents by $1$, the \emph{degrees\/}
%of the group are obtained: these being the degrees of the polynomials in a set
%of $\rk\Psi$ homogenous, algebraically independent generators for the algebra
%of $W$-invariant polynomials.
%
Every integer $1\leq m<h$ relatively prime to $h$ is an exponent
(although there may be others), 
and the order of a Weyl group is the product of the $m_i+1$.

%\begin{lemma}\label{section:result:parity}
%Let $(W,S)$ be a Coxeter group and $s\in S$ such that $m_{ss'}$ is even for
%all $s'\not= s$. Then for any $g\in W$ and words $w_1,w_2$ in the $S$ for $g$,
%the parity of the number of occurences of $s$ in $w_1$ is the same 
%as the parity of the number of
%occurences of $s$ in $w_2$.
%\end{lemma}
%
%One wouldn't expect this to be true for an arbitrary Coxeter generator. 
%For example, if $W$ is the Weyl group of type $A_2$ with Coxeter generators
%$s_1,s_2$, then the element of longest length $w_0$ has expressions
%$s_1s_2s_1$ and $s_2s_1s_2$.

%For one consequence of the Lemma, consider the Weyl group $W(B_n)$ 
%$$
%\begin{pspicture}(0,0)(14,1)
%%\showgrid
%\rput(7,.5){\BoxedEPSF{Bn.eps scaled 400}}
%%\rput(14,.5){$(\dag)$}
%\rput(3.95,.1){$s_1$}\rput(5.35,.1){$s_2$}\rput(7.2,.1){$s_{n-2}$}
%\rput(8.6,.1){$s_{n-1}$}\rput(10.1,.1){$s_n$}
%\rput(9.4,.75){$4$}
%\end{pspicture}
%$$
%Then any two expressions for the element $w_{B_n}^{}$ of longest length have the same 
%parity number of
%occurences of $s_n$. It is well known (see, eg: \cite{Haiman92}) that 
%is such an expression, hence this parity is the same as that of the rank $n$.

\subsection{Torsion in Coxeter groups}\label{subsection:torsion}

In \S\S\ref{section:homomorphisms}-\ref{section:extended}
we will want to show that certain subgroups
of Coxeter groups are torsion free, ie: contain no non-trivial elements
of finite order, and the arguments
will require a detailed knowledge of the location in
a Coxeter group of the torsion.
The first key result gives a rough description:

\begin{theorem}[\cite{Bourbaki02}*{V.4.2}, \cite{Brink93}*{Proposition 1.3}]
\label{torsion} 
Any element of finite order in a Coxeter group $W(\Gamma)$ is conjugate to
an element of a finite visible subgroup $W(\Delta)$ for $\Delta\subset\Gamma$.
\end{theorem}

See also \cite{Everitt04a}*{Theorem 1} for a geometrical proof. In particular,
if $f:W(\Gamma)\rightarrow G$ is a homomorphism, then $\ker f$
is torsion free precisely when $f$ is faithful
on the finite visible $W(\Delta)$ for $\Delta\subset\Gamma$.
Indeed, 
if $\Delta_i\subset\Gamma$ are disjoint symbols with $fW(\Delta_1)
\cap fW(\Delta_2)=\{1\}$,
then $f$ is faithful on $W(\Delta_1\coprod\Delta_2)$ if and only if $f$ is faithful
on the $W(\Delta_i)$, so it is often possible to reduce consideration to the
irreducible finite visibles.

Clearly a group is torsion free if it is free of $p$-torsion 
(elements of order $p$) for all primes $p$. 
The following is then typical of the results we will need, and
follows from \cite{Carter72}*{\S 7} or \cite{Everitt04a}*{Theorem 5}.

\begin{lemma}\label{section:preliminaries:result200}
%Let $\Psi$ be the Coxeter symbol of type $B_n$,
%$$
%\begin{pspicture}(0,0)(14,1)
%%\showgrid
%\rput(7,.5){\BoxedEPSF{Bn.eps scaled 400}}
%%\rput(14,.5){$(\dag)$}
%\rput(3.95,.1){$s_1$}\rput(5.35,.1){$s_2$}\rput(7.2,.1){$s_{n-2}$}
%\rput(8.6,.1){$s_{n-1}$}\rput(10.1,.1){$s_n$}
%\rput(9.4,.75){$4$}
%\end{pspicture}
%$$
Let $g\in W(B_n)$ be an element of $p$-torsion for $p>2$ a prime. 
Then $g$ is conjugate
to an element of the visible $W(A_{n-1})=\langle s_1\ldots,s_{n-1}\rangle$.
\end{lemma}

Because of the lemma, it will turn out that most of
our efforts will be focused on the $2$-torsion, for which
there
is a more precise
classification due to Richardson \cite{Richardson82} which we now summarize 
(see also \cites{Howlett80,Deodhar82}). 
First, if $W(\Psi)$ is an irreducible Weyl group,
define a permutation $\pi_\Psi$ of the nodes of $\Psi$ as follows: 
if $\Psi$ is of $(-1)$-type then
$\pi_\Psi$ is the identity permutation, otherwise, it is the unique symmetry of
order $2$ of the diagram.
Now let $\Gamma$ be a symbol for any $(W,S)$ and $\Delta\subset\Gamma$ a subsymbol of
$(-1)$-type with nodes $T\subset S$; %and let $\AA(\Delta)$ be those nodes 
%$s\in S\setminus T$ such that the connected component of $T\cup\{s\}$ is \emph{not\/}
%a symbol of $(-1)$-type.
for $s\in S$, let $\Delta(s)$ be the connected component of $T\cup\{s\}$ containing
$s$; let $\AA(\Delta)$ be those nodes
$s\in S\setminus T$ such that $\Delta(s)$ is \emph{not\/} a symbol of $(-1)$-type.

Two subsymbols $\Delta,\Delta'\subset\Gamma$ are connected by
an elementary $W(\Gamma)$-equivalence, written $\Delta\vdash_s\Delta'$,
if $\Delta'$ has nodes $T\cup\{s\}\setminus\{s'\}$, where 
$s\in\AA(\Delta)$ and $s'=\pi_{\Delta(s)}(s)$.
The two symbols are $W(\Gamma)$-equivalent
if there is a sequence 
$\Delta=\Delta_0\vdash_{s_1}\Delta_1\vdash_{s_2}\cdots\vdash_{s_k}\Delta_k=\Delta'$
of elementary equivalences between them.

\begin{theorem}[\cite{Richardson82}*{Theorem A}]
\label{section:preliminaries:torsion:result300}
The map $\Delta\mapsto w_\Delta$ induces a bijection from the set $\JJ_\Gamma$
of $W(\Gamma)$-equivalence classes of subsymbols $\Delta\subset\Gamma$
of $(-1)$-type to the set of conjugacy classes of involutions in
the Coxeter group $W(\Gamma)$.
\end{theorem}

For example, if $\Gamma=A_n$ ($n$ even) with
subsymbols $\Delta,\Delta'$ the nodes $\{s_i\,|\,i\text{ odd}\}$ and
$\{s_i\,|\,$ $i\text{ even}\}$ respectively, then there is a series of elementary 
$W(\Gamma)$-equivalences that successively, working from right to left,  
move the nodes of $\Delta$ one place to the right, and so realizing a 
$W(\Gamma)$-equivalence between $\Delta$ and $\Delta'$.

As a corollary, if $W(\Gamma)$ is an irreducible Weyl group, then there
is a \emph{unique\/} class 
in $\JJ_\Gamma$ containing symbols of \emph{maximal\/} rank. Indeed,
as an inspection of Table \ref{table:roots1} shows, this class contains
a single symbol, except when $\Gamma=A_n$ ($n$ even), when it contains
the two symbols $\Delta,\Delta'$ of the previous paragraph.

\begin{proposition}
\label{section:preliminaries:torsion:result400}
Let $W(\Gamma)$ be an irreducible Weyl group with Coxeter number $h$ 
even, $\xi\in W(\Gamma)$ a Coxeter element, % with even Coxeter number, and
and $\Delta\in\JJ_\Gamma$ of maximal rank.
Then $\xi^{h/2}$ is conjugate to $w_\Delta$.
\end{proposition}

\begin{proof}
By checking separately the cases, $\Gamma$ of $(-1)$-type, 
$\Gamma=A_n$ ($n$ odd), $\Gamma=D_n$ ($n$ odd),
and $\Gamma=E_6$, one can show, using the exponents in
Table \ref{table:weyl_data}, that the eigenvalues of $\xi^{h/2}$
are $-1,\ldots,-1,1\ldots,1$, with the rank of $W(\Delta)$ number of $-1$'s.
%(of an older result) Since $w_\Psi$ is central in $W(\Psi)$ 
%it suffices to show that the two elements are conjugate,
%and in particular, that all $r=\rk\Psi$ eigenvalues of $\xi^{h/2}$
%are equal to $-1$. Recall that the eigenvalues of $\xi$ are
%$\beta^{m_1},\beta^{m_2},\ldots,\beta^{m_r}$, with $\beta$ a primitive
%$h$-th root of $1$ and the $0\leq m_1\leq m_2\leq\cdots\leq m_r<h$.
%If $W(\Psi)$ is of $(-1)$-type then a quick perusal of Table
%\ref{table:weyl_data} reveals that $h$ is even and all the exponents $m_i$ are odd,
%and in particular
%$m_i\cdot h/2$ cannot be a multiple of $h$.
%For $\beta$ an $h$-th root of $1$ we have $\beta^{h/2}$ is either $1$
%or $-1$, so $(\beta^{m_i})^{h/2}=-1$.
\qed
\end{proof}

%\section{Some technical results on root lattices}\label{section:lattices}
\section{Technical results on root lattices}\label{section:lattices}

This section studies in some detail the action of a Weyl group on 
certain vector spaces over the field $\F_2$ of order two. 
The reader 
who is more interested in the broad flow of the paper than this somewhat 
intricate combinatorics (one imagines this describes most readers) can
safely skip the details: the key ideas are the sublattices $\Lambda_s\subset L$
defined after Lemma \ref{subsection:weyl_groups:result100}, 
the notions of admissibility and special admissibility after 
Proposition \ref{subsection:weyl_groups:result175},
and Theorems \ref{section:lattices:result200}-\ref{section:lattices:result300}.
These all play an essential role in \S\ref{section:homomorphisms}, 
while Theorem \ref{section:lattices:result400} and 
Proposition \ref{section:lattices:result500} are used in \S\ref{section:extended}.

\begin{lemma}\label{subsection:weyl_groups:result100}
Let %$\Psi$ be %a simply-laced connected crystallographic Coxeter symbol
%and $W(\Psi)\rightarrow\gl(V)$ the reflectional representation of the
%Weyl group $W(\Psi)$. 
%the symbol for an %simply-laced 
$W(\Psi)$ be an irreducible Weyl group of rank $n$, and
$W(\Psi)\rightarrow\gl(V)$ the reflectional representation.
Then for any node $s\in\Psi$, there is a $u_s$ in the
root lattice $L=\Z$-span$\{x_i\}$ such that 
$u_s\in
\langle x_{1}\ldots,\widehat{x}_s,\ldots,x_{n}\rangle^\perp$.%, $n=\rk\Psi$.
\end{lemma}

\pagebreak

\begin{proof}
\hspace*{1em}\begin{enumerate}
\item \emph{slick proof}: let $u_s=|\widehat{L}/L|\,\ww_s$,
with $\ww_s$ the simple weight corresponding to $s$ %via 
%$\langle v_s,\ww_s\rangle=1$ and $\langle v_t,\ww_s\rangle=0$ for $t\not=s$,
and $|\widehat{L}/L|$ the index of connection.

\item \emph{simple-minded proof (but useful for later calculations)}:
Depict 
$u_s=\sum_{t\in S}\aa_{t}x_{t}$ using the 
symbol $\Psi$ by labeling the node $t$ by $\aa_{t}$ and colouring
black the distinguished node corresponding to $s$. 
% If $m_{st}=3,4$ or $6$, let $\lambda_{st}=1,2$ or $3$ respectively.
Letting $\lambda_{st}=\langle v_s,v_t\rangle\langle v_s,v_t\rangle$, we have that 
$u_s$ has the desired properties if and only if
for every $t\in S\setminus\{s\}$, 
\begin{equation}\label{label_condition1}
2\aa_{t}=\sum \lambda_{s't}\aa_{s'},
\end{equation} 
where the sum
is over all $s'\in S$ with $m_{s't}>2$ (ie: all $s'$ connected by an edge to $t$),
\emph{except\/} if $t$ is the node in $B_n,F_4$ or $G_2$ 
with an $m_{st}=4$ or $6$, and $x_t=\sqrt{2}$ or $\sqrt{3}$ times the length
of $v_t$. In this case we require
\begin{equation}\label{label_condition2}
2\aa_t=\sum \aa_{s'}
\end{equation}
instead.
Thus it suffices to show that for any
connected $\Psi$ and distinguished black node $s$, the vertices
can be $\Z$-labeled satisfying conditions
(\ref{label_condition1}) and (\ref{label_condition2}).
We can almost always reduce to the case where $s$ lies at the end of $\Psi$,
%We can assume that the distinguished black node lies at the
%end of the symbol, 
for suppose $\Psi_1,\Psi_2$ %are
%connected crystallographic with 
have distinguished black nodes at the end,
labeled $\aa_i$ and 
satisfying 
(\ref{label_condition1}) and (\ref{label_condition2})
above. %(\ref{label_condition}),
Multiply the labels of $\Psi_i$ by $\aa_j/\text{gcd}(\aa_1,\aa_2)$ 
%$\Psi_2$ by $\aa_1/\text{gcd}(\aa_1,\aa_2)$,
and fuse the black nodes:
%\marginlabel{\small{$\aa$'s zero case}}
$$
\begin{pspicture}(0,0)(10.2,2)
%\showgrid
\rput(0,-.2){
\rput(2.5,1.75){\BoxedEPSF{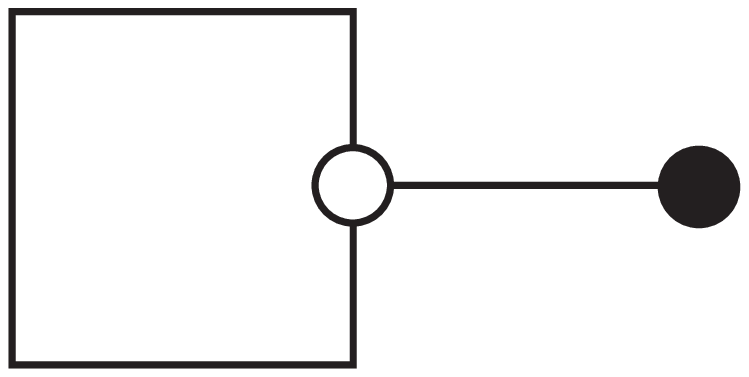 scaled 250}}
\rput(2.5,.5){\BoxedEPSF{fig2a.eps scaled 250}}
\rput(1.2,1.75){$\Psi_1$}\rput(1.2,0.5){$\Psi_2$}
\rput(3.7,1.75){$\aa_1$}\rput(3.7,0.5){$\aa_2$}
\rput{0}(4,1.125){$\left.\begin{array}{c}
\vrule width 0 mm height 17 mm depth 0 pt\end{array}\right\}$}
\rput(-1,0){\psline[linewidth=.2mm]{->}(5.4,1.125)(6.4,1.125)}
\rput(-2.4,0){\rput(10,1.125){\BoxedEPSF{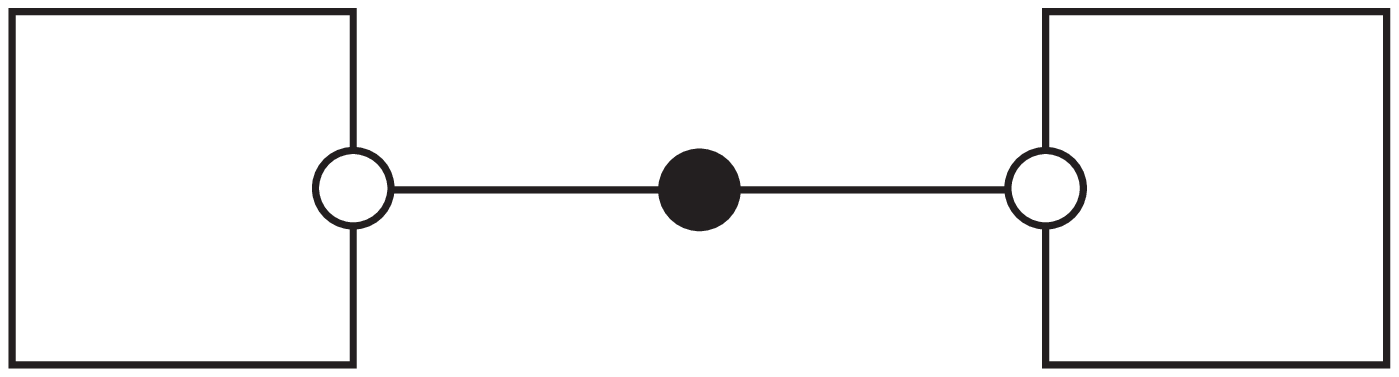 scaled 250}}
\rput(10,0.6){$\Psi$}\rput(10,2){$\text{lcm}(\aa_1,\aa_2)$}
\psline[linewidth=.2mm]{->}(10,1.8)(10,1.3)}
}
\end{pspicture}
$$
The new symbol $\Psi$ has a $\Z$-labeling 
satisfying (\ref{label_condition1}) and (\ref{label_condition2})
as long as in passing from the unfused to the
fused symbols, no node moves from having to satisfy
(\ref{label_condition1}) to having to satisfy (\ref{label_condition2})
or vice-versa. But the only time this happens is for the single pair
$(F_4,s_2)$, which we deal with separately. Thus, we just need to exhibit labelings
in the other cases, where $s$ is at the end of $\Psi$.
In types $A$ and $G$ they are,
$$
\begin{pspicture}(0,0)(12,1)
%\showgrid
\rput(0,0){%moves whole figure
\rput(0,-4.5){%moves A
\rput(1.85,5){\BoxedEPSF{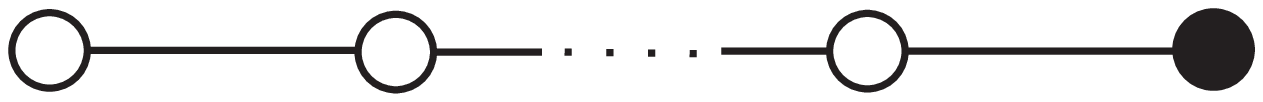 scaled 300}}
%\rput(2,5.5){$A_n$}
\rput(0.1,4.6){${\scriptstyle 1}$}\rput(1.1,4.6){${\scriptstyle 2}$}
\rput(2.6,4.6){${\scriptstyle n-1}$}\rput(3.6,4.6){${\scriptstyle n}$}
}
%\rput(5,-4.5){%moves F
%\rput(1.85,5){\BoxedEPSF{F43_dynkin.eps scaled 300}}
%%\rput(.65,5.5){$F_4$}
%%\rput(1.85,5.2){$4$}
%\rput(0.25,4.6){${\scriptstyle 2}$}\rput(1.3,4.6){${\scriptstyle 4}$}
%\rput(2.35,4.6){${\scriptstyle 3}$}\rput(3.4,4.6){${\scriptstyle -2}$}
%}
\rput(5,-4.5){%moves G
\rput(1.85,5){\BoxedEPSF{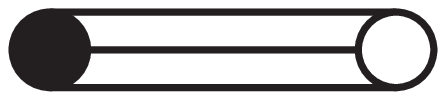 scaled 300}}
%\rput(1,5.5){$G_2$}
%\rput(1.85,5.2){$6$}
\rput(1.3,4.6){${\scriptstyle 2}$}\rput(2.4,4.6){${\scriptstyle 1}$}
}
\rput(9,-4.5){%moves G
\rput(1.85,5){\BoxedEPSF{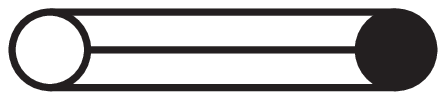 scaled 300}}
%\rput(1,5.5){$G_2$}
%\rput(1.85,5.2){$6$}
\rput(1.3,4.6){${\scriptstyle 3}$}\rput(2.4,4.6){${\scriptstyle 2}$}
}
}
\end{pspicture}
$$
The asymmetry in $G_2$, despite the symmetry of the 
symbol, is due to the asymmetry of the basis $\{x_1,x_2\}$ for the $G_2$-root
lattice. In type $F$, we have,
$$
\begin{pspicture}(0,0)(14,1)
%\showgrid
\rput(0,0){%moves whole figure
\rput(.7,-4.5){%moves F
\rput(1.85,5){\BoxedEPSF{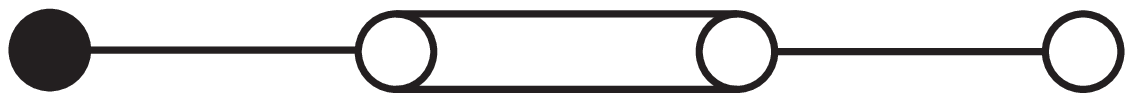 scaled 300}}
%\rput(.65,5.5){$F_4$}
%\rput(1.85,5.2){$4$}
\rput(0.25,4.6){${\scriptstyle 2}$}\rput(1.3,4.6){${\scriptstyle 3}$}
\rput(2.35,4.6){${\scriptstyle 2}$}\rput(3.4,4.6){${\scriptstyle 1}$}
}
\rput(5.2,-4.5){%moves F
\rput(1.85,5){\BoxedEPSF{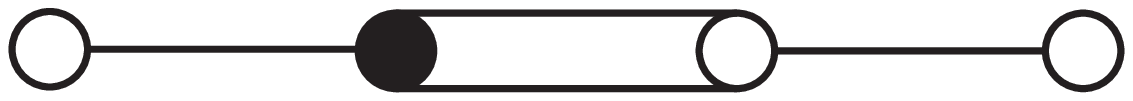 scaled 300}}
%\rput(.65,5.5){$F_4$}
%\rput(1.85,5.2){$4$}
\rput(0.25,4.6){${\scriptstyle 3}$}\rput(1.3,4.6){${\scriptstyle 6}$}
\rput(2.35,4.6){${\scriptstyle 4}$}\rput(3.4,4.6){${\scriptstyle 2}$}
}
\rput(9.65,-4.5){%moves F
\rput(1.85,5){\BoxedEPSF{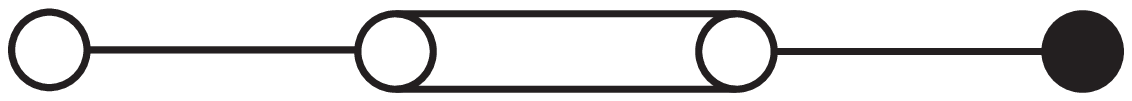 scaled 300}}
%\rput(.65,5.5){$F_4$}
%\rput(1.85,5.2){$4$}
\rput(0.25,4.6){${\scriptstyle 2}$}\rput(1.3,4.6){${\scriptstyle 4}$}
\rput(2.35,4.6){${\scriptstyle 3}$}\rput(3.4,4.6){${\scriptstyle 2}$}
}
}
\end{pspicture}
$$
In type $D$ we have,
$$
\begin{pspicture}(0,0)(14,1.75)
%\showgrid
% \rput(-4.4,-4.25){%moves D
% \rput(7,5){\BoxedEPSF{Dn3.eps scaled 300}}
% %\rput(7,5.7){$D_n$}
% \rput(8.65,4.6){${\scriptstyle 1}$}\rput(7.6,4.6){${\scriptstyle 2}$}
% \rput(6.5,4.6){${\scriptstyle n-2}$}
% \rput(4.9,4.3){${\scriptstyle \frac{n-2}{2}}$}
% \rput(5.05,5.7){${\scriptstyle \frac{n}{2}}$}
% \rput(7.5,5.5){$n$ even}
% }
\rput(2,.4){
\rput(-4.8,-4.5){%moves D
\rput(7,5){\BoxedEPSF{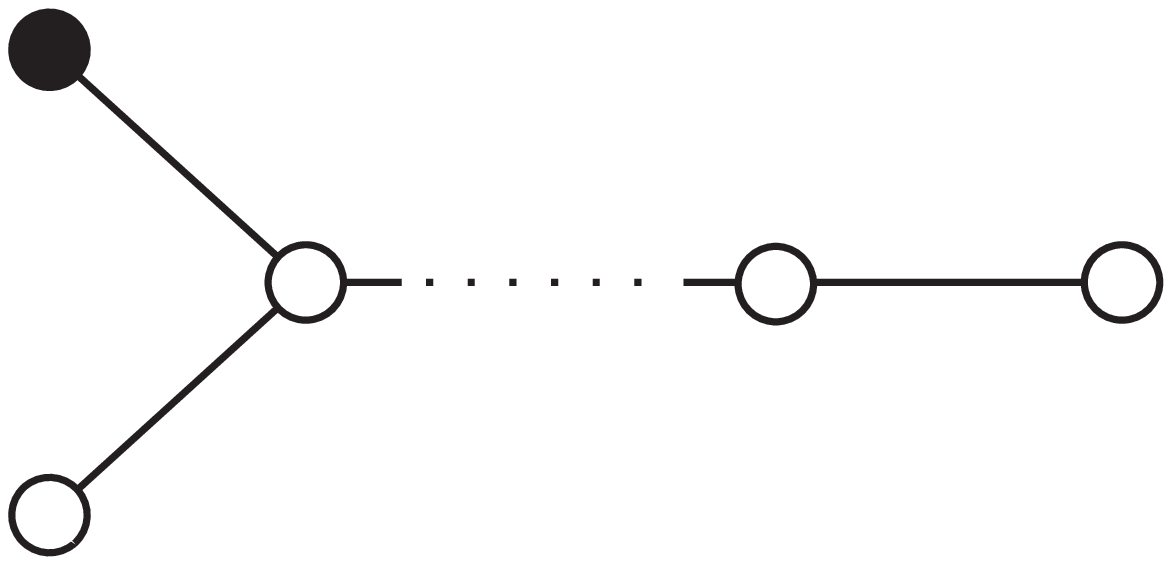 scaled 300}}
%\rput(7,5.7){$D_n$}
\rput(8.65,4.6){${\scriptstyle 2}$}\rput(7.6,4.6){${\scriptstyle 4}$}
\rput(6.5,4.6){${\scriptstyle 2n-4}$}
\rput(5.3,4){${\scriptstyle n-2}$}\rput(5.35,6){${\scriptstyle n}$}
%\rput(7.5,5.5){$n$ odd}
}
\rput(-5.55,-4.5){%moves D
\rput(12,5){\BoxedEPSF{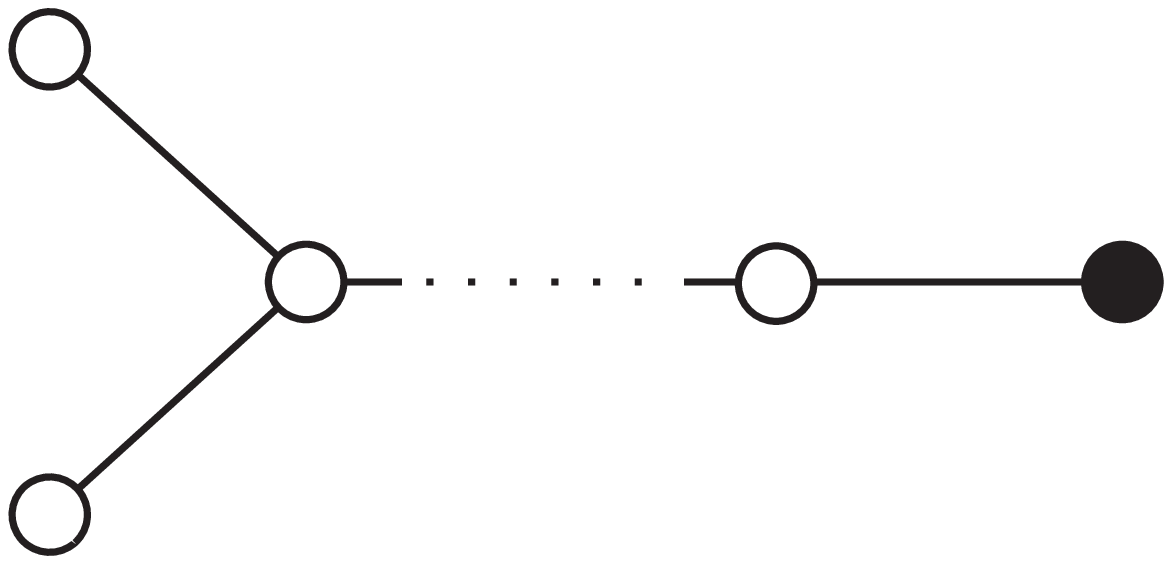 scaled 300}}
%\rput(12,5.7){$D_n$}
\rput(13.65,4.6){${\scriptstyle 2}$}\rput(12.55,4.6){${\scriptstyle 2}$}
\rput(11.1,4.6){${\scriptstyle 2}$}
\rput(10.35,4){${\scriptstyle 1}$}\rput(10.35,6){${\scriptstyle 1}$}
}
}
\end{pspicture}
$$
removing the factor of $2$ on the left when $n$ is even;
in type $B$ we have,
$$
\begin{pspicture}(0,0)(14,1)
%\showgrid
\rput(1.5,0){
\rput(0,0){%moves B
\rput(2.6,.5){\BoxedEPSF{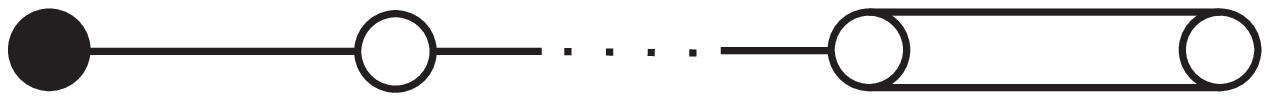 scaled 300}}
%\rput(4.9,.7){$4$}
\rput(.8,.1){${\scriptstyle 2}$}\rput(1.85,.1){${\scriptstyle 2}$}
\rput(3.35,.1){${\scriptstyle 2}$}\rput(4.4,.1){${\scriptstyle 1}$}
}
% \rput(0,0){%moves B
% \rput(8,.5){\BoxedEPSF{Bn5_dynkin.eps scaled 300}}
% %\rput(9.3,.7){$4$}
% \rput(7,1){$n\text{ even}$}
% \rput(6.2,.1){${\scriptstyle 1}$}\rput(7.25,.1){${\scriptstyle 2}$}
% \rput(8.75,.1){${\scriptstyle n-1}$}\rput(9.8,.1){${\scriptstyle \frac{n}{2}}$}
% }
\rput(-5,0){%moves B
\rput(12.5,.5){\BoxedEPSF{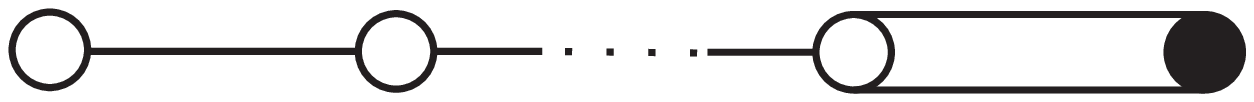 scaled 300}}
%\rput(13.75,.7){$4$}
%\rput(11.5,1){$n\text{ odd}$}
\rput(10.75,.1){${\scriptstyle 2}$}\rput(11.8,.1){${\scriptstyle 4}$}
\rput(13.2,.1){${\scriptstyle 2n-2}$}\rput(14.3,.1){${\scriptstyle n}$}
}
}
\end{pspicture}
$$
again removing factors of $2$, 
and finally in type $E$:
$$
\begin{pspicture}(0,0)(14,2)
%\showgrid
\rput(-.5,0){
\rput(-.8,-2.3){%moves E
\rput(3.45,3){\BoxedEPSF{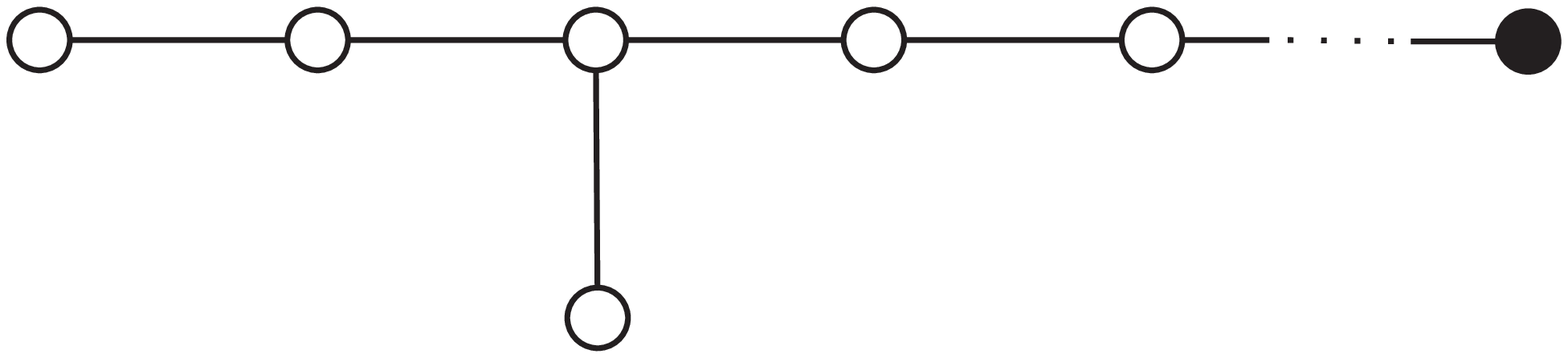 scaled 300}}
%\rput(4.2,3){$E_n$}
\rput(.6,3.9){${\scriptstyle 2}$}\rput(1.7,3.9){${\scriptstyle 4}$}
\rput(2.75,3.9){${\scriptstyle 6}$}\rput(2.4,2.45){${\scriptstyle 3}$}
\rput(3.75,3.9){${\scriptstyle 5}$}\rput(4.9,3.9){${\scriptstyle 4}$}
\rput(6.2,3.9){${\scriptstyle 10-n}$}
}
\rput(-3.35,-1){%moves E
\rput(10.35,1.975){\BoxedEPSF{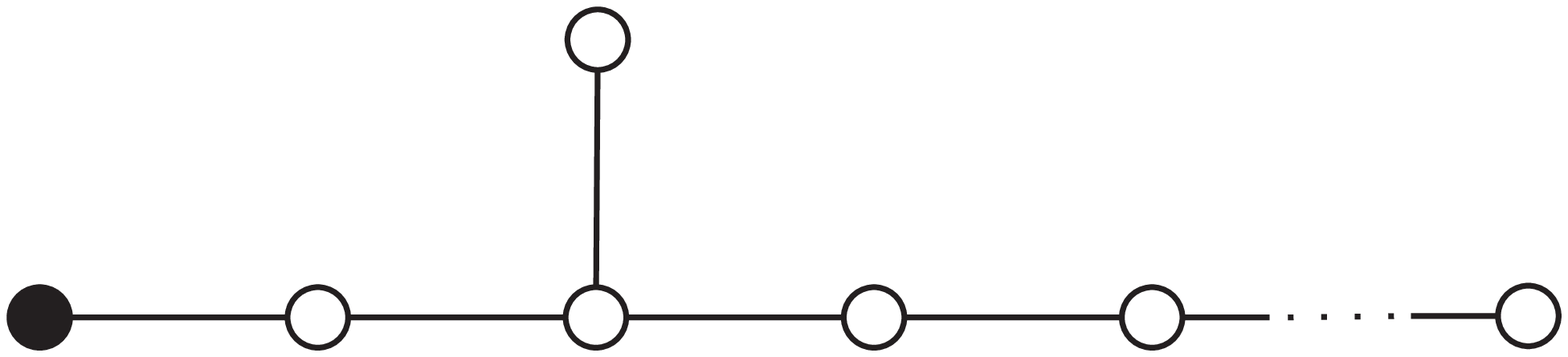 scaled 300}}
%\rput(11.1,2){$E_n$}
\rput(7.9,-3){
\rput(-.375,4.1){${\scriptstyle 4}$}\rput(.65,4.1){${\scriptstyle n-1}$}
\rput(1.7,4.1){${\scriptstyle 2n-6}$}\rput(2.3,5.5){${\scriptstyle n-3}$}
\rput(2.7,4.1){${\scriptstyle 2n-8}$}\rput(3.85,4.1){${\scriptstyle 2n-10}$}
\rput(5.3,4.1){${\scriptstyle 2}$}
}
}
\rput(6.7,0){%moves E
\rput(4.75,0.7){\BoxedEPSF{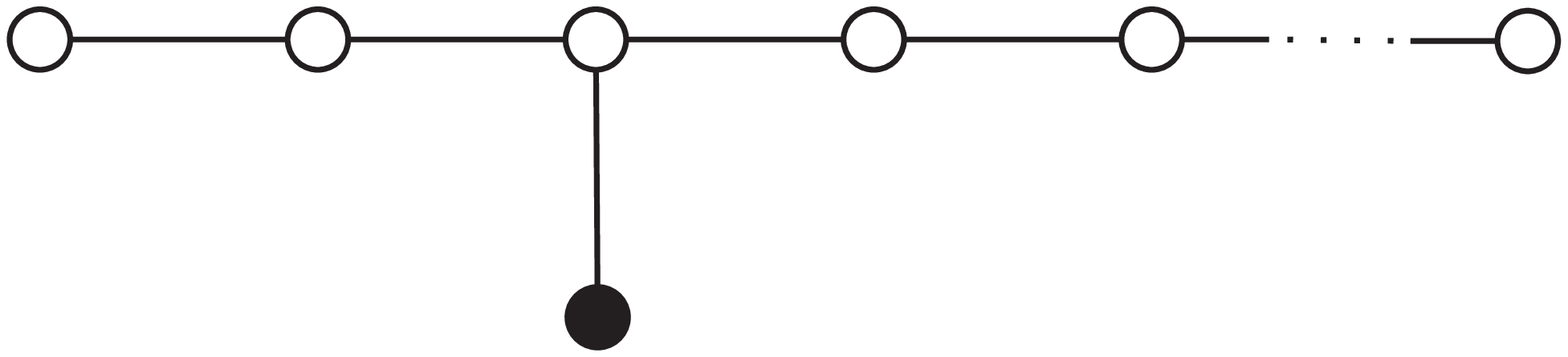 scaled 300}}
%\rput(5,.7){$E_n$}
\rput(1.3,-2.25){
\rput(.6,3.9){${\scriptstyle 1}$}\rput(1.7,3.9){${\scriptstyle 2}$}
\rput(2.75,3.9){${\scriptstyle 3}$}\rput(3.1,2.45){${\scriptstyle 0}$}
\rput(3.75,3.9){${\scriptstyle 4}$}\rput(4.9,3.9){${\scriptstyle 5}$}
\rput(6.2,3.9){${\scriptstyle n-1}$}
}
}
}
\end{pspicture}
$$
In the middle situation for $E_7$, remove the factor of $2$ that appears
in each label.
%The uniqueness upto $\Z$-multiple follows from the
%fact that $\langle v_{s_1}\ldots,$ $\widehat{v}_s,\ldots,v_{s_n}\rangle$ is a 
%hyperplane in $V$.
\end{enumerate}
\qed
\end{proof}

Now for the key player in this section. 
Let $W(\Psi)$ be an irreducible Weyl group of rank $n$,
$s\in\Psi$ and $u_s\in L$ the vector
of Lemma \ref{subsection:weyl_groups:result100}. Let 
$\Om_s\subset L$ be the orbit under $W(\Psi)$ of $u_s$ and
$\Lambda_s\subset L$ the sublattice spanned by $\Om_s$. 
Thus $U_s=\Lambda_s\otimes\R$ is a $W(\Psi)$-invariant subspace
of $V$, and as the reflectional representation is irreducible we have 
$U_s=V$, and so $\Lambda_s$ is a free $\Z$-module, of the same rank $n$ as $L$.

Consider $L/2L:=L/2$, an $n$-dimensional vector space over $\F_2$,
the subspace $\Lambda_s/2\subset L/2$,
and write $\ov{v}$ for the image of $v\in L$ via the quotient map
$L\rightarrow L/2$. When making statements about $\ov{u}_s$ in the remainder of
the paper, it is important to use the explicit $u_s$ given in the proof
of Lemma \ref{subsection:weyl_groups:result100} (particularly with factors 
of $2$ removed!).
As $2L$ is a $W(\Psi)$-invariant submodule, 
the $W(\Psi)$-actions on $L$ and $\Lambda_s$ induce
actions on $L/2$ and $\Lambda_s/2$.

Let $s,t\in\Psi$, and since the Weyl groups all have symbols that
are trees, there is a unique path containing a minimal number of edges
connecting $s$ to $t$, with nodes $s=s_1,\ldots,s_k=t$ say. Let,
$$
X^s_t=\{\ov{u}_s,s_1(\ov{u}_s),s_2s_1(\ov{u}_s),\ldots,s_k\cdots s_2s_1(\ov{u}_s)\},
$$
a collection of (at most) $k+1$ vectors in $\Lambda_s/2$. 
Observe that for $t'\in\Psi$, the sets $X^s_t$ and $X^s_{t'}$ intersect in $\ell+1$ vectors
when the paths connecting $s$ to $t$ and $t'$ have in common an initial
path containing $\ell$ edges, so in particular, $X_t^s\cap X_{t'}^s$ always 
contains at least two vectors.

%For a fixed $s$,
%we are interested in the largest number of $t$ we can take so that the $X^s_t$'s
%give linearly independent vectors. Thus, 
For $s\in\Psi$, call a set of nodes $T\subset\Psi$
\emph{independent for $s$\/} when $\bigcup_{t\in T}X^s_t$ contains
$m+1$ linearly independent vectors in $\Lambda_s/2$, where
$m$ is the number of nodes contained in the union of all the minimal paths
connecting $s$ to the $t\in T$.

Finally, the \emph{independence data\/} for $W(\Psi)$ 
is a labeling of the nodes of $\Psi$ by strings of the form,
$$
t_1,\ldots,t_k,\ov{t_{11},\ldots,t_{1n_1}},\ldots,
\ov{t_{\ell 1},\ldots,t_{\ell n_\ell}},
\text{ with the $t$'s in $\Psi$},
$$
such that if $T$ is a union of any subset of the $t_1,\ldots,t_k$
together with the $t$'s in exactly \emph{one\/} of the $\ov{t_{i1},\ldots,t_{i n_i}}$,
then $T$ is independent for the node $s$ having this label.
As always, the node labels follow the conventions in Table \ref{table:roots1}, 
and if a node $s$ is labeled by $0$, this indicates that $\{\ov{u}_s\}$ is
non-zero, but the vectors in $X^s_t$ are not independent for any $t$
(including $t=s$).

For example, in Proposition 
\ref{subsection:weyl_groups:result150} an $s$ that is $\ell\equiv 2\text{ mod }4$ nodes from the 
righthand end
of $\Psi=D_n$ has label
$2,\ov{1,n},\ov{1,n-1},\ov{n,n-1}$,
indicating that, say
\begin{align*}
X^s_2\cup X^s_1\cup X^s_{n-1}=X^s_1\cup X^s_{n-1}
=\{&\ov{u}_s,s_\ell(\ov{u}_s),s_{\ell-1}s_\ell(\ov{u}_s),\ldots,
s_1\cdots s_\ell(\ov{u}_s),\\
&s_{\ell+1}s_\ell(\ov{u}_s)
\ldots,s_{n-1}\cdots s_\ell(\ov{u}_s)\},
\end{align*}
consists of $n$ independent vectors in $\Lambda_s/2$.
The barred $t$'s are meant to 
stop us from getting too carried away: $X^s_{1}\cup X^s_{n-1}
\cup X^s_{n}$ for instance contains too many (ie: more than $n$) distinct vectors
to be independent.

\begin{proposition}%[type $\mathbf{A}$ and $\mathbf{D}$]
\label{subsection:weyl_groups:result125}
The independence data for $\Psi=A_n$ is,
$$
\begin{pspicture}(0,0)(15,1.5)
%\showgrid
\rput(1,-.5){
\rput(-0.8,0){
\rput(3.6,1.25){\BoxedEPSF{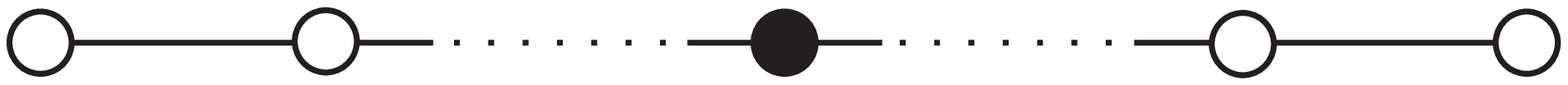 scaled 300}}
\rput(2.2,1.95){$\ell$}\rput(5,.55){$k$}
\rput(-3.7,.75){\rput{90}(5.9,.8){$\left.\begin{array}{c}
\vrule width 0 mm height 28 mm depth 0 pt\end{array}\right\}$}}
\rput(-.9,.15){\rput{90}(5.9,.8){$\left\{\begin{array}{c}
\vrule width 0 mm height 28 mm depth 0 pt\end{array}\right.$}}
}
\rput(11,1.25){$\left\{\begin{array}{l}0,\,\text{for }\ell/d,k/d\text{ odd}\\
1,n-1,\,\text{for }\ell/d\text{ odd},k/d\text{ even}\\
2,n,\,\text{for }\ell/d\text{ even},k/d\text{ odd}\end{array}\right.$}
\rput(7.8,1.25){\BoxedEPSF{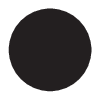 scaled 300} $=$}
}
\end{pspicture}
$$
where $d=\text{gcd}(\ell,k)$.
\end{proposition}

\begin{proof}
%Let $T=\{s_1,\ldots,s_k\}$ with $s_i$ labelling the $i$-th node from the left in 
%a type $A_k$ subsymbol, and suppose that $s=s_k$. As $u_s$ is orthogonal to
%the $v_{s'}$ for $s'\in\{s_1,\ldots,s_{k-1}\}$, 
%the $W_T$-orbit of $u_s$ %can be obtained by taking the images of $u_s$ 
%%by representatives for the cosets in 
%consists of the images of $u_s$ under representatives
%for the cosets $W_T/\langle s_1,\ldots,s_{k-1}\rangle$.
%In particular, the pair is admissible precisely when the $L\rightarrow L/2$
%images of
%$$
%u_s,s_k(u_s),s_{k-1}s_ks_{k-1}(u_s),\ldots,s_1\ldots s_{k-1}s_ks_{k-1}\ldots s_1(u_s),
%$$
%are $\Z/2$-independent. Equivalently one may take
%just $u_s,s_k(u_s),s_{k-1}s_k(u_s),\ldots,s_1\ldots s_k(u_s)$. It now reduces to 
%a case by case check: rather than try the readers patience, we reproduce the details
%just in type $A$; the others are similar (and easier). %To this end, we have the
%vector $u_s$ given by
%Suppose $(\Psi,s)$ is,
For the node $s$ indicated the vector $u_s\in L$ is given by,
$$
\begin{pspicture}(0,0)(14,1)
%\showgrid
\rput(0,0){
\rput(7,.5){\BoxedEPSF{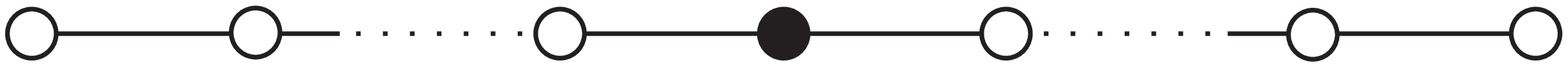 scaled 300}}
\rput(3.45,.9){${\scriptstyle k/d}$}
\rput(4.5,.9){${\scriptstyle 2k/d}$}
\rput(5.6,.9){${\scriptstyle (\ell-1)k/d}$}
\rput(7,.9){${\scriptstyle k\ell/d}$}
\rput(8.2,.1){${\scriptstyle (k-1)\ell/d}$}
\rput(9.5,.1){${\scriptstyle 2\ell/d}$}
\rput(10.55,.1){${\scriptstyle \ell/d}$}
}
\end{pspicture}
$$ 
where at least one of $\ell/d$ and $k/d$ must be odd. 
Thus $\ov{u}_s=10\cdots 101\cdots 01$ or $10\cdots 010\cdots 01$
in the first case listed at the right in the Lemma,
and $00\cdots 001\cdots 01$
or $10\cdots 100\cdots 00$ in the other two.
It's a straight out calculation from here on in.
%with the labels giving the vector $u_s$ and
%$k+\ell=n+1$. As $n+1$ is odd we have $k$, say,
%is even and $\ell$ odd. In particular, no $T$ spanning a subsymbol
%of type $A$, with $|T|$ odd, 
%and $s$ at the end, contains the node labelled $\ell$. The $L\rightarrow L/2$
%image of $u_s$ is $\ov{u}_s=0\cdots 00101\cdots 01$ (with $\ell$ $0$'s at the
%beginning), and $s_i\ldots s_k(u_s)$ has image
%$\ov{u}_s+v_i$, where $v_i$ has $1$'s in the positions corresponding to
%the $s_i,\ldots,s_k$, and $0$'s elsewhere. These are clearly $\Z/2$-independent.
\qed
\end{proof}

\begin{proposition}%[type $\mathbf{B}$]
\label{subsection:weyl_groups:result150}
The independence data for $\Psi=B_n$ is,
$$
\begin{pspicture}(0,0)(15,1.5)
%\showgrid
\rput(1,-.5){
\rput(-0.4,0){
\rput(3.6,1.25){\BoxedEPSF{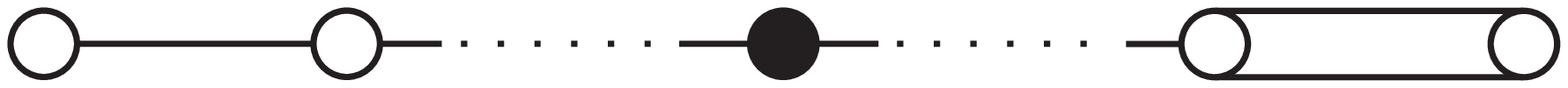 scaled 300}}
\rput(0,-.3){
\rput(-3.6,1){\rput{90}(5.9,.8){$\left.\begin{array}{c}
\vrule width 0 mm height 24 mm depth 0 pt\end{array}\right\}$}}
\rput(2.3,2.2){$\ell$}
}
%\rput(.6,0){
%\rput(-.5,.4){\rput{90}(5.9,.8){$\left\{\begin{array}{c}
%\vrule width 0 mm height 35 mm depth 0 pt\end{array}\right.$}}
%\rput(5.4,.8){$\ell$}
%}
}
\rput(11,1.25){$\left\{\begin{array}{l}0,\,\text{for }\ell\text{ odd}\\
2,n,\,\text{for }\ell\equiv 0\text{ mod }4\\
\ov{1},2,n-1,\ov{n},\,\text{for }\ell\equiv 2\text{ mod }4\end{array}\right.$}
\rput(7.9,1.25){\BoxedEPSF{admissible_vertex.eps scaled 300} $=$}
}
\end{pspicture}
$$
and the labeling of the right-most
node is $0$ for $n$ odd, and $n$ for $n$ even.
\end{proposition}

%\begin{proof}
%\qed
%\end{proof}

The proof is completely analogous to Proposition 
\ref{subsection:weyl_groups:result125},
as is the proof of:

\begin{proposition}
\label{subsection:weyl_groups:result150}
The independence data for $\Psi=D_n$ is,
$$
\begin{pspicture}(0,0)(15,1.5)
%\showgrid
\rput(0,-.3){
\rput(-0.5,0){
\rput(3.6,1){\BoxedEPSF{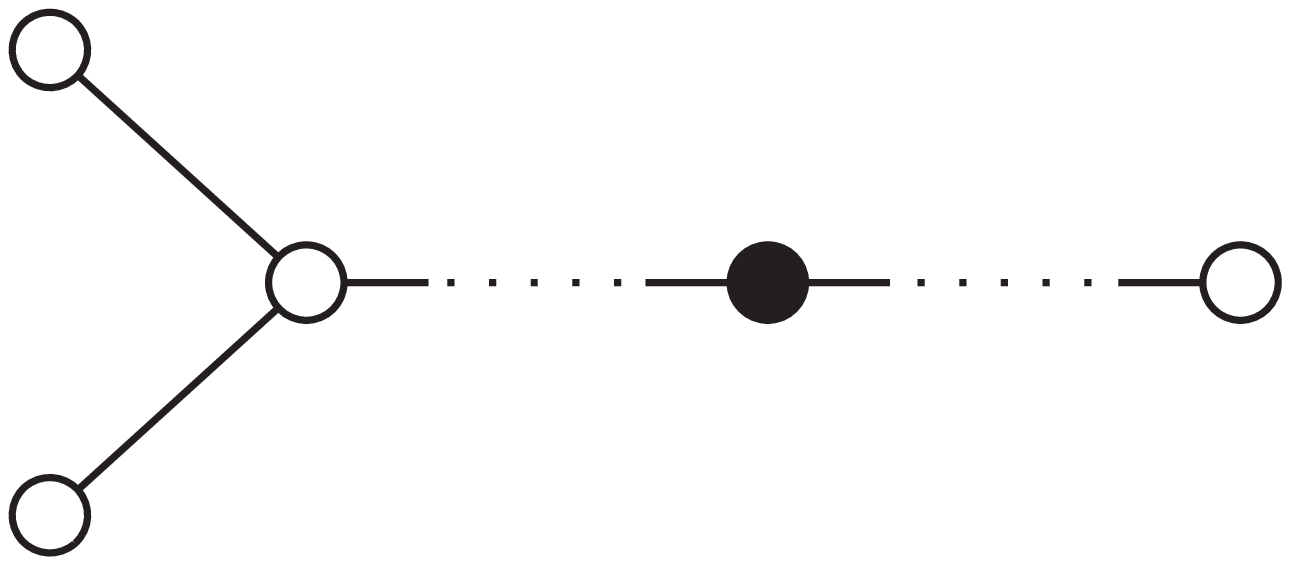 scaled 300}}
\rput(-1.2,.5){\rput{90}(5.9,.8){$\left.\begin{array}{c}
\vrule width 0 mm height 15 mm depth 0 pt\end{array}\right\}$}}
\rput(4.7,1.7){$\ell$}\rput(1.5,1.7){$0$}\rput(1.5,.3){$0$}
}
\rput(7.2,1.05){\BoxedEPSF{admissible_vertex.eps scaled 300} $=$}
\rput(11.2,1.05){$\left\{\begin{array}{l}0,\,\text{for }\ell\text{ odd}\\
2,n-1,n,\,\text{for }\ell\equiv 0\text{ mod }4\\
2,\ov{1,n},\ov{1,n-1},\ov{n,n-1},\,\text{for }
\ell\equiv 2\text{ mod }4\end{array}\right.$}
%\psline[linewidth=.2mm]{->}(7.7,1.5)(7,1.5)(7,.8)(4.1,.8)(4.1,1.2)
}
\end{pspicture}
$$
\end{proposition}

Finally, we have

\begin{proposition}%[type $\mathbf{E},\mathbf{F}$ and $\mathbf{G}$]
\label{subsection:weyl_groups:result175}
The independence data for $\Psi$ exceptional is,
$$
\begin{pspicture}(0,0)(15,3.5)
%\showgrid
\rput(0,0){
\rput(2.5,2.5){\BoxedEPSF{F4_dynkin.eps scaled 300}}
\rput(.9,2.9){${\scriptstyle 3}$}\rput(2,2.9){${\scriptstyle 4}$}
\rput(3,2.9){${\scriptstyle 3}$}\rput(4.05,2.9){${\scriptstyle 4}$}
\rput(2.5,3.3){$F_4$}
}
\rput(-5,-.4){
\rput(12,3){\BoxedEPSF{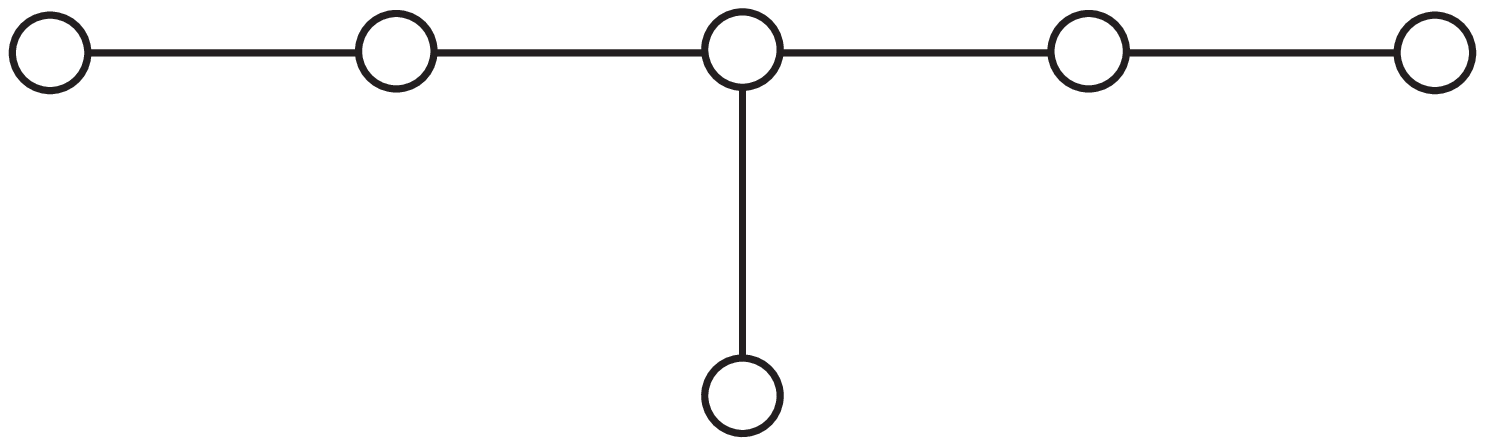 scaled 300}}
\rput(8.8,2.1){
\rput(1.1,1.8){${\scriptstyle 5}$}\rput(2.2,1.76){${\scriptstyle 5,6}$}
\rput(3.25,1.76){${\scriptstyle 1,5}$}
\rput(4.3,1.76){${\scriptstyle 1,6}$}\rput(5.3,1.8){${\scriptstyle 1}$}
\rput(2.4,.4){${\scriptstyle \ov{1},2,4,\ov{5}}$}}
\rput(13,2.8){$E_6$}
}
\rput(0,0){
\rput(12.5,3){\BoxedEPSF{G2_dynkin.eps scaled 300}}
\rput(12,2.6){${\scriptstyle 1}$}\rput(13,2.6){${\scriptstyle 2}$}
\rput(11.2,3){$G_2$}
}
\rput(0,-.4){
\rput(3.25,1){\BoxedEPSF{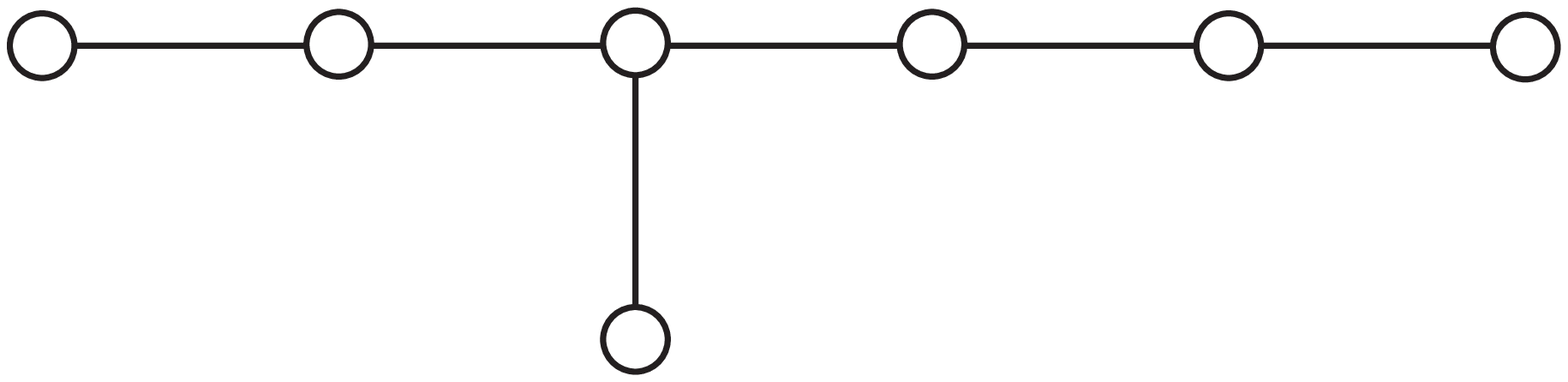 scaled 300}}
\rput(-.4,0.1){
\rput(1.05,1.76){${\scriptstyle 5,7}$}\rput(2.1,1.76){${\scriptstyle 6,7}$}
\rput(3.1,1.76){${\scriptstyle 1,5,7}$}
\rput(4.15,1.8){${\scriptstyle 0}$}\rput(5.25,1.76){${\scriptstyle 1,6}$}
\rput(6.3,1.8){${\scriptstyle 0}$}
\rput(2.6,.4){${\scriptstyle 1,5}$}
\rput(4,.8){$E_7$}
}
}
\rput(0.5,-.4){
\rput(10.25,.975){\BoxedEPSF{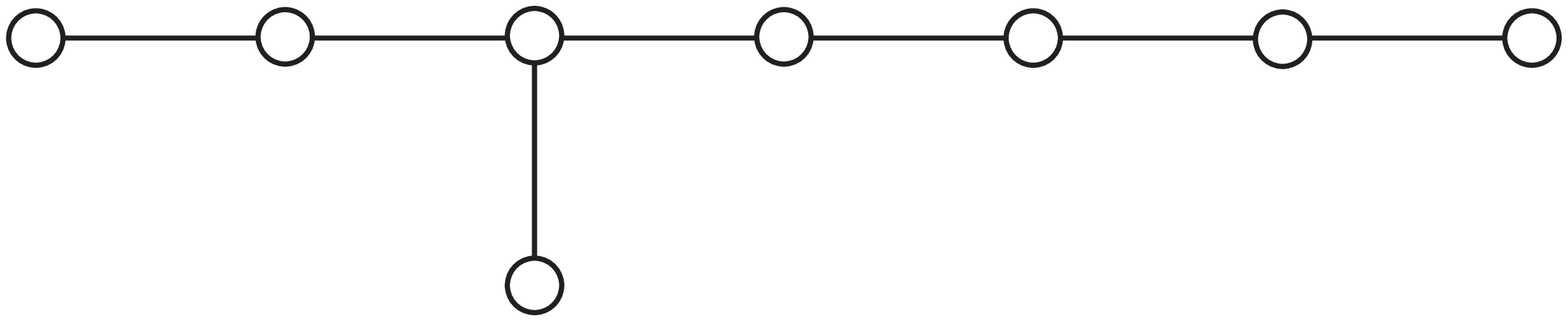 scaled 300}}
\rput(6,0.09){
\rput(1.05,1.8){${\scriptstyle 7}$}
\rput(2.15,1.76){${\scriptstyle 7,8}$}\rput(3.2,1.76){${\scriptstyle 1,7}$}
\rput(4.25,1.76){${\scriptstyle 1,6,8}$}\rput(5.3,1.76){${\scriptstyle 1,7}$}
\rput(6.35,1.8){${\scriptstyle 8}$}
\rput(7.4,1.8){${\scriptstyle 1}$}
\rput(2.4,.4){${\scriptstyle \ov{1},2,6,\ov{7}}$}}
\rput(11,.8){$E_8$}
}
\end{pspicture}
$$
\end{proposition}

%For $s\in\Psi$ a node, call the pair $(\Psi,s)$ \emph{admissible\/} whenever
%we have the following: for each $T\subset S$ such that $T$ spans a subsymbol
%$A_k\subset\Psi$ for $k=|T|$ odd, and with $s$ one of the two end vertices,
%the image under $L\rightarrow L/2$ of the $W_T$-orbit of $u_s$ spans 
%a free $\Z/2$-submodule of rank $|T|+1$, where $u_s$ is the vector
%in $L$ given by Lemma \ref{subsection:weyl_groups:result100}.
%The pair is \emph{specially admissible\/} when we have this condition
%%without the $k$ odd condition, so that 
%for all $T$ spanning a $A_k\subset\Psi$
%with no restrictions now on $k$.

Now to another important concept.
Let $W(\Psi)$ be an irreducible Weyl group and $s\in\Psi$. The pair
$(\Psi,s)$ is said to be \emph{specially admissible\/} if, 
\begin{description}
\item[(i).] in the root lattice $L$ we have $x_s=v_s$, and,
\item[(ii).] for every $t\in\Psi$
such that the minimal path from $s$ to $t$ spans a subsymbol of type $A$, we have
the set $\{t\}$ is independent for $s$. 
\end{description}
The pair is \emph{admissible\/}
if we have (i), and for every $t\in\Psi$
such that the path from $s$ to $t$ spans a subsymbol of type $A$ of \emph{odd\/}
rank, then $\{t\}$ is independent for $s$.

The first part of the definition is just an artificial device
designed to exclude the pairs $(B_n,s_n)$, $(F_4,s_3)$, $(F_4,s_4)$
and $(G_2,s_2)$, for reasons that will be explained in the proof of
Theorem \ref{section:lattices:result300} below. 
It is the second part that will prove natural and useful 
in \S\ref{section:homomorphisms},
essentially because of the following observation:
%The following observation will make these concepts
%useful later: 
suppose that 
$A_k\subset\Psi$ is the subsymbol spanned by such a path
from $s$ to $t$, with nodes $s=s_1,\ldots,s_k=t$. Then $u_s$ is 
orthogonal to $x_2,\ldots,x_k$, so is fixed by the reflections $s_2,\ldots,s_k$.
As $\{1,s_1,s_2s_1,\ldots,s_k\cdots s_1\}$ is a set of coset representatives
for $\langle s_2,\ldots,s_k\rangle$ in $W(A_k)$, we get that
$X^s_t$ constitutes the whole $W(A_k)$-orbit of $\ov{u}_s$ in $L/2$.
Thus $T=\{t\}$ is independent for $s$ is equivalent to this $W(A_k)$-orbit
spanning a subspace of dimension $k+1$.

\begin{theorem}\label{section:lattices:result200}
%Let $\Psi$ %simply-laced, connected, crystallographic, 
%be the symbol of an irreducible Weyl group. 
Let $W(\Psi)$ be a classical irreducible Weyl group. 
Then the admissible
(respectively specially admissible) pairs $(\Psi,s)$ are given by the following, where
$s=\BoxedEPSF{admissible_vertex.eps scaled 300}$ 
(resp. $s=\BoxedEPSF{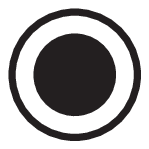 scaled 250}$\/),
$$
\begin{pspicture}(0,0)(15,3.2)
%\showgrid
\rput(0,-.65){
\rput(1.6,0.2){
\rput(2,3){\BoxedEPSF{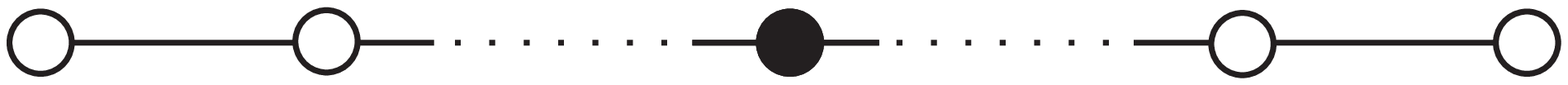 scaled 300}}
\rput(3.5,3.5){$A_n$}
\rput(-5.25,2.5){\rput{90}(5.9,.8){$\left.\begin{array}{c}
\vrule width 0 mm height 27 mm depth 0 pt\end{array}\right\}$}}
\rput(-2.5,1.9){\rput{90}(5.9,.8){$\left\{\begin{array}{c}
\vrule width 0 mm height 27 mm depth 0 pt\end{array}\right.$}}
\rput(0.7,3.7){$\ell$}\rput(3.4,2.3){$k$}
}
\rput(4.6,0.2){
\rput(7,3){\BoxedEPSF{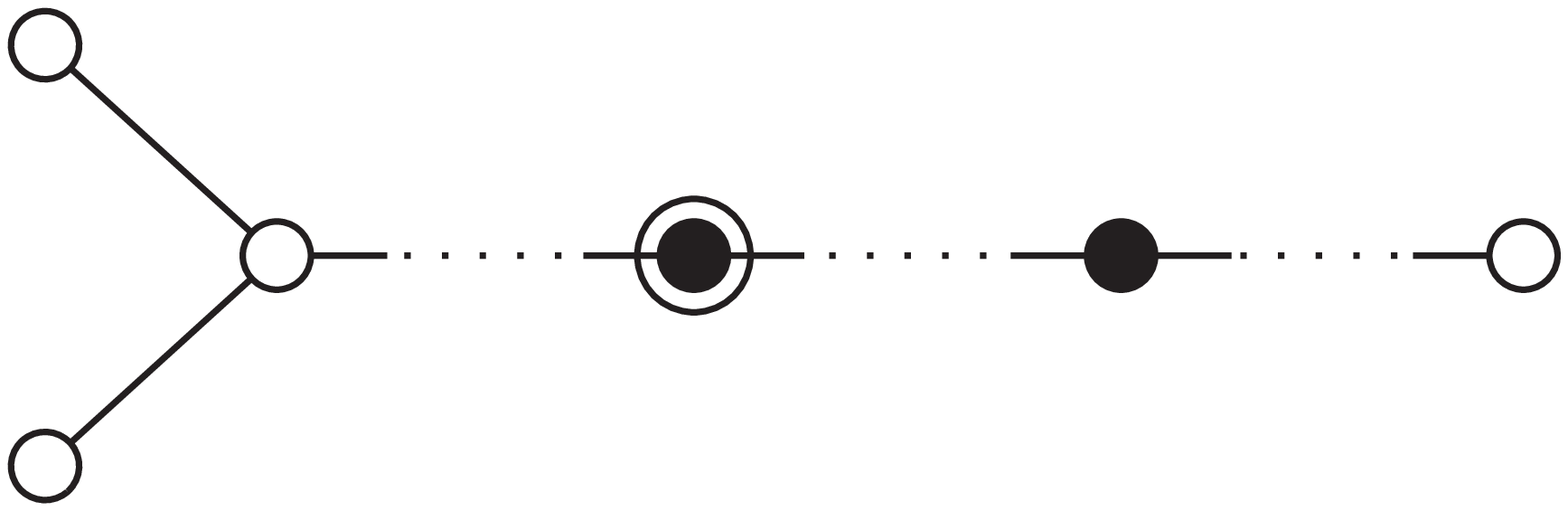 scaled 300}}
\rput(4.5,3){$D_n$}
\rput(2.15,2.5){\rput{90}(5.9,.8){$\left.\begin{array}{c}
\vrule width 0 mm height 30 mm depth 0 pt\end{array}\right\}$}}
\rput(2.9,1.9){\rput{90}(5.9,.8){$\left\{\begin{array}{c}
\vrule width 0 mm height 15 mm depth 0 pt\end{array}\right.$}}
\rput(8.05,3.7){$k$}\rput(8.8,2.3){$\ell$}
}
\rput(0,.20){
\rput(7.5,1){\BoxedEPSF{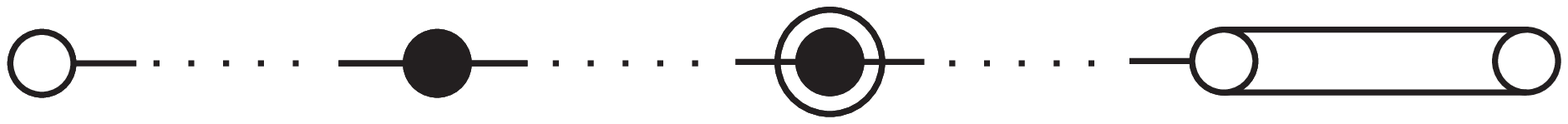 scaled 300}}
\rput(4,1){$B_n$}
\rput(-1.8,-2){
\rput(2.15,2.5){\rput{90}(5.9,.8){$\left.\begin{array}{c}
\vrule width 0 mm height 30 mm depth 0 pt\end{array}\right\}$}}
\rput(1.4,1.9){\rput{90}(5.9,.8){$\left\{\begin{array}{c}
\vrule width 0 mm height 15 mm depth 0 pt\end{array}\right.$}}
\rput(8.05,3.7){$k$}\rput(7.3,2.3){$\ell$}}
}
}
\end{pspicture}
$$
and in type $A$ we have $\ell=2^a m_1,k=2^b m_2$ with the $m_i$ odd,
and $a<b$; in types $B$ and $D$ we have $\ell\equiv 0\text{ mod }4$
and $k\equiv 2\text{ mod }4$.
\end{theorem}

The conditions in type $A$ turn out to be quite mild: if $n$ is even then
all $(A_n,s)$ are admissible; if $n$ is odd then $k$ and $\ell$ have the same
parity, and the condition becomes that they are both even, but with different
$2$-parts. The proof can be read straight off the diagrams in Propositions
\ref{subsection:weyl_groups:result125}-\ref{subsection:weyl_groups:result175},
as can the proof of,

\begin{theorem}\label{section:lattices:result250}
%Let $\Psi$ %simply-laced, connected, crystallographic, 
%be the symbol of an irreducible Weyl group. 
Let $W(\Psi)$ be an exceptional irreducible Weyl group. 
Then the admissible
(resp. specially admissible) pairs $(\Psi,s)$ are given by the following, where
$s=\BoxedEPSF{admissible_vertex.eps scaled 300}$ 
(resp $s=\BoxedEPSF{special_admissible_vertex.eps scaled 250}$\/),
$$
\begin{pspicture}(0,0)(15,3)
%\showgrid
\rput(0,3){
\rput(0,-5.5){%moves G
\rput(1.85,5){\BoxedEPSF{G23_dynkin.eps scaled 300}}
\rput(.5,5){$G_2$}
%\rput(1.85,5.2){$6$}
}
\rput(0,0){
\rput(6,-.5){\BoxedEPSF{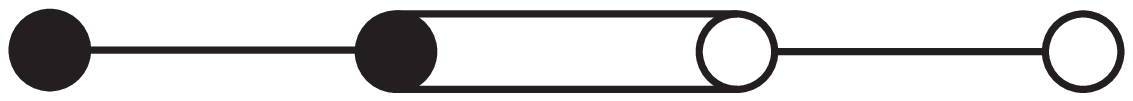 scaled 300}}
\rput(8.2,-.5){$F_4$}
}
\rput(-.5,-3.9){
\rput(12,3){\BoxedEPSF{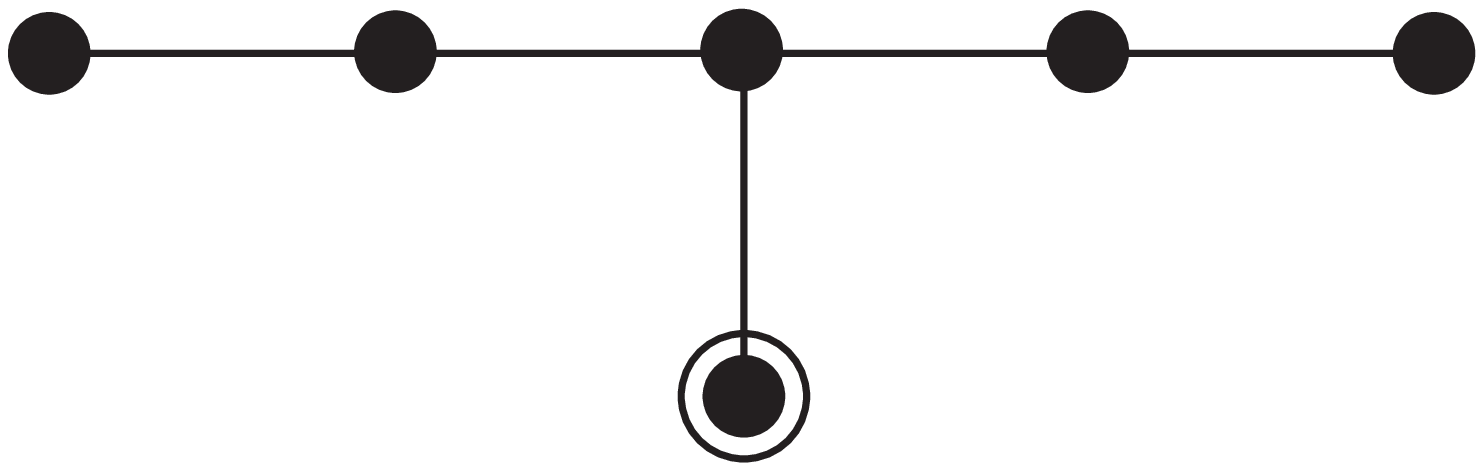 scaled 300}}
\rput(13,2.8){$E_6$}
}
\rput(0,-3.5){
\rput(3.25,1){\BoxedEPSF{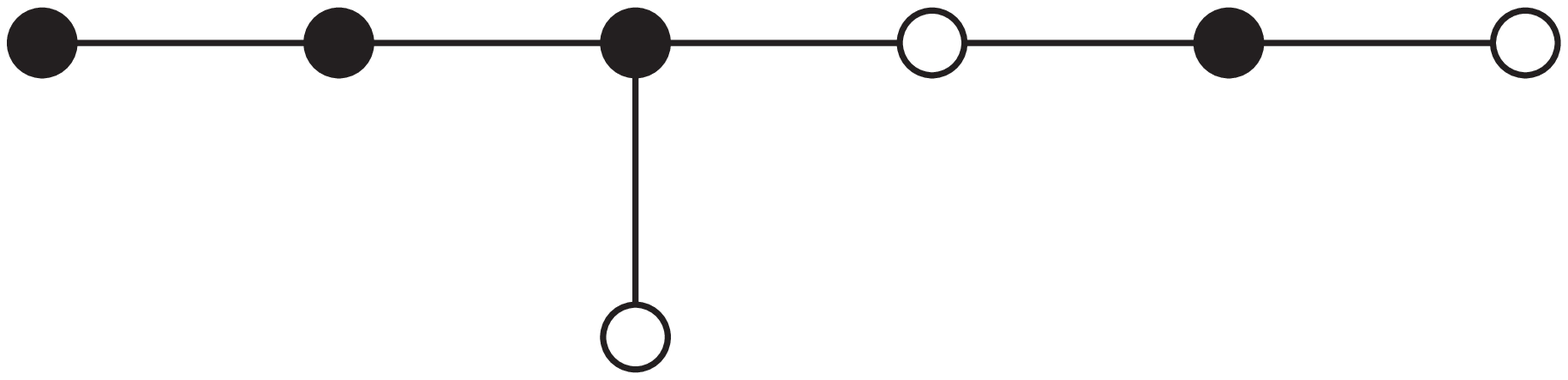 scaled 300}}
\rput(4,.8){$E_7$}
}
\rput(0,-3.5){
\rput(10.25,.975){\BoxedEPSF{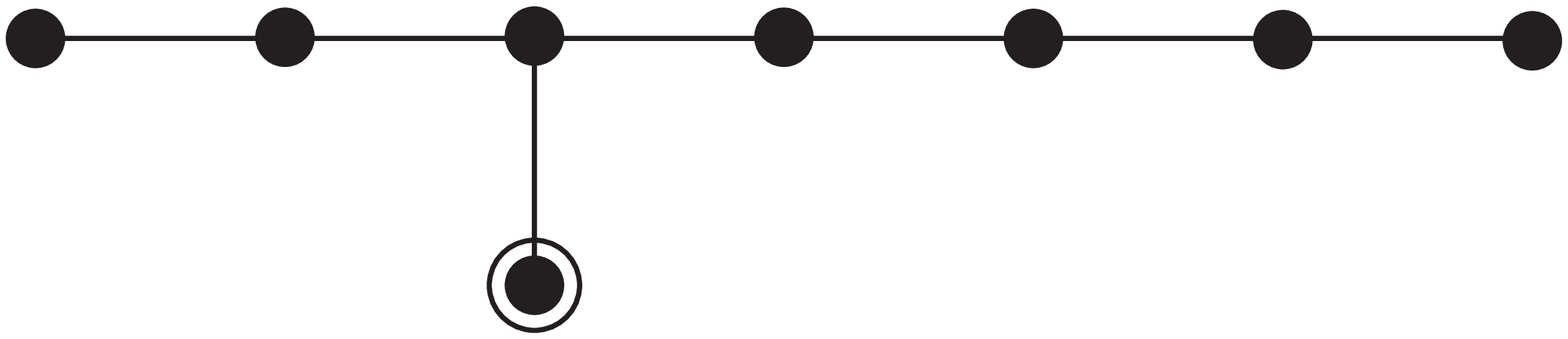 scaled 300}}
\rput(11,.8){$E_8$}
}
}
\end{pspicture}
$$
\end{theorem}

Another spin-off resulting from the independence data is the following,

\begin{theorem}\label{section:lattices:result300}
%Let $\Psi$ %simply-laced, connected, crystallographic, 
%be the symbol of an irreducible Weyl group. 
Let $W(\Psi)$ be an irreducible Weyl group of rank $n$ and  $(\Psi,s)$ admissible.
Then the subspace
$\Lambda_s/2\subset L/2$ is $n$-dimensional over $\F_2$.
\end{theorem}

\begin{proof}
The result follows as for any admissible pair $(\Psi,s)$ one can find
a $T\subset\Psi$, with $T$ independent for $s$, 
and such that the union of the nodes in the minimal paths connecting $s$ to
the $t\in T$ contains all but one of the
nodes of $\Psi$. The pairs $(B_n,s_n)$, $(F_4,s_3)$, $(F_4,s_4)$
and $(G_2,s_2)$, while %\emph{morally\/} admissible in that they 
rather trivially satisfying
condition (ii) in the definition of admissability, do not have
the property of the previous sentence, hence their exclusion.
\qed
\end{proof}

%If we consider 
As an example,
the pair $(\Psi,s)=\begin{pspicture}(0,0)(2.5,.4)
%\showgrid
\rput(1.25,.1){\BoxedEPSF{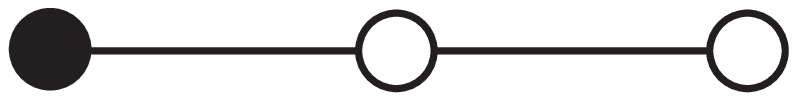 scaled 300}}
\end{pspicture}$ is not admissible by Theorem \ref{section:lattices:result200},
and for $t$ the rightmost node, $X_s^t$ is the $W(\Psi)$-orbit of $\ov{u}_s$,
and $\Lambda_s/2$ is a $1$-dimensional subspace of $L/2$. % Thus, the admissibility
% is essential to Theorem \ref{section:lattices:result300}.

From now on we will only be interested in the sublattices $\Lambda_s\subset L$
where $(\Psi,s)$ is admissible, so
that the vector spaces $\Lambda_s/2\subset L/2$ coincide by Theorem 
\ref{section:lattices:result300}.
In the rest of this section we look at some aspects of the action of the Weyl group
$W(\Psi)$ on $L/2$ as preparation for \S\ref{section:extended}.

To this end, let $\endo(L/2)$ be the endomorphism algebra of the space,
writing $g(\ov{v})$ for the image of $\ov{v}\in L/2$ under $g\in W(\Psi)$.
If $g$ is an involution then $(g+1)^2=0$ in $\endo(L/2)$, hence
$\im(g+1)\subset\ker(g+1)$;
moreover, $\ker(hgh^{-1}+1)=h\,\ker(g+1)$ and $\im(hgh^{-1}+1)=h\,\im(g+1)$.
Thus the dimension of the quotient space
$\ker(g+1)/\im(g+1)$ is an invariant
of the conjugacy class of the involution $g$, and we have proved the first half
of,

\begin{theorem}\label{section:lattices:result400}
Let $W(\Psi)$ be an irreducible Weyl group %$\Psi\not=A_n$ ($n$ even),
with even Coxeter number $h$ and $\xi$ a Coxeter element. Then the dimension 
$d_\Psi$ of
the $\F_2$ vector space 
$$
\text{ker}(\xi^{h/2}+1)/\text{im}(\xi^{h/2}+1),
$$
is independent of the choice of $\xi$, with these dimensions being,
$$
\begin{tabular}{cccccccccc}
\hline
$\Psi$&$A_n\,(n\text{ odd})$&$B_n$&$D_n\,\,(n\text{ even})$&$D_n\,(n\text{ odd})$&
$G_2$&$F_4$&$E_6$&$E_7$&$E_8$\\
$d_\Psi$&$1$&$n$&$n$&$n-2$&
$2$&$4$&$2$&$7$&$8$\\
\hline
\end{tabular}
$$
\end{theorem}

Note that the Theorem applies to all the Weyl groups except
type $A_n$ for $n$ even. Determining these dimensions
will take a little more brute force effort; indeed we will explicitly describe
the image and kernel for a particular Coxeter element, as these descriptions
will prove useful in Proposition \ref{section:lattices:result500}
below.

\begin{proof}
First and easiest, if $\Psi$ has $(-1)$-type, then 
by Proposition \ref{section:preliminaries:torsion:result400}, 
$\xi^{h/2}=w_\Psi$ for any $\xi$ (as
$w_\Psi$ is central), %the element of longest length in $W(\Psi)$,
and $w_\Psi+1$ is the zero map in $\endo(L/2)$. Hence
$\{0\}=\im(\xi^{h/2}+1)\subset\ker\xi^{h/2}+1=L/2$ and $d_\Psi=n$.

This leaves $A_n$ ($n$ odd), $D_n$ ($n$ odd) and $E_6$, for which
we will consider the Coxeter element $\xi=s_1\ldots s_n$, numbering, as always,
from Table \ref{table:roots1}.
For type $A$, let $\{u_1,\ldots,u_{n+1}\}$ be the orthonormal
basis realizing the isomorphism $W(A_n)\rightarrow\SS_{n+1}$
from \S\ref{subsection:weyl_groups}. Then $\xi$ corresponds to
the $n+1$-cycle $(u_1,\ldots,u_{n+1})$ 
and $\xi^{h/2}$ to $(u_1,u_{\ell+1})(u_2,u_{\ell+2})
\cdots(u_\ell,u_{n+1})$ for $\ell=(n+1)/2$. Thus with $\{x_i\}$ the basis for
the root lattice, $\xi^{h/2}$ interchanges $x_i$ and $x_{\ell+i}$
($i\not=\ell$) and sends $x_\ell$ to $-x_\ell$. Thus,
$\im(\xi^{h/2}+1)=\langle\ov{x}_i+\ov{x}_{\ell+i}\,|\,1\leq i<\ell\rangle\subset L/2$,
a subspace of dimension $\ell-1$, which %, by the first isomorphism theorem for
%vector spaces, 
gives a kernel of dimension $\ell$ and hence $d_\Psi=1$ as
claimed. It will be useful though
to have an explicit description: 
$\ker(\xi^{h/2}+1)=\im(\xi^{h/2}+1)\oplus\langle\ov{x}_\ell\rangle$.

\parshape=5 0pt\hsize 0pt\hsize 
0pt.55\hsize 0pt.55\hsize 0pt.55\hsize
For type $D$ and its realisation
$W(D_n)\rightarrow\SS^\pm_\circ(u_1,\ldots,u_n)$ as the group of even signed permutations
of $\{u_1,\ldots,u_n\}$, the Coxeter element $\xi$ corresponds to the
even signed permutation at right.
In particular, when $n$ is odd,
$\im(\xi^{h/2}+1)$ is the $1$-dimensional space spanned by
$\ov{x}_{n-1}+\ov{x}_n$ and 
$\ker(\xi^{h/2}+1)=\im(\xi^{h/2}+1)\oplus\langle\ov{x}_1,\ldots,\ov{x}_{n-2}\rangle$,
giving the required $d_\Psi=n-2$.
\vadjust{\hfill\smash{\lower 22pt
\llap{
\begin{pspicture}(0,0)(6,2.5)
%\showgrid
\rput(3,1.25){\BoxedEPSF{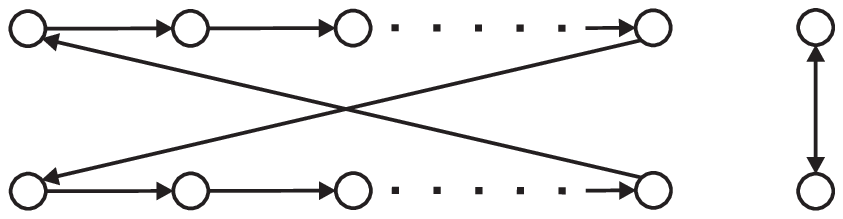 scaled 700}}
%\rput(-.75,1.25){$\zeta=$}
\rput(0,0){\rput(.2,2.15){${\scriptstyle u_1}$}\rput(.075,.35){${\scriptstyle -u_1}$}}
\rput(1.2,0){\rput(.2,2.15){${\scriptstyle u_2}$}\rput(.075,.35){${\scriptstyle -u_2}$}}
\rput(2.3,0){\rput(.2,2.15){${\scriptstyle u_3}$}\rput(.075,.35){${\scriptstyle -u_3}$}}
\rput(4.4,0){\rput(.2,2.15){${\scriptstyle u_{n-1}}$}
\rput(.075,.35){${\scriptstyle -u_{n-1}}$}}
\rput(5.7,0){\rput(.2,2.15){${\scriptstyle u_n}$}\rput(.075,.35){${\scriptstyle -u_n}$}}
\end{pspicture}
}}}\ignorespaces

\parshape=3 0pt.55\hsize 0pt.55\hsize 0pt\hsize
Finally, $E_6$ \emph{really\/} can be done by brute force, ie: by writing the generating
reflections $s_i$ as matrices in terms of the basis $\{x_i\}$ for the root
lattice and proceeding from there to get 
$\im(\xi^{h/2}+1)=\langle\ov{x}_1+\ov{x}_4,\ov{x}_2+\ov{x}_5\rangle$
and $\ker(\xi^{h/2}+1)=\im(\xi^{h/2}+1)\oplus
\langle\ov{x}_1+\ov{x}_2+\ov{x}_3,\ov{x}_6\rangle$.
\qed
\end{proof}

We will use
these descriptions in \S\ref{section:extended}, where not the quotient,
but the \emph{difference\/} $\ker(\xi^{h/2}+1)\setminus\im(\xi^{h/2}+1)\subset L/2$
will be a target which we want to ``hit'' with a certain endomorphism. As 
is so often the case in 
such situations, this will be possible if the target is big enough, in fact
as it transpires, if the dimension $d_\Psi>1$.

\begin{proposition}\label{section:lattices:result500}
Let $W(\Psi)$ be an irreducible Weyl group of rank $n$ with
the dimension $d_\Psi>1$ and 
Coxeter number $h=2^pq$ where $p>0$ and $q$ odd. Then 
for the Coxeter element $\xi=s_1\ldots s_n$, and
$$
\aa=1+\xi^q+\xi^{2q}+\cdots+\xi^{(2^{p-1}-1)q}\in\endo(L/2),
$$
there is a $\ov{u}\in L/2$ with 
$\aa(\ov{u})\in\ker(\xi^{h/2}+1)\setminus\im(\xi^{h/2}+1)$.
\end{proposition}

\begin{proof}
\parshape=5 0pt\hsize 0pt\hsize 
0pt.55\hsize 0pt.55\hsize 0pt.55\hsize
We proceed again on a case by case basis, starting with the $\Psi$ of
$(-1)$-type, for which the difference consists of all the non-zero vectors
in $L/2$ as $d_\Psi=n$. Thus it suffices to show that $\aa\not=0$ in $\endo(L/2)$, and then
we may choose any $\ov{u}\not\in\ker\aa$. Now if $p=1$, then $\aa=1$, and we are
done, leaving $\Psi=B_n$ and $F_4$ in the $(-1)$-type case to do.
\vadjust{\hfill\smash{\lower 13mm
\llap{
\begin{pspicture}(0,0)(6,2.5)
%\showgrid
\rput(3,1.25){\BoxedEPSF{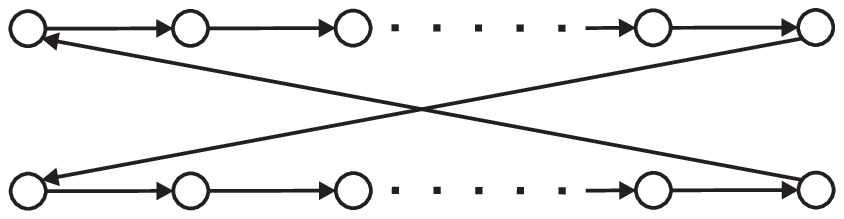 scaled 700}}
%\rput(-.75,1.25){$\zeta=$}
\rput(0,0){\rput(.2,2.15){${\scriptstyle u_1}$}\rput(.075,.35){${\scriptstyle -u_1}$}}
\rput(1.2,0){\rput(.2,2.15){${\scriptstyle u_2}$}\rput(.075,.35){${\scriptstyle -u_2}$}}
\rput(2.3,0){\rput(.2,2.15){${\scriptstyle u_3}$}\rput(.075,.35){${\scriptstyle -u_3}$}}
\rput(4.4,0){\rput(.2,2.15){${\scriptstyle u_{n-1}}$}
\rput(.075,.35){${\scriptstyle -u_{n-1}}$}}
\rput(5.7,0){\rput(.2,2.15){${\scriptstyle u_n}$}\rput(.075,.35){${\scriptstyle -u_n}$}}
\end{pspicture}
}}}\ignorespaces

\parshape=4  0pt.55\hsize 0pt.55\hsize 0pt.55\hsize 0pt\hsize
For type $B$ and the isomorphism 
from $W(B_n)$ to the group $\SS^\pm\{u_1,\ldots,u_n\}$
of signed permutations of the $\{u_i\}$ of \S\ref{subsection:weyl_groups},
the Coxeter element $\xi$
is the signed permutation at right.
The root lattice $L(B_n)$ is spanned by the $\{x_i\}$ with $x_i=v_i$, $(i<n)$ and
$x_n=\sqrt{2}v_n$, and we can then read off the diagram that
$\aa(x_1)=x_1+x_{q+1}+\cdots+x_{kq+1}$ for $k=(2^{p-1}-1)$ when $q>1$, or
$\aa(x_1)=2(x_1+\cdots+x_{n-1})+x_n$ when $q=1$. 
In either case
$\aa$ is not the zero map in $\endo(L/2)$. 
We can do $F_4$ by a similar brute force method to that employed for $E_6$
in the proof of the previous result, 
checking that the map $\aa=1+\xi^3$ is not zero in $\endo(L/2)$.
In type $D_n$ ($n$ odd) one can check that,
$$
\aa(x_{n-1})=-(2^{p-1}-1)(x_{n-1}-x_n)
-\kern-3mm\sum_{k=1}^{2^{p-1}-1} k(x_{kq}+x_{kq+1}+\cdots+x_{(k+1)q-1}),
$$
giving $\aa(\ov{x}_{n-1})\in\ker(\xi^{h/2}+1)\setminus\im(\xi^{h/2}+1)$,
while in $E_6$, we have $\aa(\ov{x}_1)=\ov{x}_2+\ov{x}_3+\ov{x}_4+\ov{x}_6$,
also hitting the target.
\qed
\end{proof}

One can also play these games using a Coxeter element for some visible subgroup
rather than all of $W(\Psi)$, and this will pay dividends later
when the visible subgroup has Coxeter number with a larger $2$-part. 
For example, $W(E_6)$ has $h=2^2 3$, whereas the visible 
$W(D_5)=\langle s_2,\ldots,s_6\rangle$ has $h=2^3$. We then have for
$\xi=s_2\ldots s_6$ that
$$
\langle\ov{x}_2+\ov{x}_3+\ov{x}_5,\ov{x}_2+\ov{x}_6\rangle
=\im(\xi^{h/2}+1)\subset\ker(\xi^{h/2}+1)
=\im(\xi^{h/2}+1)\oplus\langle\ov{x}_4,\ov{x}_5\rangle,
$$
in $L/2$ for $L$ the $E_6$ root lattice. Moreover, $\aa=1+\xi+\xi^2+\xi^3$,
gives $\aa(\ov{x}_1)
=\ov{x}_3+\ov{x}_4+\ov{x}_6\in\ker(\xi^{h/2}+1)\setminus\im(\xi^{h/2}+1)$.

%\begin{proposition}\label{section:lattices:result300}
%\marginlabel{\small{will need to be redone}}
%If $\Psi$ %non simply-laced, connected, crystallographic, 
%is the symbol of a non simply-laced irreducible Weyl group,
%then the admissible
%(resp. specially admissible) pairs $(\Psi,s)$ are given by the following, where
%$s=\BoxedEPSF{admissible_vertex.eps scaled 300}$ 
%(resp $s=\BoxedEPSF{special_admissible_vertex.eps scaled 250}$):
%$$
%\begin{pspicture}(0,0)(15,2)
%%\showgrid
%\rput(4,1){\BoxedEPSF{Bn2.eps scaled 300}}
%\rput(6.2,1.2){$4$}\rput(0.6,1){$B_n$}
%\rput(-3.9,0.5){\rput{90}(5.9,.8){$\left.\begin{array}{c}
%\vrule width 0 mm height 15 mm depth 0 pt\end{array}\right\}$}}
%\rput(2,1.7){$A_k (k\text{ even})$}
%\rput(-3.2,-.15){\rput{270}(5.9,.8){$\left.\begin{array}{c}
%\vrule width 0 mm height 31 mm depth 0 pt\end{array}\right\}$}}
%\rput(2.7,.2){$A_k (k\equiv 2\text{ mod }4)$}
%\rput(11,1){\BoxedEPSF{F42.eps scaled 300}}
%\rput(11,1.2){$4$}\rput(8.7,1){$F_4$}
%\end{pspicture}
%$$
%\end{proposition}

%\begin{proof}
%ditto.
%\qed
%\end{proof}

%We finish the section with a very elementary,

\section{Normal torsion free subgroups}\label{section:homomorphisms}

In this section we construct, for a large class of Coxeter groups, normal torsion
free subgroups whose index can be determined in terms of data associated to the
Coxeter group. 

\parshape=5 0pt\hsize 0pt\hsize 0pt.65\hsize 0pt.65\hsize 
0pt.65\hsize %0pt.7\hsize 0pt\hsize
Let $(W,S)$ be a Coxeter group with $S=S_0\cup\{t_1,\ldots,t_m\}$ and symbol 
$\Gamma$ at right, 
where $W(\Psi)=(W_{S_0},S_0)$ is an irreducible Weyl group of rank $n$,
and for each $s_i$ $(1\leq i\leq m)$, the pair
$(\Psi,s_i)$ is admissible.
%where the subsymbol
%$\Psi$ is connected,
%crystallographic. %and each pair $(\Psi,\aa_i)$ is either admissible or specially
%admissible. 
%with nodes $\{s_1,\ldots,s_n\}$. 
%Write
%$S=\{s\in\Psi,t_1,\ldots,t_m\}$ for the Coxeter generators of $W(\Gamma)$,
Let $\ss:W(\Psi)\rightarrow\gl(V)$ be the reflectional representation
of the Weyl group, $L\subset V$ the root lattice, 
and $u_i=u_{s_i}$ $(1\leq i\leq m)$ the vectors in $L$ given by Lemma
\ref{subsection:weyl_groups:result100}.
%write $v_s\in V$ for the vector
%corresponding to $s\in\Psi$, and $\ss_s=\ss(s)\in\gl(V)$ for the reflection
%in $v_s$. %(and hence the image under the reflectional representation
%of the Coxeter generator $s$).
%We will identify the Coxeter generator corresponding to $\aa\in\Psi$
%with the reflection $s_\aa\in\gl(V)$ in the root $\aa\in V$,
%and so by a slight abuse, we will think of $W(\Gamma)$ as having
%Coxeter generators $S=\{s_\aa (\aa\in\Psi),s_i (1\leq i\leq k)\}$.
\vadjust{\hfill\smash{\lower 25mm
\llap{
\begin{pspicture}(0,0)(4.5,4)
%\showgrid
\rput(-.1,.15){
\rput(2,2){\BoxedEPSF{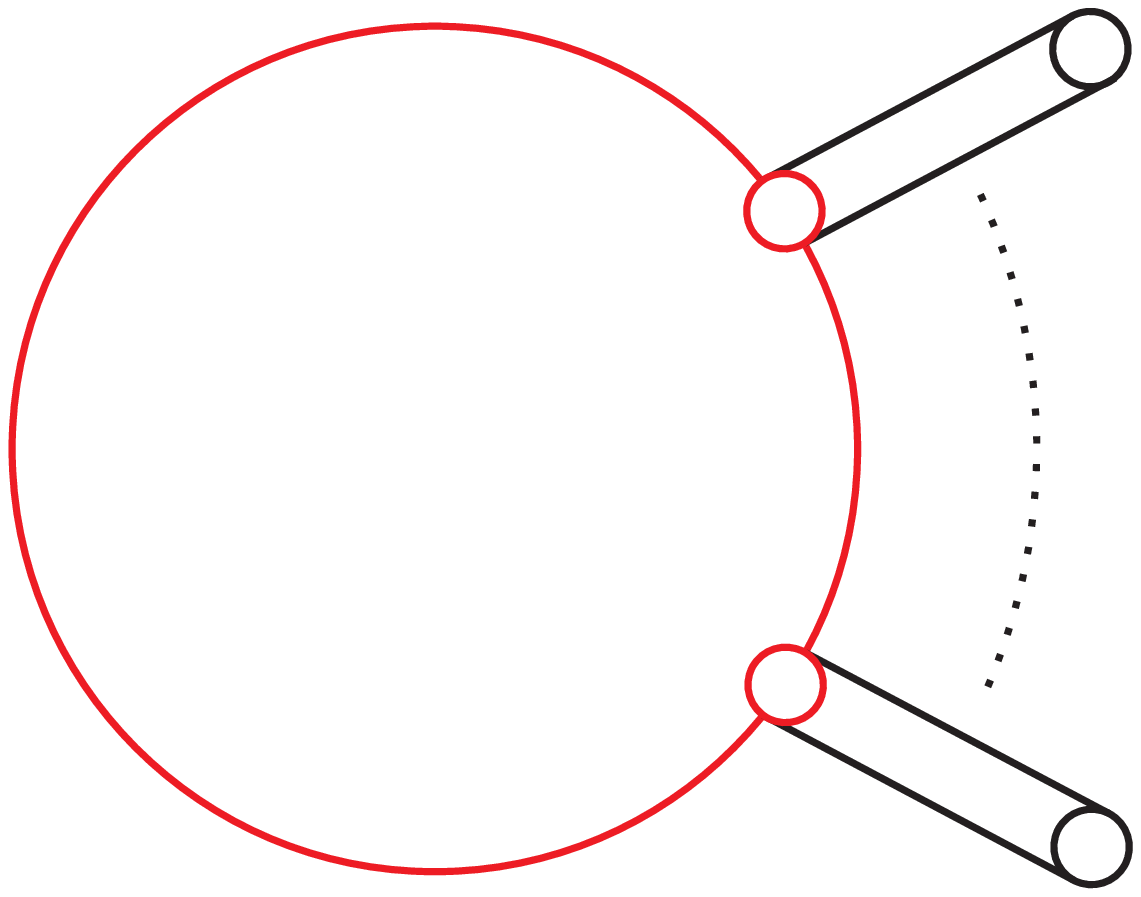 scaled 350}}
\rput(1.5,2){${\red \Psi}$}
\rput(4.4,2){$(\dag)$}
\rput(2.5,2.65){${\red s_1}$}\rput(2.5,1.45){${\red s_m}$}
%\rput(3.25,3.35){$4$}\rput(3.25,.6){$4$}
\rput(3.8,3.8){$t_1$}\rput(3.8,.2){$t_m$}
}
\end{pspicture}}}}\ignorespaces

\parshape=6 0pt.65\hsize 0pt.65\hsize 0pt.65\hsize 
0pt.65\hsize 0pt.65\hsize 0pt\hsize
Let $\Lambda_i:=\Lambda_{s_i}$ $(1\leq i\leq m)$.
We will flip between the Abelian group and vector space structures
on $\Lambda_i/2$, with the context making clear which one
we mean.
Extend the $W(\Psi)$-action on $\Lambda_i/2=L/2$,
to the direct product $\prod\Lambda_i/2$ diagonally, ie:
making $g\in W(\Psi)$ send
$(\bar{v}_1,\ldots,\bar{v}_m)$ to $(g(\bar{v}_1),\ldots,g(\bar{v}_m))$.
Form the semi-direct product 
$\prod\Lambda_i/2\rtimes W(\Psi)\cong(\Z/2)^{mn}\rtimes W(\Psi)$
induced by this action,
every element of which can be expressed 
as an ordered pair $(\v,g)$, where $g\in W(\Psi)$ and 
$\v=(\bar{v}_1,\ldots,\bar{v}_m)$,
with $g$ and $\v$ (although not the $v_i$) uniquely defined. Moreover, 
$\prod\Lambda_i/2$ is generated by the $g$-conjugates, $g\in W(\Psi)$,
of the vectors, $\u_i=(0,\ldots,\bar{u}_i,\ldots,0)$,
with the non-zero entry in the $i$-th position, and
thus $\prod\Lambda_i/2\rtimes W(\Psi)$ is generated by the $s\in S_0$ and the 
$\u_i$ (it is because of the ease of this observation that we work
here with $\Lambda_i/2$, rather than $L/2$). Finally, 
the order of %$\prod\Lambda_i/2\rtimes W(\Psi)$
the semi-direct product
is $2^{mn}|W(\Psi)|$, with the Weyl group order the product of the $m_i+1$,
for $m_i$ the exponents 
in Table \ref{table:weyl_data}.
%Define a map $f:S\rightarrow\Z_2^m$ that sends $s\in\Psi$ to 
%$\mathbf{0}\in\Z_2^m$ and
%$t_i$ to the $m$-tuple with $1$ in the $i$-th position and zeroes elsewhere.
%As the Coxeter relations are satisfied in $\Z_2^m$ by the $f(S)$,
%the map $f$ extends to a homomorphism
%$W(\Gamma)\rightarrow\Z_2^m$. If $g\in W(\Gamma)$ and
%$w$ is a word in the $S$ representing $g$ that contains $\lambda_g^i$ occurences
%of the generator $t_i$, then $f(g)$ is the $m$-tuple with %\ov{\lambda}_g^i
%$\lambda_g^i\text{ (mod }2)$ in the $i$-th position. In particular, the parity
%$\lambda_g^i\text{ (mod }2)$ of $\lambda_g^i$ is independent of the word $w$, 
%as is 
%$\lambda_g=\sum\lambda_g^i\text{ (mod }2)$.

We now define a homomorphism 
$$
\vphi:W(\Gamma)\rightarrow \prod\Lambda_i/2\rtimes W(\Psi),
$$
as follows: let $\vphi(s)=(\0,s)$ for $s\in S_0$, 
$\vphi(t_i)=(\u_i,1)$, and extend linearly. %Then $\vphi$ is a
%homomorphism if it is well defined, for which it suffices to show that the Coxeter
%relations are satisfied by the $(\0,s):=(\0,\ss(s))$ and $(\u_i,1)$.
Clearly, $\vphi(t_i)^2=(\u_i,1)^2=(\0,1)$, 
$(\vphi(t_i)\vphi(t_j))^{2}=(\u_i,1)(\u_j,1)^2=(\0,1)$
and $(\vphi(s)\vphi(s'))^{m_{ss'}}=(\0,1)$ for $s,s'\in S_0$.
Thus for $\vphi$ to be a homomorphism it remains to check
that the $(st_i)^{m_{st_i}}=(\0,1)$ relations are satisfied in the image for
$s\in S_0$. As the $u_i$ are orthogonal to the 
$v_s$ for $s\in S_0\setminus\{s_i\}$, we have
$\u_i^s:=s\u_i s^{-1}=\u_i$
in the semi-direct product, hence
$$
((\0,s)(\u_i,1))^2=(\u_i+\u_i^s,1)=(\0,1),
$$
for these $s$. Finally, for $1\leq i\leq m$,
$$
((\0,s_i)(\u_i,1))^4=(2\u_i+2\u_i^s,1)=(\0,1),
$$
and so $\vphi$ is a homomorphism, surjective
as the semi-direct
product is generated by the images of the $s\in S_0$ and the $t_i$.

We now show that when the pairs $(\Psi,s_i)$ are specially admissible,
the kernel of $\vphi$ is torsion free, otherwise, if some are merely admissible,
then it is \emph{almost\/} torsion free, and the offending torsion can be 
pinpointed pretty accurately.
The first step is to observe that if 
$\Delta_1,\Delta_2\subset\Gamma$ are disjoint subsymbols, then the
$\vphi W(\Delta_i)$ intersect trivially in $\prod\Lambda_i/2\rtimes W(\Psi)$.
For, if $\Delta\subset\Gamma$ is a subsymbol and $g\in W(\Delta)$,
then $\vphi(g)=(\v,g')$ where $g'\in W(\Delta\cap\Psi)$,
and $\v\in\prod\Lambda_i/2$ has $i$-th coordinate non-zero only if
$t_i\in\Delta$. The result now follows as the visible subgroups
$W(\Delta_i)$ intersect trivially, and by the uniqueness of the $(\v,g)$
expression for elements of $\prod\Lambda_i/2\rtimes W(\Psi)$.

In particular, by
the comments following Theorem \ref{torsion}, 
if $\vphi$ is faithful when restricted to disjoint visible subgroups 
$W(\Delta_1),$ $W(\Delta_2)$, then
it is also faithful on $W(\Delta_1\coprod\Delta_2)$, and thus to show that
$\ker\vphi$ is torsion free,
it suffices to show that $\vphi$ is faithful on the finite visible $W(\Delta)$, for
\emph{connected\/} $\Delta\subset\Gamma$.

Before proceeding, we mention in passing that one may be tempted by a simpler 
choice for $\vphi$, 
namely letting $\Lambda$
be the lattice spanned by the $W(\Psi)$-orbit of all the $\{u_1,\ldots,u_m\}$
simultaneously,
forming $\Lambda/2\rtimes W(\Psi)$ and defining $\vphi(s)=(0,s)$ and
$\vphi(t_i)=(u_i,1)$. It turns out that this $\vphi$ is also faithful
on the finite visible $W(\Delta)$ for $\Delta\subset\Gamma$ connected,
but not generally on the finite visibles, because 
$\vphi W(\Delta_i)$ with disjoint $\Delta_i$ can have non-trivial
intersection.
For example, if $\Psi$ is of type $E_6$ with the two admissible
nodes,
$$
\begin{pspicture}(0,0)(7.2,2.5)
%\showgrid
\rput(-3.4,0){
\rput(7,1.25){\BoxedEPSF{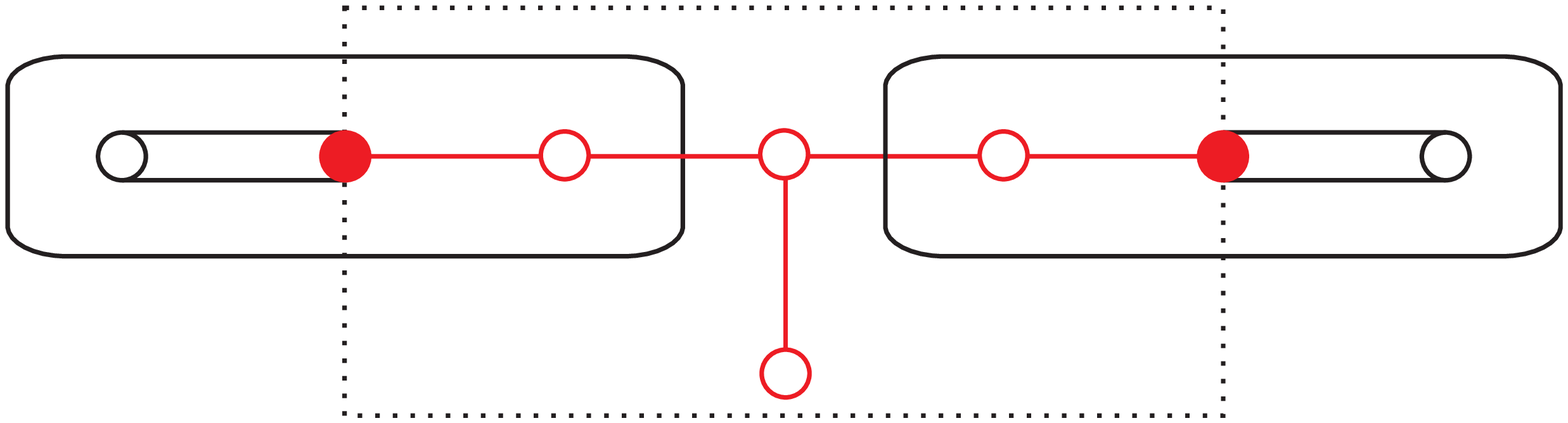 scaled 300}}
\rput(5.6,.5){${\red \Psi=E_6}$}
\rput(3.8,.7){$\Delta_1$}\rput(10.2,.7){$\Delta_2$}
\rput(3.5,1.5){$t_1$}\rput(5.1,1.75){${\red s_1}$}\rput(5.7,1.3){${\red s_3}$}
\rput(10.525,1.5){$t_2$}\rput(8.9,1.75){${\red s_2}$}\rput(8.35,1.3){${\red s_4}$}
%\rput(3.45,1.85){$4$}\rput(10.55,1.85){$4$}
}
\end{pspicture}
$$
then for $i=1,2$, both $1\not=(s_it_is_{i+2})^3\in W(\Delta_i)$
map via this $\vphi$
to the same non-trivial element of $\Lambda/2\rtimes W(\Psi)$.

\begin{lemma}\label{section:homomorphisms:result200}
Let $B_k\subset\Gamma$ 
%be a connected subsymbol with $W(\Delta)$ a Weyl group of type $B_n$ 
contain one (hence exactly one) of the $t_i$,
$$
\begin{pspicture}(0,0)(14,1)
%\showgrid
\rput(7,.5){\BoxedEPSF{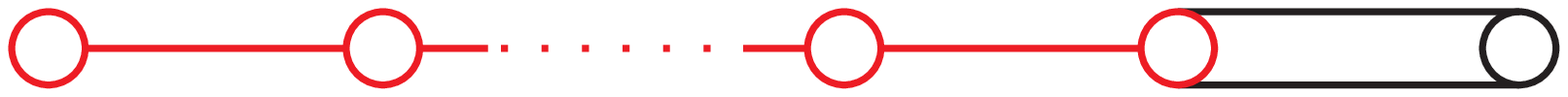 scaled 400}}
\rput(6.4,.8){${\red \Delta}$}
%\psline[linewidth=.2mm]{->}(3,.5)(3.7,.5)
%\rput(14,.5){$(\dag)$}
%\rput(3.95,.1){$x_1$}\rput(5.35,.1){$x_2$}\rput(7.2,.1){$x_{k-2}$}
%\rput(8.6,.1){$s_i$}\rput(10.1,.1){$t_i$}
%\rput(9.4,.75){$4$}
\end{pspicture}
$$
Then $\vphi$ is faithful on $W(B_k)$ if and only if
the $W(\Delta)$-orbit of $\ov{u}_i\in L/2$ spans a subspace $U_i\subset L/2$ of 
dimension $k$.
\end{lemma}

In particular, by the comments following the
definition of admissibility in \S\ref{section:lattices},
if the pair $(\Psi,s_i)$ is admissible 
(respectively specially admissible) 
then $\vphi$ is faithful on the visible $W(B_k)$
with $t_i\in B_k$ and $k$ \emph{even\/} (resp. for all $k$).
The result follows since
$W(B_k)\cong(\Z/2)^k\rtimes W(\Delta)$, whereas, 
$\vphi W(B_k)\cong U_i\rtimes W(\Delta)$, thinking of $U_i$ as an Abelian group.

We are now ready for,

\begin{theorem}\label{section:homomorphisms:result300}
Let $\Gamma$ be the symbol $(\dag)$ with $W(\Psi)$ an irreducible Weyl group
of rank $n$ and exponents $m_i$, and
the pairs $(\Psi,s_i)$
specially admissible for all $i$. Then 
$\text{ker}\,\vphi^{}\subset W(\Gamma)$ is a torsion free subgroup
of index 
$$
[W(\Gamma):\ker\vphi]=2^{mn}\prod_{i=1}^n(m_i+1).
$$
\end{theorem}

\begin{proof}
We have observed that it suffices to show that
$\vphi$ is faithful on the finite visible $W(\Delta)$ for
connected $\Delta\subset\Gamma$.
This follows immediately for $\Delta\subset\Psi$,
and if 
$\Delta=\{t_i\}$, then the faithfulness follows as 
$\bar{u}_i\not=0$. The remaining possibility is that $\Delta$ is of type
$B_k$ for $k>1$ and containing (exactly) one of the $t_i$ generators 
and this case follows by special admissibility
and the comments following Lemma \ref{section:homomorphisms:result200}.
\qed
\end{proof}

If some of the pairs are admissible but not specially admissible, we get
elements of finite order in $\text{ker}\,\vphi$, but we can say
where, and what they look like:

\begin{proposition}\label{section:homomorphisms:result400}
Let $\Gamma$ be the symbol $(\dag)$ with the pairs $(\Psi,s_i)$
admissible, and let 
$W(\Delta)$ be finite visible with $\Delta\subset\Gamma$ connected. %and
%$\Delta\subset\Gamma$ be a connected, crystallographic
%subsymbol with 
%$t_i\in\Delta$ for some (hence exactly one) $i$. 
Then 
$\text{ker}\,\vphi^{}\subset W(\Gamma)$ contains a
non-trivial element $g\in W(\Delta)$ of prime order
%$p$ 
if and only if there is $\Upsilon=\coprod\Delta_i\subset\Delta$ with the $\Delta_i$
connected, and
\begin{description}
\item[(i).] $\Delta_1$ of type $B_k$ for $k$ odd, and with a $t_i\in\Delta_1$ for some 
(hence exactly one) $i$;
\item[(ii).] $\Delta_i$ is of $(-1)$-type and contained in $\Psi$ for $i>1$;
\item[(iii).] $g$ is conjugate to 
%$\text{ker}\,\vphi$ contains 
the element of longest length in $W(\Upsilon)$.
\end{description}
\end{proposition}

\begin{proof}
As in the proof of Theorem \ref{section:homomorphisms:result300}
we can reduce to the case that $\Delta$ has type $B$ 
and contains (exactly) one of the $t_i$ generators. 
%Now, $\text{ker}\,\vphi$ contains a non-trivial element from $W(\Delta)$ 
%if and only if it contains an element of order $p\geq 2$ a prime. 
Any element of prime order $p>2$ in $W(\Delta)$ is conjugate by Lemma 
\ref{section:preliminaries:result200} to an element in $W(\Delta\cap\Psi)$.
As $\ker\vphi$ contains no such elements, it can only contain prime 
order elements of
order $2$. By Theorem \ref{section:preliminaries:torsion:result300}
any such $g$ is conjugate to the element of longest length in a finite visible
subgroup of $(-1)$-type, ie: there is $\Upsilon=\coprod\Delta_i\subset\Delta$ with 
the $\Delta_i$ connected and of $(-1)$-type, and with 
$w_{\Upsilon}^{}\in \text{ker}\,\vphi$
for $w_{\Upsilon}^{}$ the element of longest length in $W(\Upsilon)$.

If all the $\Delta_i\subset\Psi$ then $\vphi$ is faithful on the $W(\Delta_i)$,
hence on the $W(\Upsilon)$, so we must have one
of the $t_i$ contained in $\Delta_1$, say, and thus $\Delta_1$ of type $B_k$.
As $\Delta$ is connected we must also have the $\Delta_i\subset\Psi$ for 
$i>1$. If $k$ is even, then $\vphi$ is 
faithful on $W(\Delta_1)$ by the admissibility of $(\Psi,s_1)$
and Lemma \ref{section:homomorphisms:result200}, and hence on 
$W(\Upsilon)$. Thus, $\Delta_1=B_k$ for $k$ odd as required.
\qed
\end{proof}

%Comparing this to Propositions 
%\ref{section:lattices:result200}-\ref{section:lattices:result300}, we see that the
%possibility described in the Proposition always arises, that is, for any admissible
%pair $(\Psi,s)$, there is a subsymbol $B_3\subset\Gamma$ satisfying
%(i)-(iii) above.

In order to get normal torsion free subgroups when some of the pairs are
admissible but not specially admissible, 
we need to re-jig $\vphi$ a little to handle the torsion
described in the Proposition. 
The idea is simple enough: we delve deeper into $\text{ker}\,\vphi$
to avoid the elements of longest length in type $B$ groups of odd 
rank.
To this end, fix $i$ and
let $g\in W(\Gamma)$ and
$w$ a word representing it that contains $\lambda_g^i$ occurrences
of the generator $t_i$. %then $f(g)$ is the $m$-tuple with %\ov{\lambda}_g^i
%$\lambda_g^i\text{ (mod }2)$ in the $i$-th position. In particular, 
Define a map $f_i:S\rightarrow\Z/2$ that sends $t_i$ to $0$ and
$s\not=t_i$ to $1$.
As the Coxeter relations are satisfied in $\Z_2$ by the $f_i(S)$,
the map extends to a homomorphism $W\rightarrow\Z/2$, with 
$f_i(g)=\lambda_g^i\text{ (mod }2)$.
Thus the parity $\lambda_g^i\text{ (mod }2)$ of $\lambda_g^i$ is independent 
of the word $w$. %Moreover, as $\lambda_{gh}^i=\lambda_g^i+\lambda_h^i$, the
%map $W(\Gamma)\rightarrow\Z/2$ given by $g\mapsto\lambda_g^i$ is a 
%homomorphism. 
Observe that as
\begin{equation}\label{typeB_longest_length}
w_{B_n}^{}=(s_1s_2\cdots s_{n-1}s_ns_{n-1}\cdots s_2s_1)
%\cdot(s_2s_3\cdots s_n\cdots s_3s_2)
\cdots(s_{n-2}s_{n-1}s_ns_{n-1}s_{n-2})\cdot(s_{n-1}s_ns_{n-1})\cdot s_n,
\end{equation}
we have $f_i(w_{B_n}^{})=n\text{ (mod }2)$.

Let the nodes $\{s_i\}$ giving the admissible pairs be
divided $\{s_1,\ldots,s_\ell,s_{\ell+1},\ldots,s_m\}$, with
the first $\ell$ admissible but not specially admissible and the remaining 
specially admissible.
Define
$$
\varepsilon:W(\Gamma)\rightarrow\Z/2^\ell:=\prod_{\ell\text{ copies}}\Z/2,
$$
by $\varepsilon(g)=\x_g=(f_1(g),\ldots,f_\ell(g))\in\Z/2^\ell$, the $f_i$ as above, 
a product of homomorphisms,
hence a homomorphism, and moreover,
constant on each conjugacy class of $W(\Gamma)$.
Now define
$$
\wvphi:W(\Gamma)\rightarrow\Z/2^\ell\times(\prod\Lambda_i/2\rtimes W(\Psi)),
$$
by $\wvphi=\varepsilon\times\vphi$, ie: 
$\wvphi(g)=(\varepsilon(g),\vphi(g))$,
%where $\x=(\lambda_g^1,\ldots,\lambda_g^m)\in\Z_2^m$.
%Then $\wvphi$, being a product of $m+1$ homomorphisms, is also a homomorphism.
also a homomorphism.

For each $(\Psi,s_i)$, $1\leq i\leq \ell$, there must be a subsymbol
$B_k\subset\Gamma$ of odd rank $k$ containing the Coxeter generator $t_i$,
and with $\vphi$ not faithful on $W(B_k)$,
since the pair is admissible but not specially admissible. 
In particular, $\text{ker}\,\vphi$ contains the element 
$w_i=w_\Upsilon^{}$ of longest length
in $W(\Upsilon)$ with $\Upsilon=\coprod\Delta_j\subset B_k$ and 
$\Delta_1$ also a type $B$ of odd rank by Proposition
\ref{section:homomorphisms:result400}. Then 
by (\ref{typeB_longest_length}),
$f_i(w_i)=1$,
%so that $g_0\not\in\text{ker}\,\ve$.
and $f_j(w_i)=0$ for $j\not=i$. The result is that 
$\vphi(w_i)=(\0,1)$ and 
$$
\ve(w_i)=\x_i=(0,\ldots,1,\ldots,0),
$$
with the $1$ in the $i$-th position. Thus, $\wvphi(w_i)=(\x_i,\0,1)$, and we
already have that $\wvphi(s\in S_0)=(\0,\0,s)$ and $\wvphi(t_i)=(\x_i,\u_i,1)$,
so that $\wvphi(w_it_i)=(\0,\u_i,1)$. As the group 
$\Z/2^\ell\times(\prod\Lambda_i/2\rtimes W(\Psi))$ is generated by the
$(\x_i,\0,1)$, $(\0,\0,s)$ and $(\0,\u_i,1)$, we get that $\wvphi$
is surjective.

If $\Delta\subset\Gamma$, then the image under 
$\wvphi$ of $g\in W(\Delta)$ is $(\x_g,\v,g')$ with
$g'\in W(\Delta\cap\Psi)$, and the $\v\in\prod\Lambda_i/2$ and
$\x_g\in\Z/2^\ell$ having $i$-th coordinate non-zero only if
$t_i\in\Delta$. Thus, $\wvphi$, like $\vphi$, sends trivially
intersecting visibles to trivially intersecting subgroups,
so that
if $\wvphi$ is faithful on
$W(\Delta_i)$ with the $\Delta_i\subset\Gamma$ disjoint, then
it is also faithful on $W(\Delta_1\coprod\Delta_2)$.

\begin{theorem}\label{section:homomorphisms:result500}
Let $\Gamma$ be the symbol $(\dag)$ with 
$W(\Psi)$ an irreducible Weyl group
of rank $n$ and exponents $m_i$, and with
the pairs $(\Psi,s_i)$, $(1\leq i\leq\ell)$
admissible but not specially admissible and the remaining pairs
specially admissible. Then 
$\text{ker}\,\wvphi\subset W(\Gamma)$ is a torsion free subgroup 
of index 
$$
[W(\Gamma):\ker\wvphi]=2^{mn+\ell}\prod_{i=1}^n(m_i+1).
$$
\end{theorem}

\begin{proof}
It suffices as before to show that $\wvphi$ is faithful on the finite
visible $W(\Delta)$ with $\Delta$ connected, which in turn follows if
$\text{ker}\,\wvphi$ contains no elements of prime order from such subgroups. 
%Let $g$ be such an 
%element of such a subgroup contained in $\text{ker}\,\wvphi$. 
Such an element $g\in\text{ker}\,\wvphi$
is contained in both $\text{ker}\,\ve$ and $\text{ker}\,\vphi$, and the second
of these, combined with Proposition \ref{section:homomorphisms:result400},
gives that $g$ is conjugate to the element $w_\Upsilon^{}$ of longest length
in $W(\Upsilon)$ where the $\Upsilon=\coprod\Delta_j$ have the properties 
described in the Proposition. As $\ve$ is constant on the conjugacy classes
we have $w_\Upsilon^{}\in\text{ker}\,\ve$ also. But $\lambda_{w_\Upsilon}^i=1$, 
as the element
of longest length in a type $B$ group of odd rank contains an odd number 
of occurrences of $t_i$ and thus $\varepsilon(w_\Upsilon^{})\not=\0$. 
This contradiction gives the result.
\qed
\end{proof}

%Comment that all of the above still works, and so $\text{ker}\,\wvphi$ is
%torsion free, when $W(\Psi)$ is an arbitrary crystallographic group (but
%the kernel doesn't now have finite index).

%\section{The image and its conjugacy classes}
%\section{The finite group $W(\Gamma)/\text{ker}\,\wvphi$}
%\section{Conjugacy in $\mathbf{\Z}/2^\ell\times(\prod\Lambda_i/2\rtimes W(\Psi))$}
\section{Even bigger, non-normal, torsion free subgroups}
\label{section:extended}

In the last section we constructed normal torsion free subgroups
in our family of Coxeter groups; in this, we extend them to larger
non-normal torsion free subgroups. The strategy is simple: we look 
for subgroups
$H$ in the finite image $W(\Gamma)/\text{ker}\theta$,
where $\theta=\vphi$ or $\wvphi$, and with the 
property that no non-trivial element of finite order in $W(\Gamma)$ maps via 
$\theta$ into $H$. 
Thus, $\theta^{-1}H\subset W(\Gamma)$ will be a torsion free subgroup
containing $\ker\theta$, but with $1/|H|$-th the index. Moreover, both the
normality of the original group and the construction of the new
one will have geometrical interpretations in the next section.

For simplicity we restrict our attention to cyclic $H$ of order
$2^p$ for some $p$, and work just with $\wvphi$, pointing out where appropriate
the corresponding results for $\vphi$.
Thus,
with $\wvphi:W(\Gamma)\rightarrow\Z/2^\ell\times(\prod\Lambda_i/2\rtimes W(\Psi))$
from \S\ref{section:homomorphisms}, we seek an
$$
H=\Z/(2^p)\subset\Z/2^\ell\times(\prod\Lambda_i/2\rtimes W(\Psi)),
$$ 
with $\wvphi^{-1}H\subset W(\Gamma)$ torsion free.

Firstly, if $\wvphi^{-1}H$ contains non-trivial torsion, then it contains
$2$-torsion:
if $g\in\wvphi^{-1}\Z/(2^p)$ has finite order then 
as $\wvphi$ kills no torsion in $W(\Gamma)$, 
the image $\wvphi(g)$ has the same order in $\Z/(2^p)$,
and thus $g$ has even order; some power of $g$ thus has
order $2$, with order $2$ image in $\Z/(2^p)$. We conclude that
$\wvphi^{-1}\Z/(2^p)$ is torsion free if and only if 
there are no elements of order $2$ in $W(\Gamma)$ mapping via $\wvphi$
to $H$; alternatively, $H$ intersects
trivially any conjugacy class of order two 
in $\Z/2^\ell\times(\prod\Lambda_i/2\rtimes W(\Psi))$ that is the image of a 
conjugacy class of elements of order $2$ in $W(\Gamma)$.

Thus, our efforts will be focused on the conjugacy classes of elements of
order $2$ in $W(\Gamma)$ and $W(\Gamma)/\text{ker}\wvphi$, as well as the relationship
between them given by $\wvphi$. 
Let $\JJ_\Gamma$ be the $W(\Gamma)$-equivalence classes of subsymbols of 
$(-1)$-type from \S\ref{subsection:torsion}, and for $\Delta\in\JJ_\Gamma$,
let $w_{\kern-1.5pt\Delta}^{}$ 
be the element of longest length
in $W(\Delta)$. Although the empty symbol is not of $(-1)$-type, we 
adopt the convention $w_\varnothing=1$.
If $\Delta\in\JJ_\Gamma$, write $\CC(\Delta)$
for the conjugacy class of $w_{\kern-1.5pt\Delta}^{}$ in $W(\Delta)$.

%Now to the image group. 
%Let $A=\prod\Lambda_i/2$ and 
%$\endo(A)$ its endomorphism algebra. We have the 
%$W(\Psi)\rightarrow\aut(A)\hookrightarrow\endo(A)$,
%giving rise to the semi-direct product
%$A\rtimes W(\Psi)$. Write $1\in\aut(A),0\in\endo(A)$ for the identity automorphism
%and zero endomorphism. If $g\in W(\Psi)$ is an involution, observe that
%$(g+1)^2=0$ in $\endo(A)$, so that $\im(g+1)\subset\ker(g+1)\subset A$.
%Let $A_g$ be the quotient space
%$\text{\ker}(g+1)/\text{im}(g+1)$, and in particular, write $A_\Delta$
%if $g=w_{\kern-1.5pt\Delta}^{}$.

We can get a very rough parametrization of the conjugacy classes 
$2$-torsion in the image group. First, we extend the notation of
\S\ref{section:lattices}, so that if $\aa\in\endo(\prod\Lambda_i/2=\prod L/2)$,
then $\ker\aa$ and $\im\aa$ are the direct sums of the
kernels and images in $\endo(\Lambda_i/2)$. 
Let $\TT$ be the set of triples 
$(\x,\v,\Delta)$ where $\x\in\Z/2^\ell$, 
$\v\in\prod\Lambda_i/2$
and $\Delta\in\JJ_\Psi\cup\{\varnothing\}$.

\begin{lemma}\label{section:extended:result100}
The map 
$(\x,\v,\Delta)\mapsto(\x,\v,w_{\kern-1.5pt\Delta}^{})$ 
induces a surjection from $\TT$ to the set of conjugacy classes of involutions
in the group $\Z/2^\ell\times(\prod\Lambda_i/2\rtimes W(\Psi))$.
Writing $\CC(\x,\v,\Delta)$ for the class that is the image of 
$(\x,\v,\Delta)\in\TT$, we have
\begin{description}
\item[(i).] if $\CC(\x_1,\v_1,\Delta_1)=\CC(\x_2,\v_2,\Delta_2)$ then
$\x_1=\x_2$ and $\Delta_1=\Delta_2$;
%\item[(ii).] If $\Psi$ itself is of $(-1)$-type then 
%$\CC(\x_1,\v_1,\Psi)=\CC(\x_2,\v_2,\Psi)$ if and only if 
%$\x_1=\x_2$ and $\v_1,\v_2$ lie in the same $W(\Psi)$-orbit on 
%$\ker(w_\Psi+1)/\im(w_\Psi+1)=L/2$.
\item[(ii).] the class
$\CC(\0,\0,\Delta)=\bigcup\{(\0,\u,g)\,|\,\u\in\im(g+1)\}$,
where the union is over all $W(\Psi)$-conjugates $g$ of $w_\Delta$.
\end{description}
\end{lemma}

\begin{proof}
If $(\x,\v,g)\in\Z/2^\ell\times(\prod\Lambda_i/2\rtimes W(\Psi))$, then $(\x,\v,g)^2
=(\0,(g+1)\v,g^2)$, so that the element is an involution if and only if
either $g$ is the identity or $g$ is an involution and $\v\in\text{ker}(g+1)$. 
%Moreover, if $\v_1,\v_2\in A$
%differ by $\w\in\text{im}(g+1)$, then conjugating $(\x,\v_1,g)$ 
%by $(\0,\u,1)$, for $\w=(g+1)\u$, gives $(\x,\v_2,g)$.
Thus every involution $(\x,\u,g)$ is conjugate to 
a $(\x,\v,w_{\kern-1.5pt\Delta}^{})$ with the corresponding
triple $(\x,\v,\Delta)\in\TT$, thus giving the surjectivity of the map.
Clearly elements $(\x_1,\v_1,g_1)$ and $(\x_1,\v_2,g_2)$ are conjugate
only if $\x_1=\x_2$ and the $g_i$ are conjugate in $W(\Psi)$, and this gives
the first part.
%
%If $\Psi$ itself has $(-1)$-type then $w_{\Psi}^{}$ is the identity
%on $A$ and so $A_\Psi=A$. Then for any $h\in W(\Psi)$, the conjugate
%of $(\x,\v,w_{\Psi}^{})$ by $(\y,\u,h)$ is
%$(\x,h(\v),w_{\Psi}^{})$, since $w_{\Psi}^{}$ is central
%and $w_{\Psi}^{}+1$ is the zero map in $\endo(A)$,
%so that the $W(\Psi)$-orbits on $A$
%parametrise the conjugacy classes as claimed.
The second part follows as $(\x,\u,h)(\0,\0,w_\Delta)(\x,\u,h)^{-1}
=(\0,(1+hw_\Delta h^{-1})(\u),hw_\Delta h^{-1})$.
\qed
\end{proof}

When the image group is $W(\Gamma)/\ker\vphi$, the above goes straight
through with pairs $(\v,\Delta)$ instead.
%Which of the classes in $W(\Gamma)/\ker\wvphi$ lie in the image of the classes
%in $W(\Gamma)$? 
%If $g\in\CC(\Theta)$
%then $\wvphi(g)\in\CC(\x,\v,\Delta)$, 
%where $\Delta$ is a subsymbol of $(-1)$-type in $\Theta\cap\Psi$.
%
%In particular, if $\Psi$ is of $(-1)$-type and $\Theta\in\JJ_\Gamma$
%is not $\Psi$, then $\Delta$ is a proper subsymbol of $\Psi$, and so not
%$W(\Psi)$-equivalent to $\Psi$ (equivalent symbols having the
%same rank). Thus,
%$\wvphi\,\CC(\Theta)\not=\CC(\x,\v,\Psi)$ for
%any $\x,\v$ in this case, whereas $\wvphi\,\CC(\Psi)=\CC(\0,\0,\Psi)$. Summarizing,

\begin{lemma}\label{section:extended:result200}
Let %$\Delta\in\JJ_\Psi$ be of maximal rank and 
$\Theta\in\JJ_\Gamma$.
%If $\Psi$ is a symbol of $(-1)$-type then exactly one of the conjugacy classes
%$\CC(\x,\v,\Psi)$ lies in the $\wvphi$ image of the classes $\CC(\Delta)$
%for $\Delta\in\JJ_\Gamma$.
\begin{description}
\item[(i).] $\Theta$ is the disjoint union of a $\Delta'\subset\Psi$
with symbols of type $B$ containing a $t_i$, and
$\wvphi\,\CC(\Theta)=\CC(\x,\v,\Delta')$ for some $\x,\v$.
\item[(ii).] Let $\Delta\in\JJ_\Psi$ be of maximal rank.
If $\wvphi\,\CC(\Theta)=\CC(\x,\v,\Delta)$ for some $\x,\v$,
then $\Theta$ is a disjoint union of a $\Delta'$ that is $W(\Psi)$-equivalent
to $\Delta$ together with some of the nodes $t_i$ where $s_i\not\in\Delta'$.
\item[(iii).] If $\CC(\0,\v,\Delta)=\wvphi\,\CC(\Theta)$ 
for $\Delta\in\JJ_\Psi$ of maximal rank,
then $\CC(\0,\v,\Delta)=\CC(\0,\0,\Delta)$.
\end{description}
\end{lemma}

It is the last part that is the really useful bit: only one of the classes of the form
$\CC(\0,\v,\Delta)$, namely
$\CC(\0,\0,\Delta)$, lies in the image of the $2$-torsion of $W(\Gamma)$, when
$\Delta$ is of maximal rank.

\begin{proof}
If $\Theta\in\JJ_\Gamma$ is connected and $\not\subset\Psi$ then it contains
one (hence exactly one) of the $t_i$ and so is of type $B$, giving the first
part of (i). Observe that the image under $\wvphi$ of the $w_{B_n}$ given
by (\ref{typeB_longest_length}) 
collapses to the identity in the third
component, and since $w_\Theta$ is the product of the elements of longest length
in the irreducible components of $W(\Theta)$, we get
$\wvphi(w_\Theta)=(\x,\v,w_{\Delta'})$ for some $\x$ and $\v$ as claimed.

In part (ii), $\Theta$ must now be the disjoint union of a $\Delta'$ 
that is $W(\Psi)$-equivalent to
$\Delta$ and some type $B$'s that contain various $t_i$'s. By checking through the
possible $\Psi$ in Table \ref{table:roots1}
and the 
$\Delta\subset\Psi$ with $\Delta\in\JJ_\Psi$ of maximal rank,
we can see that every $s\in\Psi$ not already in $\Delta$ is joined by an edge to it. 
Thus, the
only type $B$'s that remain disjoint from $\Delta$ are the single nodes
$t_i$ with $s_i\not\in\Delta'$.

Now, for $\Theta$ of the form given in part (ii), the element $w_\Theta$
is the product of $w_{\Delta'}$ with some of the $t_i$'s such that $s_i\not\in\Delta'$. 
These $s_i$ must be admissible but not specially admissible, as  
for every specially admissible pair $(\Psi,s)$ we have $s\in\Delta'$ by
Theorems \ref{section:lattices:result200}-\ref{section:lattices:result300}.
Thus $\wvphi(w_\Theta)=(\x,\v,w_{\kern-1.5pt\Delta}^{})$, where $\x$ 
contains a $1$ in the $i$-th coordinate for each $t_i$ present, by the
comments preceeding Theorem \ref{section:homomorphisms:result500}.
Since by Lemma \ref{section:extended:result100},
$\CC(\x,\v,\Delta)\not=\CC(\0,\v,\Delta)$ when $\x\not=\0$, the
only $\CC(\Theta)$ mapping to a $\CC(\0,\v,\Delta)$ is $\CC(\Delta')=\CC(\Delta)$
itself, for which we have $\wvphi\,\CC(\Delta)=\CC(\0,\0,\Delta)$.
\qed
\end{proof}

%The lemma indicates that the $\Psi$ is of $(-1)$-type case is the most promising
%place to start looking for whatever it is we are looking for.

%\begin{lemma}\label{section:extended:result300}
%Let $W(\Psi)$ be an irreducible Weyl group of $(-1)$-type with Coxeter
%number $h$. If $\xi$ is any Coxeter element of $W(\Psi)$ then $\xi^{h/2}=w_{\Psi}^{}$.
%\end{lemma}

If we have $\vphi$ instead, then the pairs $(\Psi,s_i)$ must all be
specially admissible, so that if $\vphi\CC(\Theta)=\CC(\v,\Delta)$ for $\Delta$
maximal, then again $\CC(\Theta)=\CC(\Delta)$, for which $\vphi\,\CC(\Delta)=
\CC(\0,\Delta)$.
Here is the main result of the section:

\begin{theorem}\label{section:extended:result400}
Let $W(\Gamma)$ be the Coxeter group $(\dag)$, where
$W(\Psi)$ 
%an irreducible Weyl group of $(-1)$-type with Coxeter
%number $h=2^pq$, for $q$ odd. 
is an irreducible Weyl group of rank $n$, with exponents $m_i$,
the dimension $d_\Psi>1$,
Coxeter number $h=2^pq$ where $p>0$ and $q$ odd,
the first $\ell$ pairs $(\Psi,s_i)$ admissible and the remaining ones
specially admissible.
Then there is a Coxeter element $\xi\in W(\Psi)$,
an $\x\in\Z/2^\ell$ and a $\u\in\prod\Lambda_i/2$, such that
$\zeta=(\x,\u,\xi^q)$ generates a subgroup,
$$
\langle\zeta\rangle\cong
\Z/(2^p)\subset
W(\Gamma)/\ker\wvphi,
%\cong\Z/2^\ell\times(\prod\Lambda_i/2\rtimes W(\Psi)),
$$ 
with $\wvphi^{-1}\langle\zeta\rangle\subset W(\Gamma)$ torsion free
of index $2^{mn+\ell-p}\prod_{i=1}^n(m_i+1)$.
\end{theorem}

The evenness of the Coxeter number rules out the type $A$ of even rank
and the dimension $d_\Psi>1$ rules out the type $A$'s of odd rank. Everything
else is ruled in.

\begin{proof}
We have 
$\zeta^k=(k\x,\u+\xi^q(\u)+\xi^{2q}(\u)+\cdots+\xi^{(k-1)q}(\u),\xi^{kq})$.
Let $\xi=s_1\ldots s_n$, $\x=\0$ and $\u=(\ov{u},\ldots,\ov{u})$
where $\ov{u}\in L/2$ is given by Proposition \ref{section:lattices:result500}.
Then $\zeta^{2^{p-1}}$ has middle component $\aa(\u)$ for the 
$\aa$ of Proposition \ref{section:lattices:result500}
and last component
$\xi^{(2^{p-1})q}=\xi^{h/2}$. Since $\aa(\u)\in\ker(\xi^{h/2}+1)$
we have by the proof of Lemma \ref{section:extended:result100}, that
$\zeta^{2^{p-1}}$ is an involution, and so $\langle\zeta\rangle\cong\Z/(2^p)$.
Moreover, by Proposition \ref{section:preliminaries:torsion:result400},
this involution lies in $\CC(\0,\aa(\u),\Delta)$, where $\Delta\in\JJ_\Psi$
is of maximal rank, and by Lemma \ref{section:extended:result200}(iii),
$\wvphi^{-1}\langle\zeta\rangle$ is torsion free if
$\CC(\0,\aa(\u),\Delta)\not=\CC(\0,\0,\Delta)$. But one last application
of Proposition \ref{section:lattices:result500}, together with
Lemma \ref{section:extended:result100}(ii), gives
$\zeta^{2^{p-1}}\not\in\CC(\0,\0,\Delta)$, and we are done.
%and in particular, by Lemma \ref{section:extended:result300},
%$\zeta^{2^p}=(\0,(\xi^{h/2}+1)\aa(\v),1)=
%(\0,(w_\Psi+1)\aa(\v),1)$, for
%$$
%\aa=1+\xi^q+\xi^{2q}+\cdots+\xi^{(2^{(p-1)}-1)q}\in\endo(A).
%$$
%Thus $\zeta$ has order $2^p$, as $w_\Psi+1=0$ in $\endo(A)$.
%The group $\langle\zeta\rangle$ has the single involution
%$\zeta^{2^{(p-1)}}=(\0,\aa(\v),\xi^{h/2})=(\0,\aa(\v),w_\Psi)$, so that
%$\zeta^{2^{(p-1)}}
%\in\CC(\0,\aa(\v),\Psi)$. By Lemmas
%\ref{section:extended:result100}-\ref{section:extended:result200},
%$\wvphi^{-1}\langle\zeta\rangle$ is 
%torsion free if and only if $\CC(\0,\aa(\v),\Psi)$ $\not=\CC(\0,\0,\Psi)$,
%which happens precisely when $\aa(\v)$ is not in the same $W(\Psi)$-orbit
%as $\0$, ie: when $\aa(\v)\not=\0$.
%
%Thus we require only that $\aa$ is not the zero map in $\endo(A)$, and then
%we can take any $\v\not\in\ker\aa$. This leaves us to 
%go through the $\Psi$ of $(-1)$-type on a case by case
%basis, with the $2^p$ given by,
%$$
%\begin{tabular}{cccccccc}
%\hline
%$\Psi$&$A_1$&$B_n$&$D_n\,(n\text{ even})$&$G_2$&$F_4$&$E_7$&$E_8$\\
%$2$-part of $h$&$2$&$2\cdot 2$-part of $n$&$2$&$2$&$4$&$2$&$2$\\
%\hline
%\end{tabular}
%$$
%for each $\Psi$. When $p=1$ we have $\aa=1$ and so $\zeta=(\x,\v,\xi^q)$
%generates a subgroup $\cong\Z/2$ that does the trick for any 
%$\xi$, $\x$ and $\v$ whatsoever. %leaving the $B_n$ and $F_4$ to do.
\qed
\end{proof}

%Observe that in all of the above the $\x\in\Z/2^\ell$ plays virtually no role,
%so that if the pairs $(\Psi,s_i)$ are all specially admissible, then
%everything above goes straight through with all references to $\x$ and
%$\Z/2^\ell$ removed. Thus we have $\TT$ the pairs now of $(\v,\Delta)$,
%$\zeta$ the element $(v,\xi^q)\in\prod\Lambda_i/2\rtimes  W(\Psi)$ and
%$\vphi^{-1}\langle\zeta\rangle$ torsion free of index $2^{m\,\rk\Psi-p}|W(\Psi)|$
%in $W(\Gamma)$.

In the case
where all the $(\Psi,s_i)$ are specially admissible,
the Theorem remains unchanged if $\wvphi$ is replaced by $\vphi$ 
(and the references to $\x$ and $\Z/2^\ell$ removed of course) 
Finally, we wrap-up with an omnibus result,

\begin{theorem}\label{section:extended:result500}
Let $W(\Gamma)$ be the Coxeter group with symbol $(\dag)$, where
\begin{description}
\item[(i).] $\Psi=A_n$ with $n$ odd. Then the element $\zeta=(\0,\u,\xi^{h/2})$ 
for some $\u$ and $\xi=s_1\ldots s_n$, generates
a subgroup $\cong\Z/2$ with $\wvphi^{-1}\langle\zeta\rangle$ torsion free.
\item[(ii).] $\Psi=E_6$. Then if $\xi=s_2\ldots s_n$ is a Coxeter element of the
visible $W(D_5)=\langle s_2,\ldots,s_6\rangle$, the element
$\zeta=(\0,\u,\xi)$ for some $\u$ generates
a subgroup $\cong\Z/(2^3)$ with $\wvphi^{-1}\langle\zeta\rangle$ torsion free.
\end{description}
\end{theorem}

Notice that while $E_6$ does in fact fall under the auspices of Theorem
\ref{section:extended:result400}, the torsion free
subgroup given there has twice the index as that
given above. The proof of Theorem \ref{section:extended:result500} is
entirely analogous to Theorem \ref{section:extended:result400}, and 
uses the relevant calculations from the end of \S\ref{section:lattices}.

\section{Geometric manifolds and their symmetries}
\label{section:manifolds}

Let $X$ be one of the three 
simply connected Riemannian $n$-manifolds of constant sectional 
curvature, ie: $X=$ the $n$-sphere $S^n$, the Euclidean space
$\E^n$, or the % unit\footnote{\emph{unit\/} in this context means the sphere
% of squared radius $-1$ in the $n+1$-dimensional Minkowski space-time.}
$n$-dimensional hyperbolic space $\H^n$. Let $\isom X$ be the group of isometries
of $X$ and $G\subset\isom X$ a discrete subgroup. Then $G$ acts 
properly discontinuously on $X$, and in the presence of
curvature (ie: when $X\not=\E^n$), the volume of a fundamental region
is an invariant for $G$, the \emph{covolume\/} $\covol(G)$.

If $G_1\subset G_2\subset\isom X$ are discrete, then $M=X/G_1$ is a 
complete Riemannian $n$-manifold with constant sectional curvature
that of $X$ if and only
if $G_1$ acts freely. In this case
$G_1\cong\pi_1(M)$ and $M$ has volume $\vol(M)=[G_2:G_1]\,\covol(G_2)$.
Moreover, we have a free action whenever $G$ is torsion free (the converse
is not true when $\isom X$ is compact, ie: when $X=S^n$, where 
discrete subgroups are necessarily finite, but there still exist free actions 
of non-trivial groups).

An automorphism of $M=X/G$ is an isometry $M\rightarrow M$, and the group
$\aut(M)$ of all such is isomorphic to $N(G)/G$, where $N(G)$ is the normalizer
in $\isom X$ of $G$. A subgroup $H\subset\aut(M)$ acts freely on $M$ precisely when
$G\lhd \eta^{-1}H$ is torsion free, for $\eta:N(G)\rightarrow N(G)/G$ the
quotient homomorphism. Letting $\widehat{M}=X/\eta^{-1}H$, then gives a Galois 
covering of manifolds $M\rightarrow\widehat{M}$ with group of deck transformations $H$.

Suppose now that $G$ is (isomorphic to) a Coxeter group $W(\Gamma)$. Then $W(\Gamma)$
is finite when $X=S^n$, 
an affine Weyl group (as in eg: \cite{Humphreys90}*{Chapter 4}) when $X=\E^n$, and
hyperbolic in the sense of \S\ref{subsection:coxeter_groups} when $X=\H^n$.
%As it is only in the hyperbolic case that we get both an exact correspondence
%between torsion free goup asnd groups acting freely, as well as a notion
%of covolume, we will restrict to this case from now on. 

Looking more closely at $X=\H^n$, there is no general method for computing 
the covolumes of hyperbolic Coxeter groups.
If the fundamental region is a polytope with a particularly simple structure, such 
as a simplex, then the covolumes have been computed
\cite{Johnson99}. If $n$ is even, one can instead appeal to homological
considerations, for the Gauss-Bonnet(-Hirzebruch)
formula \cites{Gromov82,Spivak79}, gives
\begin{equation}\label{euler_char}
\covol W(\Gamma)=\kappa_n\,\chi\,W(\Gamma),\text{ where }
\kappa_n=2^n (n!)^{-1} (-\pi)^{n/2} (n/2)!
\end{equation}
and $\chi\,W(\Gamma)=
\sum_i(-1)^i\text{rank}_{\Z}H_i\,W(\Gamma)$ is the (group) Euler characteristic of $W(\Gamma)$.
By \cites{Chiswell92}, 
the computation of the Euler characteristic of a Coxeter group is 
a simple combinatorial exercise, requiring no more information than the symbol
$\Gamma$ and the orders of the finite visible subgroups.

The upshot of all this is that Theorems 
\ref{section:homomorphisms:result300}-\ref{section:extended:result500}
have direct geometric analogues. Here, for instance, is
the one for Theorem \ref{section:extended:result400}:

\begin{theorem}\label{manifolds:result100}
Let $W(\Gamma)$ be a hyperbolic Coxeter group with the symbol
$(\dag)$ of \S$\ref{section:homomorphisms}$, where $W(\Psi)$ is an
irreducible Weyl group of rank $n$, with exponents $m_i$, dimension
$d_\Psi>1$, Coxeter number $h=2^p\,q$ for $p>0$ and $q$ odd, 
the first $\ell$ pairs $(\Psi,s_i)$ admissible and the remaining
ones specially admissible. Then
there is a Galois covering,
$$
M\rightarrow\widehat{M},
$$
of $N$-dimensional hyperbolic manifolds, with deck transformation group
$\Z/(2^p)$, and 
$$
\vol(\widehat{M})=2^{mn+\ell-p}\prod(m_i+1)\covol W(\Gamma),
$$
where $N\geq n$ is the largest rank of a finite visible subgroup.
\end{theorem}

\parshape=7 0pt\hsize 0pt\hsize 0pt\hsize 0pt\hsize 
0pt.55\hsize 0pt.55\hsize 0pt.55\hsize
We finish with an explicit example.
A $(n+1)$-dimensional \emph{Lorentzian\/} lattice
is an $(n+1)$-dimensional 
free $\Z$-module equipped with a $\Z$-valued bilinear form of signature $(n,1)$. 
For each $n$, there is a unique such, denoted
$I_{n,1}$, that is odd and self-dual (see \cite[Theorem V.6]{Serre73},
or \cite{Milnor73, Neumaier83}).
By \cite{Borel62}, the
group $\text{PO}_{n,1}\Z$ of (projectivised) automorphisms of $I_{n,1}$ acts 
properly discontinuously by isometries, and with finite covolume, 
on the hyperbolic space $\H^n$ obtained by 
projectivising the negative norm vectors
in $I_{n,1}\otimes\R$. %(to get a faithful 
%action one normally
%passes to the centerless version $\po_{n,1}\Z$). 
\vadjust{\hfill\smash{\lower 7mm
\llap{
\begin{pspicture}(0,0)(6.5,2)
%\showgrid
\rput(3.5,1){\BoxedEPSF{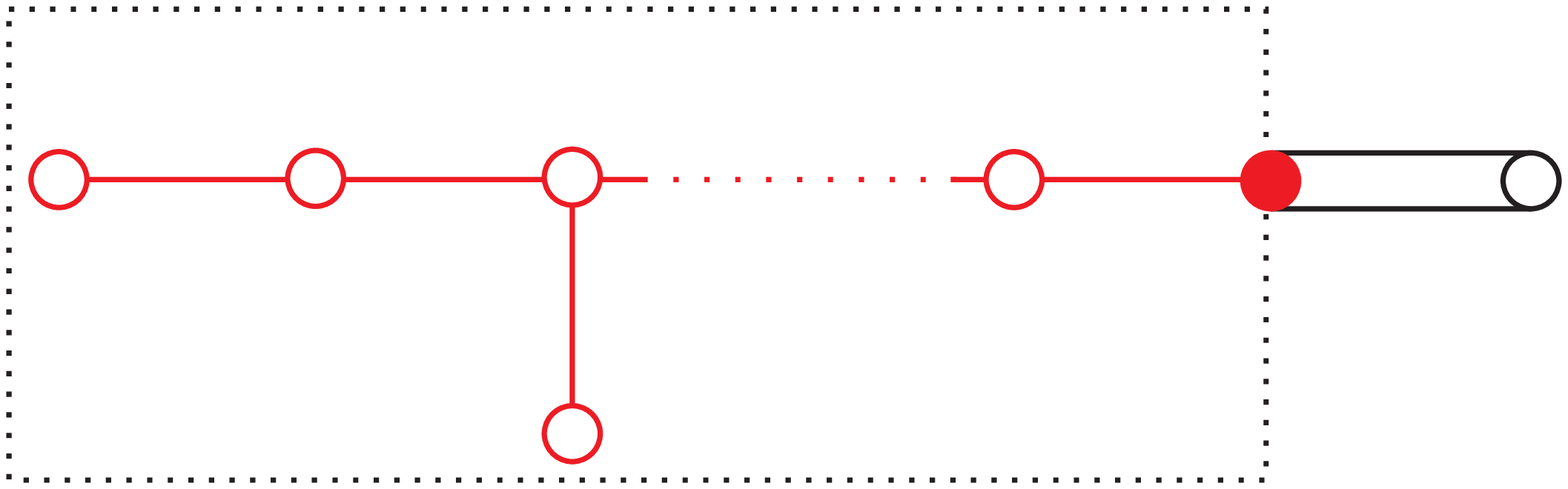 scaled 300}}
\rput(.6,.3){${\red \Psi}$}
\end{pspicture}
}}}\ignorespaces

\parshape=3 0pt.55\hsize 0pt.55\hsize 0pt\hsize
Vinberg and Kaplinskaja showed \cite{Vinberg78,Vinberg72} that
the subgroup $W(\Gamma)$ of $\po_{n,1}\Z$ generated by reflections 
in positive norm vectors
has finite index if 
and only if $n\leq 19$, thus yielding a family of 
finite covolume $n$-dimensional hyperbolic Coxeter groups
for $2\leq n\leq 19$. 
Indeed, Conway and Sloane \cite{Conway99}*{Chapter 28}
have shown
that for $n\leq 19$ the quotient 
of $\po_{n,1}\Z$ by $W(\Gamma)$ is a subgroup of the automorphism 
group of the Leech lattice.
%Borcherds \cite{Borcherds87} showed that the (non self-dual) even 
%sublattice of $I_{21,1}$
%also acts cofinitely, yielding the highest dimensional example 
%known of a Coxeter group
%acting cofinitely on hyperbolic space.

The group $W(\Gamma)$ has symbol 
$
\begin{pspicture}(0,0)(2.5,.4)
%\showgrid
\rput(1.25,.1){\BoxedEPSF{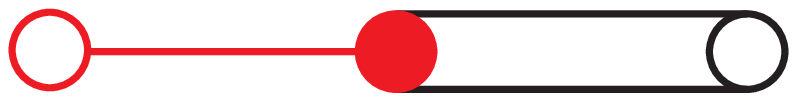 scaled 300}}
\rput(.7,.3){${\red \infty}$}
\end{pspicture}
$ 
and 
$\BoxedEPSF{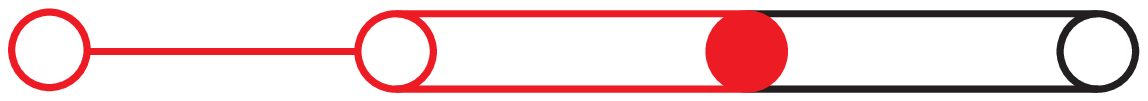 scaled 300}$ when $n=2,3$, and
when $4\leq n\leq 9$, the symbol 
shown above,
with $n+1$ nodes and the fundamental region a non-compact, finite volume
$n$-simplex in $\H^n$ (for $9<n\leq 19$, the region has a more complicated 
structure). The symbol $\Psi$ is,
$$
\begin{tabular}{ccccccc}
\hline
$n$&$4$&$5$&$6$&$7$&$8$&$9$\\
$\Psi$&$A_4$&$D_5$&$E_6$&$E_7$&$E_8$&$\widetilde{E}_8$\\
\hline
\end{tabular}
$$
Thus, we have $W(\Gamma)$ hyperbolic with the symbol $(\dag)$
of \S\ref{section:homomorphisms}, and $W(\Psi)$ an irreducible Weyl group
when $3\leq n\leq 8$ (but for $n=2$ and $9$, $W(\Psi)$ is the \emph{affine\/}
Weyl group $\widetilde{A}_1$ and $\widetilde{E}_8$). By Theorems
\ref{section:lattices:result200}-\ref{section:lattices:result250}
the pair $(\Psi,s)$, with $s$ the red node, is admissible for
$n=4,6$ and $8$. 

Thus, let $\Psi=A_4, E_6$ or $E_8$, with $m=\ell=1$, $h=2^p\,q$ the Coxeter
number of these three (with $p=0,2,1$ respectively) and $m_i$ the exponents 
(which are the positive integers relatively prime to the Coxeter number for 
$A_4$ and $E_8$,
and these, together with $4$ and $8$ for $E_6$).
The covolume of $W(\Gamma)$ can be computed
using the Euler characteristic (\ref{euler_char}), but in this particular case
Siegel \cite{Siegel36} gave the volume in terms of a limit
which was evaluated in \cite{Ratcliffe97} to give the following remarkable 
formula,
$$
\covol W(\Gamma)=\frac{(2^{\frac{n}{2}}\pm 1)\pi^{\frac{n}{2}}}{n!}
\prod_{k=1}^{\frac{n}{2}} |B_{2k}|,
$$
where $B_{2k}$ is the $2k$-th Bernoulli number, and we have the minus sign
for $n=4,6$ and the plus sign for $n=8$. 

A direct application of Theorem \ref{manifolds:result100} and the geometric 
analogue of Theorem \ref{section:homomorphisms:result500}
now gives a Galois covering,
$M\rightarrow\widehat{M}$,
of $4,6$ and $8$-dimensional hyperbolic manifolds, with deck transformation group
$\Z/(2^p)$, and 
$$
\vol(\widehat{M})=2^{n-p+1}
\frac{(2^{\frac{n}{2}}\pm 1)\pi^{\frac{n}{2}}}{n!}
\prod_{i=1}^n(m_i+1)
\prod_{k=1}^{\frac{n}{2}} |B_{2k}|.
$$
In $4$-dimensions this translates to $\frac{8}{3}\pi^2$, or a manifold
of Euler characteristic $\chi=2$. In $6$-dimensions, we do even better
if we use instead the geometric analogue of Theorem \ref{section:extended:result500}(ii),
obtaining the volume $\frac{16}{15}\pi^3$, or a manifold of Euler
characteristic $\chi=-2$.

%%%%%%%%%%%%%%%%%%%%%%%%%%%%%%%%%%%%%%%%%%%%%%%%%%%%%%%%%%%%%%%%%%%%%%%%%

%\bibliography{brent}

\section*{References}

\begin{biblist}

%\bib{Borcherds87}{article}{
%   author={Borcherds, Richard},
%   title={Automorphism groups of Lorentzian lattices},
%   journal={J. Algebra},
%   volume={111},
%   date={1987},
%   number={1},
%   pages={133--153},
%   issn={0021-8693},
%%   review={\MR{913200 (89b:20018)}},
%}

\bib{Borel62}{article}{
   author={Borel, Armand},
   author={Harish-Chandra},
   title={Arithmetic subgroups of algebraic groups},
   journal={Ann. of Math. (2)},
   volume={75},
   date={1962},
   pages={485--535},
   issn={0003-486X},
   %review={MR{0147566 (26 #5081)}},
}

\bib{Bourbaki02}{book}{
    author={Bourbaki, Nicolas},
     title={Lie groups and Lie algebras. Chapters 4--6},
    series={Elements of Mathematics (Berlin)},
      note={Translated from the 1968 French original by Andrew Pressley},
 publisher={Springer-Verlag},
     place={Berlin},
      date={2002},
     pages={xii+300},
      isbn={3-540-42650-7},
%    %review={MR1890629 (2003a:17001)},
}

\bib{Brink93}{article}{
   author={Brink, Brigitte},
   author={Howlett, Robert B.},
   title={A finiteness property and an automatic structure for Coxeter
   groups},
   journal={Math. Ann.},
   volume={296},
   date={1993},
   number={1},
   pages={179--190},
   issn={0025-5831},
%   %review={MR{1213378 (94d:20045)}},
}

\bib{Carter72}{article}{
   author={Carter, R. W.},
   title={Conjugacy classes in the Weyl group},
   journal={Compositio Math.},
   volume={25},
   date={1972},
   pages={1--59},
   issn={0010-437X},
%   %review={MR{0318337 (47 \#6884)}},
}

\bib{Chiswell92}{article}{
   author={Chiswell, I. M.},
   title={The Euler characteristic of graph products and of Coxeter groups},
   conference={
      title={Discrete groups and geometry},
      address={Birmingham},
      date={1991},
   },
   book={
      series={London Math. Soc. Lecture Note Ser.},
      volume={173},
      publisher={Cambridge Univ. Press},
      place={Cambridge},
   },
   date={1992},
   pages={36--46},
%   review={\MR{1196914 (94a:05090)}},
}

%\bib{Chiswell76}{article}{
%   author={Chiswell, Ian M.},
%   title={Euler characteristics of groups},
%   journal={Math. Z.},
%   volume={147},
%   date={1976},
%   number={1},
%   pages={1--11},
%   issn={0025-5874},
%%   review={\MR{0396785 (53 \#645)}},
%}

\bib{Conway99}{book}{
   author={Conway, J. H.},
   author={Sloane, N. J. A.},
   title={Sphere packings, lattices and groups},
   series={Grundlehren der Mathematischen Wissenschaften [Fundamental
   Principles of Mathematical Sciences]},
   volume={290},
   edition={3},
   note={With additional contributions by E. Bannai, R. E. Borcherds, J.
   Leech, S. P. Norton, A. M. Odlyzko, R. A. Parker, L. Queen and B. B.
   Venkov},
   publisher={Springer-Verlag},
   place={New York},
   date={1999},
   pages={lxxiv+703},
   isbn={0-387-98585-9},
 %  review={\MR{1662447 (2000b:11077)}},
}

%\bib{Conway82}{article}{
%   author={Conway, J. H.},
%   author={Sloane, N. J. A.},
%   title={Leech roots and Vinberg groups},
%   journal={Proc. Roy. Soc. London Ser. A},
%   volume={384},
%   date={1982},
%   number={1787},
%   pages={233--258},
%   issn={0080-4630},
% %  review={\MR{684311 (84b:10048)}},
%}

\bib{Coxeter88}{article}{
   author={Coxeter, H. S. M.},
   title={Regular and semi-regular polytopes. III},
   journal={Math. Z.},
   volume={200},
   date={1988},
   number={1},
   pages={3--45},
   issn={0025-5874},
 %  review={\MR{972395 (89m:52005)}},
}

%\bib{Coxeter35}{article}{
%   author={Coxeter, H. S. M.},
%   title={The complete enumeration of finite groups of the form $R_i^2=
%          (R_iR_k)^{k_{ij}}=1$},
%   journal={J. London Math. Soc.},
%   volume={10},
%   date={1935},
%   pages={21--25},
%}
%
%\bib{Coxeter34}{article}{
%   author={Coxeter, H. S. M.},
%   title={Discrete groups generated by reflections},
%   journal={Ann. of Math. (2)},
%   volume={35},
%   date={1934},
%   number={3},
%   pages={588--621},
%   issn={0003-486X},
% %  review={\MR{1503182}},
%}
%
%\bib{Davis85}{article}{
%   author={Davis, Michael W.},
%   title={A hyperbolic $4$-manifold},
%   journal={Proc. Amer. Math. Soc.},
%   volume={93},
%   date={1985},
%   number={2},
%   pages={325--328},
%   issn={0002-9939},
%%   %review={MR{770546 (86h:57016)}},
%}

\bib{Deodhar82}{article}{
   author={Deodhar, Vinay V.},
   title={On the root system of a Coxeter group},
   journal={Comm. Algebra},
   volume={10},
   date={1982},
   number={6},
   pages={611--630},
   issn={0092-7872},
 %  review={\MR{647210 (83j:20052a)}},
}

\bib{Everitt04a}{article}{
    author={Everitt, Brent},
     title={Coxeter groups and hyperbolic manifolds},
   journal={Math. Ann.},
    volume={330},
      date={2004},
    number={1},
     pages={127\ndash 150},
      issn={0025-5831},
%    %review={$\,$MR2091682 (2005m:20090)},
%    eprint={arXiv:math.GT/0205157}
}

\bib{Everitt04b}{article}{
    author={Everitt, Brent},
     title={3-manifolds from Platonic solids},
   journal={Topology Appl.},
    volume={138},
      date={2004},
    number={1-3},
     pages={253\ndash 263},
      issn={0166-8641},
%    %review={$\,$MR2035484 (2004m:57031)},
%    eprint={arXiv:math.GT/0104182}
}

\bib{Everitt05}{article}{
    author={Everitt, Brent},
    author={Ratcliffe, John},
    author={Tschantz, Steven},
     title={The smallest hyperbolic 6-manifolds},
   journal={Electron. Res. Announc. Amer. Math. Soc.},
    volume={11},
      date={2005},
     pages={40\ndash 46},
      issn={1079-6762},
%    %review={$\,$MR2150943 (2006a:57019)},
%    eprint={arXiv:math.GT/0410474}
}

\bib{Everitt05}{article}{
    author={Everitt, Brent},
    author={Ratcliffe, John},
    author={Tschantz, Steven},
     title={Arithmetic hyperbolic 6-manifolds of smallest volume},
   note={(in preparation)},
}

\bib{Everitt00}{article}{
   author={Everitt, Brent},
   author={Maclachlan, Colin},
   title={Constructing hyperbolic manifolds},
   conference={
      title={Computational and geometric aspects of modern algebra
      (Edinburgh, 1998)},
   },
   book={
      series={London Math. Soc. Lecture Note Ser.},
      volume={275},
      publisher={Cambridge Univ. Press},
      place={Cambridge},
   },
   date={2000},
   pages={78--86},
 %  review={\MR{1776768 (2001i:57022)}},
}

\bib{Gabai07}{article}{
    author={Gabai, David},
    author={Meyerhoff, Robert},
    author={Milley, Peter},
     title={Minimum volume cusped hyperbolic $3$-manifolds},
    eprint={arXiv:0705.4325 [math.GT]}, 
}


\bib{Gromov82}{article}{
   author={Gromov, Michael},
   title={Volume and bounded cohomology},
   journal={Inst. Hautes \'Etudes Sci. Publ. Math.},
   number={56},
   date={1982},
   pages={5--99 (1983)},
   issn={0073-8301},
%   review={\MR{686042 (84h:53053)}},
}

%\bib{Haiman92}{article}{
%   author={Haiman, Mark D.},
%   title={Dual equivalence with applications, including a conjecture of
%   Proctor},
%   journal={Discrete Math.},
%   volume={99},
%   date={1992},
%   number={1-3},
%   pages={79--113},
%   issn={0012-365X},
%   %review={MR{1158783 (93h:05173)}},
%}

\bib{Howlett80}{article}{
   author={Howlett, Robert B.},
   title={Normalizers of parabolic subgroups of reflection groups},
   journal={J. London Math. Soc. (2)},
   volume={21},
   date={1980},
   number={1},
   pages={62--80},
   issn={0024-6107},
 %  review={\MR{576184 (81g:20094)}},
}

\bib{Humphreys90}{book}{
    author={Humphreys, James E.},
     title={Reflection groups and Coxeter groups},
    series={Cambridge Studies in Advanced Mathematics},
    volume={29},
 publisher={Cambridge University Press},
     place={Cambridge},
      date={1990},
     pages={xii+204},
      isbn={0-521-37510-X},
%    %review={MR1066460 (92h:20002)},
}

\bib{Johnson99}{article}{
   author={Johnson, N. W.},
   author={Ratcliffe, J. G.},
   author={Kellerhals, R.},
   author={Tschantz, S. T.},
   title={The size of a hyperbolic Coxeter simplex},
   journal={Transform. Groups},
   volume={4},
   date={1999},
   number={4},
   pages={329--353},
   issn={1083-4362},
%   review={\MR{1726696 (2000j:20070)}},
}

\bib{Kane01}{book}{
    author={Kane, Richard},
     title={Reflection groups and invariant theory},
    series={CMS Books in Mathematics/Ouvrages de Math\'ematiques de la SMC,
            5},
 publisher={Springer-Verlag},
     place={New York},
      date={2001},
     pages={x+379},
      isbn={0-387-98979-X},
%    %review={MR1838580 (2002c:20061)},
}

%\bib{Mihalik}{article}{
%   author={Mihalik, Michael},
%   author={Ratcliffe, John},
%   author={Tschantz, Steven},
%   title={Quotient isomorphism invariants of a finitely generated Coxeter group},
%   eprint={arXiv:math.GR/0607428}
%}

\bib{Milnor73}{book}{
   author={Milnor, John},
   author={Husemoller, Dale},
   title={Symmetric bilinear forms},
   note={Ergebnisse der Mathematik und ihrer Grenzgebiete, Band 73},
   publisher={Springer-Verlag},
   place={New York},
   date={1973},
   pages={viii+147},
   %review={MR{0506372 (58 #22129)}},
}

\bib{Neumaier83}{article}{
   author={Neumaier, A.},
   author={Seidel, J. J.},
   title={Discrete hyperbolic geometry},
   journal={Combinatorica},
   volume={3},
   date={1983},
   number={2},
   pages={219--237},
   issn={0209-9683},
   %review={MR{726460 (85h:51025)}},
}

%\bib{Ratcliffe04}{article}{
%   author={Ratcliffe, John G.},
%   author={Tschantz, Steven T.},
%   title={Integral congruence two hyperbolic 5-manifolds},
%   journal={Geom. Dedicata},
%   volume={107},
%   date={2004},
%   pages={187--209},
%   issn={0046-5755},
%%   %review={MR{2110762 (2005m:57026)}},
%}

\bib{Ratcliffe00}{article}{
   author={Ratcliffe, John G.},
   author={Tschantz, Steven T.},
   title={The volume spectrum of hyperbolic 4-manifolds},
   journal={Experiment. Math.},
   volume={9},
   date={2000},
   number={1},
   pages={101--125},
   issn={1058-6458},
%   %review={MR{1758804 (2001b:57048)}},
}

\bib{Ratcliffe97}{article}{
   author={Ratcliffe, John G.},
   author={Tschantz, Steven T.},
   title={Volumes of integral congruence hyperbolic manifolds},
   journal={J. Reine Angew. Math.},
   volume={488},
   date={1997},
   pages={55--78},
   issn={0075-4102},
%   review={\MR{1465367 (99b:11076)}},
}

\bib{Richardson82}{article}{
   author={Richardson, R. W.},
   title={Conjugacy classes of involutions in Coxeter groups},
   journal={Bull. Austral. Math. Soc.},
   volume={26},
   date={1982},
   number={1},
   pages={1--15},
   issn={0004-9727},
%   %review={MR{679916 (84c:20056)}},
}

% \bib{Selberg60}{article}{
%    author={Selberg, Atle},
%    title={On discontinuous groups in higher-dimensional symmetric spaces},
%    conference={
%       title={Contributions to function theory (internat. Colloq. Function
%       Theory, Bombay, 1960)},
%    },
%    book={
%       publisher={Tata Institute of Fundamental Research},
%       place={Bombay},
%    },
%    date={1960},
%    pages={147--164},
%    review={\MR{0130324 (24 \#A188)}},
% }
	

\bib{Serre73}{book}{
   author={Serre, J.-P.},
   title={A course in arithmetic},
   note={Translated from the French;
   Graduate Texts in Mathematics, No. 7},
   publisher={Springer-Verlag},
   place={New York},
   date={1973},
   pages={viii+115},
   %review={MR{0344216 (49 #8956)}},
}

\bib{Siegel45}{article}{
   author={Siegel, Carl Ludwig},
   title={Some remarks on discontinuous groups},
   journal={Ann. of Math. (2)},
   volume={46},
   date={1945},
   pages={708--718},
   issn={0003-486X},
 %  review={\MR{0014088 (7,239c)}},
}

\bib{Siegel36}{article}{
   author={Siegel, Carl Ludwig},
   title={\"Uber die analytische Theorie der quadratischen Formen. II},
   language={German},
   journal={Ann. of Math. (2)},
   volume={37},
   date={1936},
   number={1},
   pages={230--263},
   issn={0003-486X},
 %  review={\MR{1503276}},
}

\bib{Spivak79}{book}{
   author={Spivak, Michael},
   title={A comprehensive introduction to differential geometry. Vol. V},
   edition={2},
   publisher={Publish or Perish Inc.},
   place={Wilmington, Del.},
   date={1979},
   pages={viii+661},
   isbn={0-914098-83-7},
   review={\MR{532834 (82g:53003e)}},
}

\bib{Springer78}{article}{
   author={Springer, T. A.},
   title={A construction of representations of Weyl groups},
   journal={Invent. Math.},
   volume={44},
   date={1978},
   number={3},
   pages={279--293},
   issn={0020-9910},
 %  review={\MR{0491988 (58 \#11154)}},
}

\bib{Vinberg85}{article}{
   author={Vinberg, {\`E}. B.},
   title={Hyperbolic groups of reflections},
   language={Russian},
   journal={Russian Math. Surveys},
   volume={40},
   date={1985},
   number={1(241)},
   pages={31--75},
   issn={0042-1316},
%   %review={\MR{783604 (86m:53059)}},
}

\bib{Vinberg72}{article}{
   author={Vinberg, {\`E}. B.},
   title={The groups of units of certain quadratic forms},
   language={Russian},
   journal={Mat. Sb. (N.S.)},
   volume={87(129)},
   date={1972},
   pages={18--36},
   %review={MR{0295193 (45 #4261)}},
}

\bib{Vinberg78}{article}{
   author={Vinberg, {\`E}. B.},
   author={Kaplinskaja, I. M.},
   title={The groups $O\sb{18,1}(Z)$ and $O\sb{19,1}(Z)$},
   language={Russian},
   journal={Dokl. Akad. Nauk SSSR},
   volume={238},
   date={1978},
   number={6},
   pages={1273--1275},
   issn={0002-3264},
   %review={MR{0476640 (57 #16199)}},
}

\end{biblist}
%\bibliographystyle{plain}

%\end{bibdiv}

\end{document}